\documentclass[preprint, 9pt]{elsarticle}

\usepackage{jasr}
\usepackage[utf8]{inputenc}
\usepackage{amssymb}
\usepackage{bm}
\usepackage{amsmath}
\usepackage{color,soul}
\usepackage{todonotes}
\usepackage{colortbl}
\definecolor{myblue}{rgb}{0, 0.23, 0.64}
\definecolor{myred}{HTML}{8e0000}

\usepackage{breakurl}
\usepackage[breaklinks]{hyperref}
\hypersetup{
    colorlinks=true,
    linkcolor=myblue,
    citecolor=myblue,
    urlcolor=myred,
}

\makeatletter
\AtBeginDocument{%
    \def\@citecolor{myblue}
    \def\@linkcolor{myblue}
    \def\@urlcolor{myred}
}
\def\ps@pprintTitle{%
  \let\@oddhead\@empty
  \let\@evenhead\@empty
  \let\@oddfoot\@empty
  \let\@evenfoot\@oddfoot
}
\makeatother

\usepackage{graphicx}
\usepackage[version=4]{mhchem}
\usepackage[group-separator={,}]{siunitx}
\usepackage{longtable,tabularx}
\usepackage{textcomp}
\usepackage{gensymb}
\usepackage{mathtools}
\usepackage{float}
\usepackage[ruled, vlined, linesnumbered]{algorithm2e}
\SetKwComment{Comment}{$\triangleright$\ }{}
\usepackage{subcaption}
\usepackage{mwe}
\usepackage{enumitem}
\usepackage{multirow}
\usepackage{booktabs}
\usepackage{lscape}
\usepackage{mdframed}
\usepackage{adjustbox}
\usepackage{dsfont}
\usepackage{xspace}
\usepackage{tcolorbox}
\usepackage[flushleft]{threeparttable}
\usepackage{xurl}

\journal{Advances in Space Research}

\begin{document}

\begin{frontmatter}

\title{Optimal Placement and Coordinated Scheduling of Distributed Space-Based Lasers for Orbital Debris Remediation}

\author[WVU]{David O. \snm{Williams Rogers}\fnref{fn1}}
\ead{davidowen.williamsrogers@mail.wvu.edu}
\author[WVU]{Matthew C. \snm{Fox}}
\ead{mcf0014@mail.wvu.edu}
\author[NASA]{Paul R. \snm{Stysley}}
\ead{paul.stysley@nasa.gov}
\author[WVU]{Hang Woon \snm{Lee}\corref{cor}}

\affiliation[WVU]{organization={Department of Mechanical, Materials and Aerospace Engineering, West Virginia University},
            addressline={1306 Evansdale Drive}, 
            city={Morgantown},
            state={WV},
            postcode={26506}, 
            country={USA}}

\affiliation[NASA]{organization={Lasers and Electro-optics Branch, NASA Goddard Space Flight Center},
            addressline={8800 Greenbelt Road}, 
            city={Greenbelt},
            state={MD},
            postcode={20771},
            country={USA}}

\cortext[cor]{Corresponding author: Email: \href{mailtohangwoon.lee@mail.wvu.edu}{hangwoon.lee@mail.wvu.edu}.}

\begin{abstract}
The significant expansion of the orbital debris population poses a serious threat to the safety and sustainability of space operations. This paper investigates orbital debris remediation through a constellation of collaborative space-based lasers, leveraging the principle of momentum transfer onto debris via laser ablation. A novel delta-v vector analysis framework quantifies the cumulative effects of multiple concurrent laser-to-debris (L2D) engagements by utilizing the vector composition of the imparted delta-v vectors. The paper formulates the Concurrent Location-Scheduling Optimization Problem (CLSP) to optimize the placement of laser platforms and the scheduling of L2D engagements, aiming to maximize debris remediation capacity. Given the computational intractability of the CLSP, a decomposition strategy is employed, yielding two sequential subproblems: (1) determining optimal laser platform locations via the Maximal Covering Location Problem, and (2) scheduling L2D engagements using a novel integer linear programming approach to maximize debris remediation capacity. Computational experiments evaluate the efficacy of the proposed framework across diverse mission scenarios, demonstrating critical constellation functions such as collaborative and controlled nudging, deorbiting, and just-in-time collision avoidance. A sensitivity analysis further explores the impact of varying the number and distribution of laser platforms on debris remediation capacity, offering insights into optimizing the performance of space-based laser constellations.
\end{abstract}

\begin{keyword}
\KWD Space-based lasers\sep Orbital debris remediation\sep Laser ablation\sep Constellation configuration optimization\sep Constellation scheduling\sep Integer linear programming
\end{keyword}
\end{frontmatter}

% \linenumbers

\section{Introduction} \label{sec:intro}
The number of resident space objects is rapidly increasing, largely due to the emergence of small satellites, owing to standardized manufacturing processes, advancements in technology, and lower costs for sharing space on launch vehicles. This increase in the number of satellites is correlated with the significant growth in orbital debris, posing a serious risk to both manned and unmanned missions as well as operational satellites that form critical infrastructure, yielding an escalation in conjunction alerts (\textit{e.g.,} satellites approaching within \SI{3}{km} \citep{russian}) and an increased probability of in-space collisions, which can create debris clouds of varying sizes.

Large orbital debris, objects with characteristic lengths greater than \SI{10}{cm}, is identified as the primary source of new debris \citep{mcknight2021} and consists mostly of defunct satellites, payloads, and rocket stages. Its presence in orbit increases the cost of operations, requiring operating satellites to perform collision avoidance maneuvers, thereby demanding additional fuel consumption, the loss of mission objectives due to maneuvering, and increased labor costs for planning these maneuvers. Additionally, large debris restricts access to space for new satellites or payloads, affecting activities that heavily rely on space infrastructure, such as telecommunications, financial exchanges, and climate monitoring~\citep{colvin2023cost}.

Conversely, small debris comprises objects with characteristic lengths between 1 and \SI{10}{cm}. The challenges of tracking small debris pieces, coupled with their significant population close to \num{500000} objects~\citep{debris_whitehouse}, present a continuous threat of collisions in space. Even with their relatively small mass, these fragments have the potential to inflict hypervelocity impacts, penetrate spacecraft shielding, and risk the mission's success. The origins of small debris include loose parts from operational or defunct satellites, payloads, or rocket stages, such as sodium-potassium droplets, solid rocket motor slag and dust, multi-layer insulation, ejecta, and paint flakes~\citep{MASTER8}. Notably, a significant part of the population has been generated by anti-satellite tests~\citep{antisatall}, accidental explosions of satellites and rocket bodies~\citep{antisatCHINA}, and debris-to-debris collisions~\citep{CLASH}.
 
Considering the substantial size of the debris population and the looming threat of triggering the Kessler syndrome~\citep{kessler1978collision}, strategies limited to the prevention of new debris formation, such as implementing regular conjunction assessments for active spacecraft~\citep{JOHNSON2010362} and enforcing restrictions on mission operations and end-of-life procedures~\citep{debris_whitehouse} are insufficient for addressing the debris problem comprehensively. In response to this, a range of promising and innovative debris remediation technologies has been proposed to reduce the debris population.

The literature on debris remediation introduces contact-based methods as strategies for eliminating orbital debris. Chaser spacecraft equipped with tethered nets are considered an effective debris remediation mechanism~\citep{sharf2017experiments} for addressing irregularly shaped debris.
Space balls and space winches are mechanisms designed to tackle the problem of debris with high angular momentum. These methods are capable of exerting a retarding torque on debris, aiding in its stabilization. However, they require an accurate activation at the time of the engagement~\citep{mcknight2013detumbling}.
For relatively stable (\textit{i.e.,} detumbled)  objects, debris remediation mechanisms such as grappling can be implemented for targeted and controlled removal. For instance, spacecraft equipped with robotic arms can accomplish this task~\citep{shan2016review} by leveraging a rendezvous maneuver to capture the target and subsequently relocating it into a disposal orbit or into an orbit where atmospheric effects induce reentry. However, this technique is highly dependent on the physical properties of the debris, such as its shape and surface texture. It also necessitates precise attitude adjustment to counteract the forces exerted by the robotic arm's movements~\citep{pulliam2011catcher}. Hooks and harpoons are penetrating mechanisms used for debris remediation; however, they require proximity operations, endangering the mission's success, and during penetration, they can trigger an explosion~\citep{mcknight2013detumbling}, generating additional debris pieces.

Unlike the contact-based debris remediation mechanisms described above, which are targeted at large debris, remediating small debris requires specific methods compatible with their size. Novel tethered plate systems, capable of remediating debris pieces without damaging the plate's structure, have been proposed~\citep{takeichi2021tethered}. Furthermore, perimeter-ring-truss systems can handle small debris while ensuring scalability and practicality in the design phase~\citep{foster2022practical}. However, the cost of remediation per debris piece is analyzed to be the highest in comparison with other methods~\citep{colvin2023cost}.

Ground-based lasers have emerged as a promising, cost-effective debris remediation solution to address the growing orbital debris population while circumventing the challenges and limitations of contact-based debris remediation methods. 
The change in velocity, $\Delta v$, required to reduce the target's orbit altitude can be imparted by a ground-based laser leveraging photon pressure or laser ablation mechanisms. Photon pressure-based lasers exert a small force on the target object sufficient to induce small orbit changes. \citep{mason2011orbital} present the use of two photon pressure ground-based laser systems to engage debris multiple times and change its orbit due to the applied perturbation, thereby ultimately avoiding collisions with operational satellites.
However, given the laser's wavelength, debris temperature, and material properties, there is a threat of generating specular reflections on the debris surface (\textit{i.e.,} Iridium flares), which scale linearly with the laser power (excluding thermal properties)~\citep{scharring2021ground}. In light of the small forces that photon pressure imparts on the target, its application is constrained to collision avoidance only \citep{scharring2021ground}. 
Alternatively, laser ablation mechanisms rely on the laser's high energy to rapidly melt the target's surface and generate a material jet that produces a reactive momentum, typically higher than the one produced by photon pressure-based lasers~\citep{scharring2021ground}. On the one hand, continuous wave heating produces irregular melt ejection, which can generate more debris. In the case of a tumbling target, the average momentum transfer can be nullified due to the slow heating and decay characteristic of the imparted thrust~\citep{phipps2014laser}.
On the other hand, pulsed lasers appear as a mechanism suitable for collision avoidance and deorbiting debris, given that the per-pulse energy delivered from the ground-based station is enough to change debris orbit \citep{esmiller2014space}. The Laser-Orbital-Debris-Removal, a ground-based laser system, is first introduced in~\cite{phipps_asr_2012} to deorbit small debris and nudge large debris. Further, in~\citep{phipps2014laser}, the ground-based laser system is proposed to be located in polar or equatorial zones to maximize the effectiveness of the remediation mission. To enlarge the number of encounters with small debris, and consequently increase the number of deorbited debris, a system of ground-based lasers is proposed in~\citep{scharring_mdpi_2023}.

Despite their promise as debris remediation solutions, ground-based laser debris remediation systems face significant technical challenges. These systems are subject to numerous atmospheric perturbations, including aerosol attenuation, cloud cover, scintillation effects, and turbulence, all of which affect beam quality~\citep{mason2011orbital}. Moreover, ground-based lasers have limitations in range and angles and thus require strategic positioning of their ground stations to maximize efficiency while considering various civil and operational constraints (\textit{i.e.,} ideally situated away from airports and air routes)~\citep{esmiller2014space}. 
Ultimately, given their ground-based nature, laser engagement opportunities depend on debris passing over the systems' operational range, constraining ground-based lasers to preventive just-in-time collision avoidance without the ability to tackle immediate collision threats~\citep{schall2002laser}.

Against this backdrop, the idea of space-based lasers has garnered increased attention over the past decade, owing to their capacity to address debris of varying sizes (both large and small) and their potential to overcome the inherent challenges associated with ground-based lasers~\citep{schall2002laser, soulard2014ican, phipps2014adroit, bondarenko1997prospects, fang2019effects, bonnal2020just}. Space-based lasers possess several advantages compared to their ground-based counterparts, including more efficient energy delivery and beam quality~\citep{schall2002laser}, increased access to debris~\citep{phipps2014adroit}, application of control laws to adopt a chaser-target formation and maximize laser-to-debris (L2D) engagement effectiveness~\citep{isobe_asr_2024}, and reduced risk of collateral damage due to more precise control enabled by shorter ranges~\citep{scharring2017irregular}. 
According to the NASA Cost and Benefit Analysis of Orbital Debris Remediation report~\citep{colvin2023cost}, ground-based and space-based laser systems are debris remediation methods that can handle both trackable and non-trackable debris pieces with the best cost-to-benefit remediation ratio. Moreover, when not used for remediation purposes, the operator can rely on them to track and characterize debris objects~\citep{colvin2023cost}.
The literature proposes several mission concepts that adopt a single laser platform to remediate orbital debris. A laser system capable of nudging large debris and deorbiting small debris in low Earth orbit (LEO) and geosynchronous Earth orbits (GEO) is proposed in~\citep{phipps2014adroit}. A laser platform is proposed to protect the International Space Station from small debris assumed to be in circular orbits~\citep{schall2002laser}. Further, ~\citep{fang2019effects} incorporates in the $\Delta v$ calculation the inertia matrix of spherical debris objects. Aiming to effectively tackle the large population of debris and its rapid expansion, the literature proposes to use multiple space-based laser platforms. A two-laser system is adopted in~\citep{phippsbonnal} to increase the access rate and decrease the propellant dedicated to debris-chasing maneuvers. Further, a system of multiple lasers is proposed to tackle small debris at different orbital altitudes~\citep{Gambi_AA_2024}.

The existing literature on space-based lasers successfully validates the concept of debris remediation for different sizes; however, several research gaps need to be addressed to make such systems more attractive. First, there is a lack of rigorous mathematical justification for the location of single or multiple laser platform systems. The effectiveness of laser ablation is heavily influenced by the irradiation distance, angle, and revisit time, posing a challenge for laser platforms arbitrarily located in space. Second, the literature lacks an a posteriori analysis of the L2D engagement considering the valuable assets operating near the platforms and the debris field. Nudging debris or placing it on a descent trajectory can alleviate immediate threats to specific valuable assets, however, it may inadvertently create hazardous conjunction events for other objects in space.

In response to these challenges, we propose an optimal constellation of interconnected, collaboratively working space-based lasers, optimally designed to maximize the debris remediation capacity, that is, the ability to nudge and deorbit debris and to perform just-in-time collision avoidance. To materialize this concept, we address the following research questions in this paper:
\begin{enumerate}
\item \textit{``Where do we optimally locate the space-based laser platforms with respect to each other, the debris field, and the valuable assets}?''
\item \textit{``How do we optimally schedule L2D engagements such that the debris remediation capacity is maximized?''}
\end{enumerate}

\begin{figure}[htpb]
    \centering
    \includegraphics[width=0.7\linewidth]{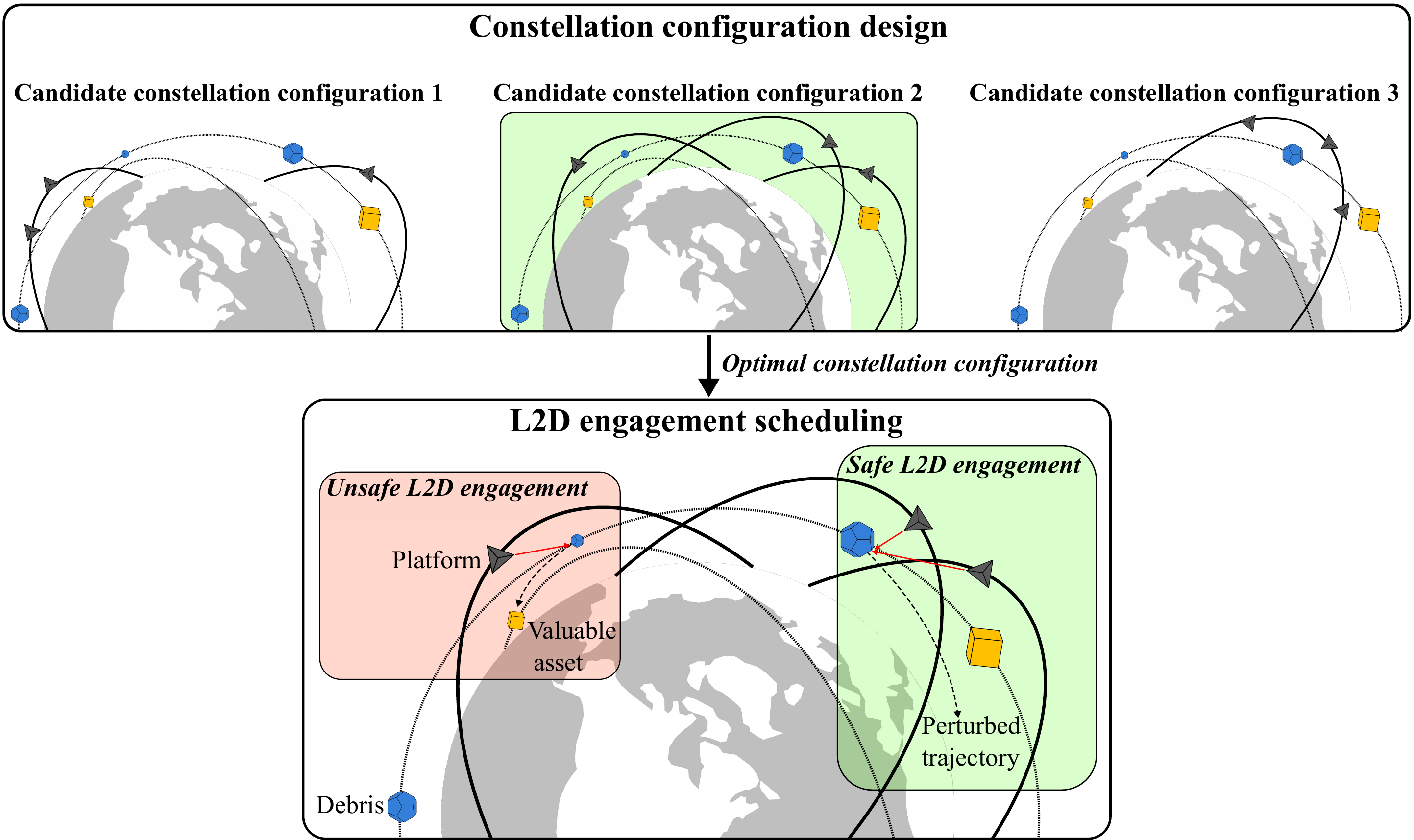}
    \caption{Illustration of the proposed optimization framework.}
    \label{fig:general_ecf}
\end{figure}

Aiming to answer the research questions, we propose to advance the state-of-the-art literature in orbital debris remediation with the following four contributions:
\begin{itemize}
    \item The \textit{$\Delta v$ Vector Analysis} (DVA) is a vector analysis of the $\Delta v$'s imparted by multiple laser platforms. Employing multiple lasers increases the degrees of freedom in the debris control mechanism for remediation missions. This allows for finer control over the target debris resultant $\Delta v$ magnitude and direction, offering a wider range of potential paths for its subsequent descent or collision avoidance trajectory.
    \item The \textit{Concurrent Location-Scheduling Optimization Problem} (CLSP) is an integer linear programming (ILP) formulation that seeks to simultaneously determine the optimal constellation configuration (\textit{i.e.,} the location for a set of laser platforms) while scheduling the L2D engagements leveraging the concept of DVA to maximize the debris remediation capacity. We decompose the CLSP into two problems to tackle the rapid expansion in its solution space.
    \item  The location problem is solved by leveraging the formulation of the \textit{Maximal Covering Location Problem} (MCLP) \citep{church1974maximal,lee2023regional,williams2023FLP}, typically found in the literature on facility location problems to design an optimal constellation configuration that maximizes laser engagement rewards with the debris field. In its original domain, MCLP aims to find the optimal locations for a set of facilities to maximize coverage over a set of customers. 
    We exploit the similarities between satellite constellation pattern design and facility location problems, as discussed in~\citep{williams2023FLP}, interpreting facilities as laser platforms and customers as debris. 
    \item The scheduling problem is addressed by proposing a novel ILP-based formulation, referred to as the \textit{L2D engagement scheduling problem} (L2D-ESP). This formulation determines which combination of lasers should engage with which debris, considering a reward function that accounts for debris mass, change in periapsis radius, and the risk of potential conjunctions with operational satellites before and after the engagements. Third, each optimization formulation takes as input parameters the debris field, valuable assets, the risk of potential conjunctions between them, the number of platforms, and their laser specifications, including pulse energy, pulse length, and operational range. 
\end{itemize}
It is important to note that all optimization formulations presented in this paper are agnostic to the characteristics of the laser platform adopted. This is because the laser specifications serve as an input parameter for each formulation. Further, we do not seek to endorse any specific laser parameters present in the literature but rather outline the value of our contributions. Figure~\ref{fig:general_ecf} outlines the sequential design of an optimal constellation of space-based lasers for orbital debris remediation. First, the constellation configuration design is tackled through the MCLP, which determines, among many combinations, the orbital distribution that maximizes the rewards. The optimal constellation configuration is highlighted in green. Second, the L2D engagement scheduling adopts the optimal constellation configuration and schedules the L2D engagements that maximize the rewards. In addition, it determines which L2D engagements are safe, highlighted in green, and which impose a threat to valuable assets, highlighted in red.

Case studies considering (1) small debris only, (2) large debris only, and (3) a mixed debris field along with 10 valuable assets are presented to illustrate the extension of the formulations. For all three case studies, the formulations presented in this paper outperform Walker-Delta~\citep{walker1984satellite}-based constellations with the same number of platforms in terms of the debris remediation capacity. Additionally, we present a sensitivity analysis to characterize the impact of varying the number of laser platforms on debris remediation capacity.
For the case studies presented in this paper, we adopt the L'ADROIT laser platform~\citep{phipps2014adroit} given its ability to tackle both large and small debris objects. Additionally, to demonstrate the flexibility of the optimization formulations we adopt in~\ref{appendix:new_case_study} a different set of laser parameters, and obtain the optimal constellation configurations for the defined debris field. It is noteworthy to mention that the design and validation of a specific laser system are outside the scope of this paper.

The rest of the paper is structured as follows: Sec.~\ref{sec:parameters_and_framework} describes the technical details of the laser ablation process, DVA, and optimization parameters. Section~\ref{sec:opt_models} presents the three mathematical optimization formulations. Section~\ref{sec:case_studies} illustrates the case studies with debris fields of diverse sizes. Section~\ref{sec:sensitivity} proposes a sensitivity analysis varying the number of platforms. Section~\ref{sec:limitations} outlines the limitations of this paper, and Sec.~\ref{sec:conclusions} presents the conclusions of this paper. 

\section{Laser-to-debris assumptions, modeling, and optimization parameters}\label{sec:parameters_and_framework}
In this paper, space-based laser platforms engage debris by irradiating it with a laser beam, inducing a laser ablation mechanism that alters the debris's orbit. The L2D framework is built upon several assumptions. First, debris is considered a spherical object, with uniform mass distribution, which is always engaged at its center of mass. Second, all L2D engagements are deterministic, meaning that the position and velocity of debris after the engagement are known. Third, each laser platform has complete knowledge of the characteristics of debris (\textit{e.g.,} material, mass, orbit) and whether they endanger the safe operation of valuable assets in space, such as the International Space Station, and in-service telecommunications and Earth-observation satellites. Additionally, we assume that the space-based laser constellation operates as a centralized system with global knowledge. In this system, each platform is aware of the $\Delta v$ magnitude and direction that all other platforms in the constellation can impart over debris, as well as their orbital states during the mission.

The debris remediation capacity of space-based constellations is correlated with the location of laser platforms. Given the distribution of the debris field and its concentration over specific altitude bins, the constellation's target access rate depends on the orbits taken by the platforms. Further, the laser ablation mechanism depends on the relative geometry between the platforms and debris. Hence, an optimal constellation configuration can increase the effectiveness of laser ablation, leading to a higher debris remediation capacity. During the debris remediation mission, each platform can present multiple L2D engagement opportunities. The dynamic scheduling of resources enables the constellation to select the platforms that present the most favorable relative geometry with respect to debris, such that the debris remediation capacity is maximized due to the L2D ablation mechanism.

The mission time horizon is defined as the set of uniformly discretized time steps $\mathcal{T} = \{t_0, \ldots, t_{T-1}\}$, with index $t$ and cardinality $T$, where $t_0$ corresponds to the epoch.
The laser platform set is defined as $\mathcal{P}=\{p_1,\ldots,p_P\}$ with index $p$ and cardinality $P$, and the debris set $\mathcal{D}=\{d_1,\ldots,d_D\}$ with index $d$ and cardinality $D$.

The remainder of Sec.~\ref{sec:parameters_and_framework} is organized as follows. Section~\ref{sec:laser_ablation} describes the laser ablation principle, and Sec.~\ref{sec:DVA} the DVA framework. In Sec.~\ref{sec:network_parameters}, we exhibit the generation of parameters involved in the optimization, and in Sec.~\ref{sec:reward}, we introduce the debris remediation capacity reward.

\subsection{Laser ablation principle}\label{sec:laser_ablation}
Pulsed laser ablation is the process by which a solid ejects plasma from its surface due to the action of a short, intense laser pulse~\citep{lorazo2003short}. At regimes of high irradiance, the vaporized surface material becomes ionized and begins to absorb the incident laser beam, leading to vapor breakdown and plasma formation~\citep{amoruso1999modeling}. The ablation rate, defined as the thickness of the ablated material per laser pulse, is inversely proportional to the square root of the beam diameter, and negatively correlated with laser wavelength given the reduced optical absorptivity and high reflectivity of targets at large wavelengths~\citep{stafe2014pulsed}. Given the velocity of the ejected material, and by the principle of momentum conservation, the ejected mass of the plasma generates a net impulse on the object~\citep{battocchio2020ballistic}. 

Pulsed laser ablation can be leveraged to impart an impulse to debris, causing a change in velocity that alters its orbit. The magnitude of this perturbation is determined by laser parameters such as wavelength, pulse energy, beam quality, and the length and frequency of pulses \citep{settecerri1993laser,wilken2015modelling}, as well as by the properties of debris, including its mass, density, and surface material composition~\citep{scharring2017irregular,liedahl2010momentum, liedahl2013pulsed}. Further, the geometric relationship between the laser and debris determines the direction and magnitude of the velocity perturbation applied to the debris, influencing the final orbit.

The momentum coupling factor $c_\text{m}$ relates the imparted impulse on debris with the used laser's pulse optical energy~\citep{scharring2016numerical}, and is defined as \citep{liedahl2013pulsed}:
\begin{equation}
     m\Delta v = c_\text{m}E
\end{equation}
where $E$ is the laser energy delivered to the debris by a single pulse and $m$ is the debris mass. The on-debris delivered laser fluence $\varphi$ is defined as \citep{phipps2014adroit}:
\begin{equation}
    \varphi = \frac{4ED_{\text{eff}}^2T_{\text{tot}}}{\pi B^4\zeta^2\lambda^2u^2}\label{eq:fluence}
\end{equation}
with $D_{\text{eff}}^2$ being the effective beam diameter, $T_{\text{tot}}$ the total system loss factor, $B$ the beam quality factor, $\zeta$ a constant that regulates diffraction, $\lambda$ the wavelength, and $u$ the range between the platform and debris. Consequently, the per-pulse L2D $\Delta v$ delivered on debris with surface mass $\rho$ is given as:
\begin{equation}
    \Delta v = \eta \frac{ c_\text{m}\varphi}{\rho}\label{eq:DeltaV}
\end{equation}
where $\eta$ is the impulse transfer efficiency, which takes into account shape effects, tumbling, improper thrust direction on debris, and other factors \citep{phipps2014adroit}.

\subsection{$\Delta v$ vector analysis} \label{sec:DVA}
The DVA is a framework that quantifies the effects of multiple simultaneous L2D engagements on target debris by leveraging a vector composition of the imparted $\Delta \bm{v}$ vectors. The applicability of DVA is constrained to two necessary conditions. First, the magnitude of the L2D $\Delta v$ is deterministic and can be computed exactly. Second, the vector components of the $\Delta v$ can be obtained with respect to the adopted reference frame. We let $\bm{r}_{td}^-$ and $\bm{v}_{td}^-$ indicate the position and velocity vectors, respectively, of debris $d$ at time step $t$ immediately before an L2D engagement. Similarly, $\bm{r}_{td}^+$ and $\bm{v}_{td}^+$ are the position and velocity vectors, respectively, of debris $d$ at time step $t$ immediately after an L2D engagement.

We assume that at time step $t$ an L2D engagement induces an instantaneous change in the velocity of target debris, but its position remains unchanged. Hence, we have:
\begin{subequations}
    \begin{alignat}{2}
        &\bm{v}_{td}^+ = \bm{v}_{td}^- + \Delta \bm{v}_{tpd}\\
        &\bm{r}^-_{td} = \bm{r}^+_{td}=\bm{r}_{td}
    \end{alignat}
\end{subequations}
where $\Delta \bm{v}_{tpd}$ represents the change in velocity experienced by debris $d$ due to an L2D engagement from laser platform $p$.

At time step $t$, given position vectors $\bm{r}_{td}$ and $\bm{r}_{tp}$ of debris $d$ and platform $p$, respectively, the relative position vector that points from platform $p$ to debris $d$ is given as $\bm{u}_{tpd} = \bm{r}_{td} - \bm{r}_{tp}$; its unit vector is given as $\hat{\bm{u}}_{tpd} = \bm{u}_{tpd} /{u}_{tpd}$, with ${u}_{tpd}=\left \lVert \bm{u}_{tpd} \right \rVert_2$. The total change in velocity imparted from laser platform $p$ to debris $d$ at time step $t$ is calculated by generalizing Eq.~\eqref{eq:DeltaV} for every time step $t$, laser platform $p$, and debris $d$. Equation~\eqref{eq:DeltaVsd} presents the vector $\Delta \bm{v}_{tpd}$ from laser platform $p$ to debris $d$ at time step $t$ given unit vector $\hat{\bm{u}}_{tpd}$.
\begin{equation}
    \Delta \bm{v}_{tpd} = N_d\frac{\eta_{pd} \varphi_{pd} c_{\text{m},pd}}{\rho_d} \hat{\bm{u}}_{tpd} \label{eq:DeltaVsd}
\end{equation}
where $N_d$ is the number of laser pulses per time step $t$ on debris $d$ computed from the time step size and the laser pulse repetition frequency (PRF), which denotes the number of laser pulses per second.

The DVA framework captures multiple, simultaneous L2D engagements and represents them as a total effective single L2D engagement. The total $\Delta \bm{v}_{td}$ over debris $d$ at time step $t$ due to multiple L2D engagements can be represented as follows:
\begin{equation}
    \Delta \bm{v}_{td} = \sum_{p\in\mathcal{P}_{td}} \Delta \bm{v}_{tpd}\label{eq:dva_full}
\end{equation}
where $\mathcal{P}_{td}$ represents the set of laser platforms that engage debris $d$ at time step $t$. To determine how $\Delta \bm{v}_{td}$ affects debris trajectory, the new debris orbital parameters (\textit{i.e.,} semi-major axis, eccentricity) can be analytically computed using $\bm{r}_{td}$ and $\bm{v}_{td}^+ = \bm{v}_{td}^- + \Delta \bm{v}_{td}$ and compared with those before the L2D engagement.

\begin{figure}[htb]
    \centering
    \includegraphics[width=0.3\textwidth]{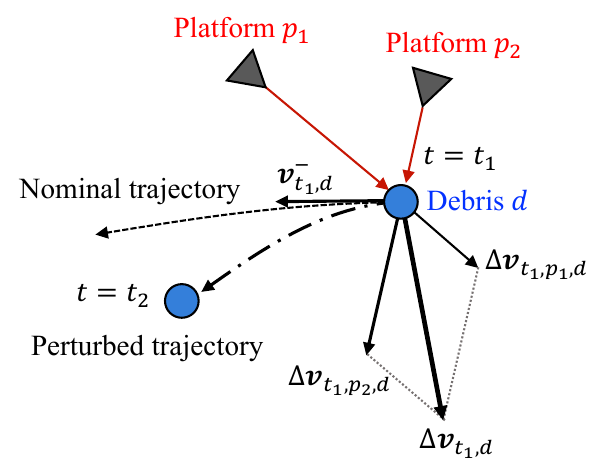}
    \caption{DVA illustrative example with two laser platforms $p_1$ and $p_2$ and debris $d$.}
    \label{fig:DVA}
\end{figure}

The significance of DVA lies in its ability to enable a higher degree of control over debris dynamics by controlling the imparted $\Delta \bm{v}$ in both magnitude and direction, thus leading to a new debris orbit. Figure~\ref{fig:DVA} illustrates an example involving two laser platforms acting on the same debris at a given time step. In this specific scenario, there are three different perturbed orbits available for the debris: (1) one acted solely by platform $p_1$, (2) another acted solely by platform $p_2$, and (3) one acted collaboratively by both platforms simultaneously, resulting from a combined $\Delta \bm{v}_{{t_1},d} = \Delta \bm{v}_{{t_1},{p_1},d} + \Delta \bm{v}_{{t_1},{p_2},d}$. The availability of multiple L2D engagement opportunities is advantageous, providing the space-based laser constellation with greater flexibility to achieve more effective debris remediation.

\subsection{Constellation optimization parameters}\label{sec:network_parameters}
As a consequence of an L2D engagement, debris can relocate to a new orbit. We define $\mathcal{J}_{td}$ as the set of orbital slots, with index $j$ and cardinality $J_{td}$, for debris $d$ at time step $t$. Each orbital slot $j \in \mathcal{J}_{td}$ is associated with the position and velocity vectors defined at time step $t$ to fully define the state of debris orbit. Laser platform orbital slots are defined in set $\mathcal{S}$, indexed by $s$ and with cardinality $S$. Each laser platform is assumed to take an orbital slot in $\mathcal{S}$ and maintain it during the entire mission. Moreover, we define a set of valuable assets $\mathcal{K}$ with cardinality $K$.

The feasibility of the L2D engagements is encoded with the Boolean parameter $W_{tsd}$, given as:
\begin{equation}
    W_{tsd}=\begin{cases}
        1, &\text{if a platform located at orbital slot $s$ can engage debris $d$ at time step $t$}\\
        0, &\text{otherwise}
    \end{cases}
\end{equation}

\begin{figure}[htb]
    \centering
    \includegraphics[width=0.25\linewidth]{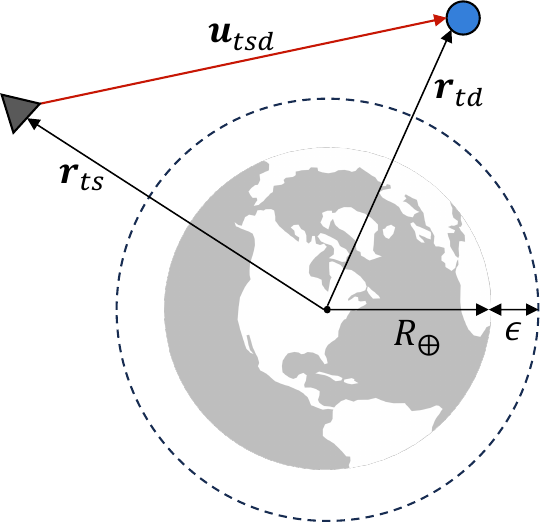}
    \caption{Illustration of Eq.~\eqref{eq:arg_geom} when debris is in the line-of-sight of a platform.}
    \label{fig:engagement_condition}
\end{figure}

To determine if an L2D engagement is feasible, two conditions must be satisfied. First, debris $d$ has to be in the line-of-sight of a laser platform located at orbital slot $s$ at time step $t$. Given the relative position vector $\bm{u}_{tsd}$ and the range $u_{tsd} = ||\bm{u}_{tsd}||_2$, Eq.~\eqref{eq:arg_geom} computes the line-of-sight indicator \citep{williams2023FLP}. 
\begin{equation}
q_{tsd} = \left((r_{ts})^2-(R_{\oplus}+\epsilon)^2\right)^{1/2}+\left((r_{td})^2-(R_{\oplus}+\epsilon)^2\right)^{1/2}-u_{tsd}\label{eq:arg_geom}
\end{equation}
where $R_{\oplus}$ is the radius of the Earth and $\epsilon$ is a bias parameter. If $q_{tsd}>0$, debris is in the line-of-sight of a platform, as illustrated in Fig.~\ref{fig:engagement_condition}. Conversely, for $q_{tsd}=0$, $u_{tsd}$ is tangential to the dashed line-sphere, and for $q_{tsd}<0$, it intersects the dashed line-sphere in two points; consequently, in both cases, debris is not in the line-of-sight of the platform.
Second, an L2D engagement can happen if $u_{tsd}$ lies within the maximum and minimum operational ranges defined as $u_{\max}^d$ and $u_{\min}^d$, respectively. The given laser specifications determine the upper bound, and the lower bound is set to achieve a safe operation of the system. Leveraging parameter $W_{tsd}$, it is possible to generate for each time step $t$ the set of laser platform orbital slots $\mathcal{S}_{tdj}$, with index $s$ and cardinality $S_{tdj}$, that engage debris $d \in \mathcal{D}$ and generate debris orbital slot $j \in \mathcal{J}_{t+1,d}$.

\subsection{Debris remediation capacity reward}\label{sec:reward}
We introduce debris remediation capacity reward $R_{tdij}$ to quantify the value of relocating debris $d$ from its current orbital slot $i \in \mathcal{J}_{td}$ to a new orbital slot $j \in \mathcal{J}_{t+1,d}$ at time step $t$, as a consequence of an L2D engagement while considering potential conjunctions with the set of valuable assets $\mathcal{K}$. The debris remediation capacity reward is defined as:
\begin{equation}
    R_{tdij} = {C}^0_{td} + {C}_{tdij} + \alpha\Delta h_{tdij}  + \beta{M}_{td} 
    \label{eq:reward}
\end{equation}

\begin{figure}[htbp]
    \centering
    \begin{subfigure}[b]{.49\textwidth}
        \centering
        \includegraphics[width=.6\linewidth]{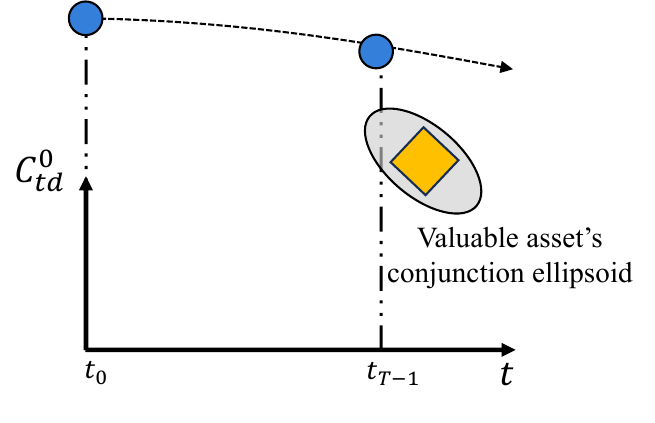}
        \caption{Mission scenario with ${C}^0_{td}=0$.}
        \label{fig:c0tdi0}
     \end{subfigure}
     \begin{subfigure}[b]{.49\textwidth}
        \centering
        \includegraphics[width=.6\linewidth]{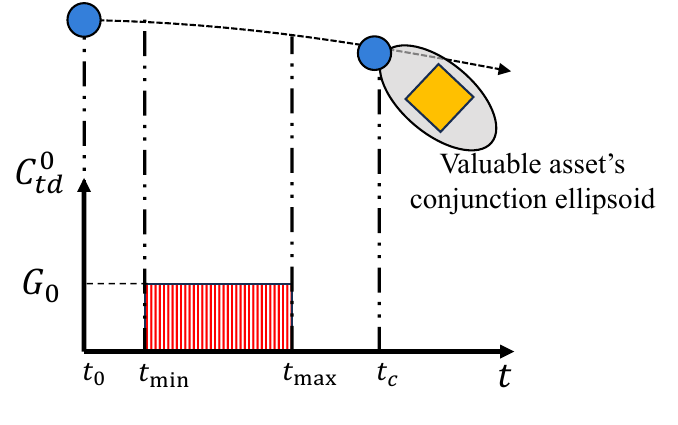}
        \caption{Mission scenario with ${C}^0_{td}=G_0$ for all $ t \in [t_{\min}, t_{\max}]$.}
        \label{fig:c0tdiM}
    \end{subfigure}
    \caption{Generation of ${C}^0_{td}$ reward for two different scenarios.}
    \label{fig:c0tdi}
\end{figure}

First, ${C}^0_{td}$ is an incentive term that accounts for the conjunction analysis between debris $d$ and the set of valuable assets $\mathcal{K}$ during the entire mission scenario, assuming no L2D engagements. Debris orbit is examined if it lies inside the deterministic conjunction ellipsoid of any valuable asset in $\mathcal{K}$. We define ${C}^0_{td}$ as:
\begin{equation}
    {C}^0_{td} = 
    \begin{cases}
        G_0, &\text{if $t \in [t_{\min}, t_{\max}]$}\\
        0, &\text{otherwise}
    \end{cases}
\end{equation}
If a conjunction, or multiple ones, is feasible between debris $d$ and any valuable asset in $\mathcal{K}$,  ${C}^0_{td} = G_0$ for all $ t \in [t_{\min}, t_{\max}]$ where $G_0$ is a large positive constant that incentivizes the constellation to engage debris during the time window $[t_{\min}, t_{\max}]$. The length of the interval and its starting point are parameters to be defined by a user considering the characteristics of the space-based laser, its capability for changing debris orbit, and the desired minimum miss distance. Given $t_c$, the time step at which the first conjunction occurs (if debris $d$ is to have multiple conjunctions with the same valuable asset or with several ones), $t_{\max}$ should be defined such that $t_{\max} \le t_c$. Figure~\ref{fig:c0tdi0} illustrates a scenario where reward ${C}^0_{td}$ is not activated since there is no feasible conjunction. Figure~\ref{fig:c0tdiM} outlines the case where a conjunction is imminent and the reward is activated.

\begin{figure}[htbp]
    \centering
    \begin{subfigure}[b]{0.45\textwidth}
        \centering
        \includegraphics[width=0.7\linewidth]{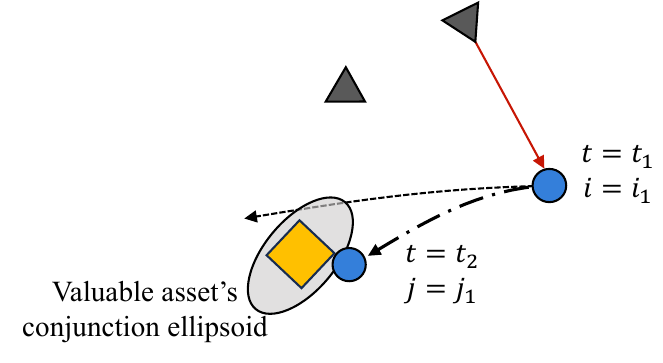}
        \caption{Mission scenario with activated penalty ${C}_{tdij}=-G$.}
        \label{fig:ctdji_M}
    \end{subfigure}
    \begin{subfigure}[b]{0.45\textwidth}
        \centering
        \includegraphics[width=0.7\linewidth]{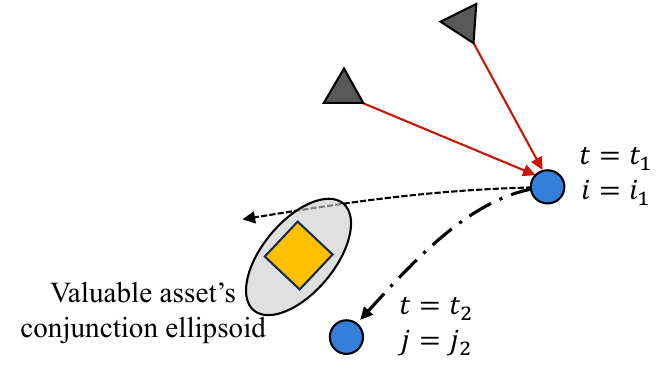}
        \caption{Mission scenario with inactive penalty ${C}_{tdij}=0$.}
        \label{fig:ctdji_0}
     \end{subfigure}
    \caption{Generation of ${C}_{tdij}$ penalty for two different engagement scenarios.}
    \label{fig:ctdji}
\end{figure}

Second, ${C}_{tdij}$ is a look-a-head penalty term that accounts for whether relocating debris $d$ from orbital slot $i$ to orbital slot $j$ due to an L2D engagement at time step $t$ triggers at least one conjunction with any valuable asset in $\mathcal{K}$. Once all debris candidate orbital slots are generated for step $t + 1$, each orbital slot is propagated from $t+1$ to $t+1+\tau$, where $\tau$ represents the number of look-ahead time steps. We compare for every $t \in [t+1, t+1 +\tau]$ if the range between debris $d$ at slot $j$ with any valuable asset is less than the range threshold. If a specific L2D engagement with debris generates a conjunction with any valuable asset defined in $\mathcal{K}$, we set ${C}_{tdij}=-G$ to ensure that the constellation is discouraged from relocating debris $d$ to orbital slot $j$, as illustrated in Fig.~\ref{fig:ctdji_M}. Conversely, Fig.~\ref{fig:ctdji_0} highlights two simultaneous L2D engagements where the resultant debris orbital slot does not invade the valuable asset's conjunction ellipsoid; hence ${C}_{tdij}$ is not activated.

Third, $\Delta h_{tdij}$ is a look-a-head incentive term that compares the periapsis radius of debris $d$ orbital slot $j$ after an L2D engagement at time step $t$ with the periapsis radius threshold $h^*$, at which we consider that debris is governed by atmospheric drag and subject to induced natural orbital decay. Once this decay is initiated, the debris is considered deorbited within the framework. Given $h_{t+1,dj}$, the periapsis radius of the new orbital slot $j$ for debris $d$ after an L2D engagement, the reward term is defined as:
\begin{equation}
    \Delta h_{tdij} = \gamma_{tdij}\left(\frac{h^*}{h_{t+1,dj}}\right)^3
\end{equation}
where $\gamma_{tdij}$ is a variable that compares the periapsis radius of debris $d$ between initial orbital slot $i$ at time step $t$ and after the L2D engagement at final orbital slot $j$ at time step $t+1$, defined as:
\begin{equation}
    \gamma_{tdij} = 
    \begin{cases}
        -G_h, &\text{if $h_{t+1,dj} > h_{tdi}$}\\
        1, &\text{otherwise}
    \end{cases}
\end{equation}
If an L2D engagement increases the periapsis radius, then a large negative constant $G_h$ is assigned to $\gamma_{tdij}$ to discourage the engagement. Furthermore, given that $\Delta h_{tdij}$ is a parameter that captures the ratio between the periapsis radius threshold and debris orbit periapsis radius after the engagement, a value of $\Delta h_{tdij} \ge 1$ implies that the object is deorbited; hence, we constrain it to be at most one to avoid obtaining large rewards due to the nature of the terms' formulation. Additionally, we scale it by $\alpha$ in Eq.~\eqref{eq:reward} to give more flexibility to the user.

Lastly, ${M}_{td}$ accounts for debris $d$ mass, which is activated if debris is engaged at time step $t$, and whose objective in the reward formulation is two-fold. First, debris of large mass is identified as the main source of new debris in case of a collision, and therefore in alignment with Refs.~\citep{mcknight2021} and~\citep{liou2011active}, which assign a debris remediation capacity reward to debris proportional to its mass, we aim to assign a bigger reward for engaging debris of larger mass. Second, given Eq.~\eqref{eq:DeltaVsd}, it is clear that for the same laser and $c_\text{m}$ parameters, debris of larger mass will have a smaller $\Delta v$ than those of smaller mass, making it harder to achieve significant changes in its trajectory, resulting in a smaller $\Delta h_{tdij}$. To address the two key points stated above, we introduce the mass reward term, defined as:
\begin{equation}
    {M}_{td} = 
    \begin{cases}
        m_d/m_{\max}, &\text{if debris $d$ is engaged at step $t$}\\
        0, &\text{otherwise}
    \end{cases}
\end{equation}
where $m_{\max}$ is the maximum debris mass in $\mathcal{D}$. We normalize the mass of debris $d$ by dividing it by the maximum debris mass to be consistent in terms of magnitude with $\Delta h_{tdij}$. Further, we introduce $\beta$ in Eq.~\eqref{eq:reward} to scale it with respect to other parameters.

\section{Constellation configuration design and scheduling optimization} \label{sec:opt_models}

The placement and scheduling of an optimal constellation of space-based lasers for debris remediation is tackled with an ILP formulation. The CLSP concurrently optimizes the location and scheduling of the laser platforms given a specific debris field; however, due to the rapid expansion of its solution space is further decomposed into two problems. First, the location of the laser platforms is tackled with the MCLP; second, the scheduling of the L2D engagements is tackled with the L2D-ESP.

This section is structured as follows. Section~\ref{sec:allinone} presents the CLSP formulation with an illustrative example. Section~\ref{sec:MCLP} outlines the application of the MCLP to locate the laser platforms. Section~\ref{sec:scheduling} introduces the L2D-ESP.

\subsection{The concurrent location-scheduling optimization problem}\label{sec:allinone}
The optimal location and scheduling of a constellation of space-based lasers is concurrently tackled by proposing the CLSP. This problem seeks to determine, for the entire mission planning horizon, the optimal location for a set of platforms while considering an optimal L2D engagement scheduling to maximize the debris remediation capacity. We report all sets, parameters, and variables used in the CLSP mathematical formulation in Table~\ref{table:clsp}.

\begin{table*}[htb]
\caption{Sets, parameters, and variables for the CLSP formulation.}
\centering
\begin{tabularx}{\textwidth}{llX}
\hline
Type & Symbol & Description \\
\hline
Sets &$\mathcal{T}$ & Mission planning horizon (index $t$; cardinality $T$) \\
    &$\mathcal{S}$ & Set of orbital slots for laser platforms (index $s$; cardinality $S$)\\
    &$\mathcal{D}$ & Set of debris (index $d$; cardinality $D$) \\
   &$\mathcal{J}_{td}$ & Set of debris $d$ orbital slots at time step $t$ (index $j$; cardinality $J_{td}$) \\
    &$\mathcal{S}_{tdj}$ & Set of laser platform orbital slots that engage debris $d$ at time step $t$ and generate orbital slot $j$ (index $s$; cardinality ${S}_{tdj}$) \\
Parameters     &$R_{tdij}$ & Debris remediation capacity reward for relocating debris $d$ from orbital slot $i$ to orbital slot $j$ at time step $t$ \\
                &$W_{tsd}$ & $\begin{cases}
                     1, &\text{if debris $d$ can be engaged by a platform located at orbital slot $s$ at time step $t$} \\
                    0, &\text{otherwise}
                    \end{cases}$ \\
                &$P$&Number of platforms in the constellation\\
Decision variables &  $z_{s}$ & $\begin{cases}
                     1, &\text{if a platform is located at orbital slot $s$} \\
                    0, &\text{otherwise}
                    \end{cases}$ \\
                    &  $y_{tsd}$ & $\begin{cases}
                     1, &\text{if a platform located at orbital slot $s$ engages debris $d$ at time step $t$} \\
                    0, &\text{otherwise}
                    \end{cases}$ \\
                    &  $x_{tdij}$ & $\begin{cases}
                     1, &\text{if debris $d$ relocates from orbital slot $i$ to orbital slot $j$ at time step $t$} \\
                    0, &\text{otherwise}
                    \end{cases}$ \\
\hline
\label{table:clsp}
\end{tabularx}
\end{table*}

\subsubsection{Decision variables}
The placement of the laser platforms in the candidate orbital slots is modeled using the following binary \textit{platform location decision variables}:
\begin{equation}
z_{s} = \begin{cases}
    1, & \text{if a laser platform is located at orbital slot $s$} \\
    0, & \text{otherwise}
\end{cases}
\end{equation}
Similarly, the L2D engagements for each platform during the mission are scheduled with the following binary \textit{L2D engagement decision variables}:
\begin{equation}
y_{tsd} = \begin{cases}
    1, & \text{if a laser platform located at orbital slot $s$ engages debris $d$ at time step $t$} \\
    0, & \text{otherwise}
\end{cases}
\end{equation}

A debris orbit remains unchanged for successive time steps if no L2D engagement occurs. Conversely, the ablated debris takes a new orbit as a consequence of an L2D engagement by a single or multiple platforms, as the engagement defines the new debris' state vector. In CLSP, the change in debris orbit during the mission is encoded with binary \textit{debris relocation decision variables}, given as:
\begin{equation}
x_{tdij} = \begin{cases}
    1, & \text{if debris $d$ relocates from orbital slot $i$ to orbital slot $j$ at time step $t$} \\
    0, & \text{otherwise}
\end{cases}
\end{equation}
At each time step $t$, debris $d$ is enforced to take a new orbital slot $j \in \mathcal{J}_{t+1,d}$. The set of debris orbital slots can be subdivided into two subsets. First, $\Tilde{\mathcal{J}}_{t+1,d} \subseteq \mathcal{J}_{t+1,d}$ with index $j$ and cardinality $\Tilde{J}_{t+1,d}$ contains the orbital slots that encode no changes in the current orbit for debris $d$ since it is not engaged, and consequently they have no associated reward. Second, the set of new orbits originated as a consequence of L2D engagements, where each one has its corresponding reward $R_{tdij}$, encoded as $\mathcal{J}_{t+1,d} \setminus \Tilde{\mathcal{J}}_{t+1,d}$ with index $j$. 

\subsubsection{Constraints and objective function}
We introduce constraints~\eqref{local:location_coupling} to enforce for each time step that a laser platform located at orbital slot $s$ engages debris $d$ only if it satisfies the L2D engagement requirements encoded in parameter $W_{tsd}$, and if the orbital slot is effectively taken. Constraints~\eqref{local:eng_limit} limit each platform to engage at most one debris per time step.
\begin{subequations}
\begin{alignat}{2}
&W_{tsd}z_s \ge y_{tsd}, &\quad \forall t\in\mathcal{T}, \forall s \in \mathcal{S},\forall d\in\mathcal{D}\label{local:location_coupling}\\
&\sum_{d\in\mathcal{D}} y_{tsd} \le 1, &\quad \forall t\in\mathcal{T}, \forall s \in \mathcal{S} \label{local:eng_limit}
\end{alignat}
\end{subequations}

We define path contiguity constraints to ensure the correct flow of debris, \textit{i.e.,} a feasible evolution of debris orbit, throughout the entire mission horizon. First, we require all debris to occupy a new orbital slot at time step $t_1$ with constraints~\eqref{local:flow_relocate}. Either indicating that debris has not changed its orbit and relocates to an orbital slot $j \in \Tilde{\mathcal{J}}_{t_1,d}$, or indicating that it has been engaged by at least one platform and relocates to an orbital slot $j \in \mathcal{J}_{t_1,d}\setminus\Tilde{\mathcal{J}}_{t_1,d}$ which encodes the state vector of a new orbit. Second, constraints~\eqref{local:flow_balance} ensure that debris $d$ follows a continuous orbital slot path between time steps $t-1$ and $t+1$ with respect to the intermediate reference slot $i \in \mathcal{J}_{td}$ at time step $t$.
\begin{subequations}
    \begin{alignat}{2}
&\sum_{j\in\mathcal{J}_{t_1,d}} x_{t_{0},d,i_{0},j} = 1, &\quad \forall d \in \mathcal{D}\label{local:flow_relocate}\\
&\sum_{j\in\mathcal{J}_{t+1,d}} x_{tdij} -  \sum_{\upsilon\in\mathcal{J}_{t-1,d}} x_{t-1,d\upsilon i} = 0,&\quad \forall t\in\mathcal{T}\setminus\{t_0,t_{T-1}\}, \forall d \in \mathcal{D}, \forall i \in \mathcal{J}_{td}\label{local:flow_balance}
\end{alignat}
\end{subequations}

Constraints~\eqref{local:engagements} are linking constraints that couple debris relocation decision variables $x_{tdij}$ and L2D engagement decision variables $y_{tsd}$. A relocation to orbital slot $j \in \mathcal{J}_{t+1,d}\setminus\Tilde{\mathcal{J}}_{t_1,d}$ (that encodes a new orbit) is only possible if debris $d$ is engaged by the set of laser platforms located in orbital slots $\mathcal{S}_{tdj}$. We relax $\Tilde{\mathcal{J}}_{t+1,d}$ given that it has an associated $\mathcal{S}_{tdj}=\emptyset$ and debris is enforced to occupy only one new orbital slot $j$ as imposed by constraints~\eqref{local:flow_relocate} and~\eqref{local:flow_balance}.
\begin{equation}
    \sum_{s\in\mathcal{S}_{tdj}} y_{tsd} \ge S_{tdj}x_{tdij}, \quad  \forall t\in\mathcal{T}\setminus\{t_{T-1}\},\quad \forall d\in\mathcal{D}, \forall i \in \mathcal{J}_{td},\forall j \in \mathcal{J}_{t+1,d}\setminus{\Tilde{\mathcal{J}}_{t+1,d}} \label{local:engagements}
\end{equation}

Additionally, constraint~\eqref{local:cardinality} enforces the number of platforms in the constellation to be exactly $P$.
\begin{equation}
\sum_{s\in\mathcal{S}} z_{s} = P \label{local:cardinality}
\end{equation}

The domains of all decision variables are given as follows:
\begin{subequations}
    \begin{alignat}{2}
        & z_s \in \{0,1\}, &\quad \forall s\in \mathcal{S} \label{local:x_s} \\
        & y_{tsd} \in \{0,1\}, &\quad\forall t\in \mathcal{T},\forall s\in \mathcal{S}, \forall d \in \mathcal{D} \label{local:y_sd}\\
        & x_{tdij} \in \{0,1\}, &\quad\forall t\in \mathcal{T}\setminus\{t_{T-1}\},\forall d \in \mathcal{D}, \forall i \in \mathcal{J}_{td},\forall j \in \mathcal{J}_{t+1,d}  \label{local:x_dkj}
\end{alignat}
\end{subequations}

Lastly, the formulation's objective is to maximize the remediation capacity of the constellation, which is given as:
\begin{equation}
V=\sum_{t\in\mathcal{T}\setminus\{t_{T-1}\}}\sum_{d\in\mathcal{D}}\sum_{i\in\mathcal{J}_{td}}\sum_{j\in\mathcal{J}_{t+1,d}} R_{tdij}x_{tdij} \label{local:obj}
\end{equation}

\subsubsection{Mathematical formulation}
Piecing everything together, the mathematical formulation for CLSP is given as follows:
\begin{alignat*}{2}
\max \quad & \text{Objective function~\eqref{local:obj}} \\
\text{s.t.} \quad & \text{Constraints \eqref{local:location_coupling} to \eqref{local:x_dkj}}
\end{alignat*}

\subsubsection{Illustrative example}
\begin{figure}[htp]
    \centering
    \includegraphics[width=0.4\textwidth]{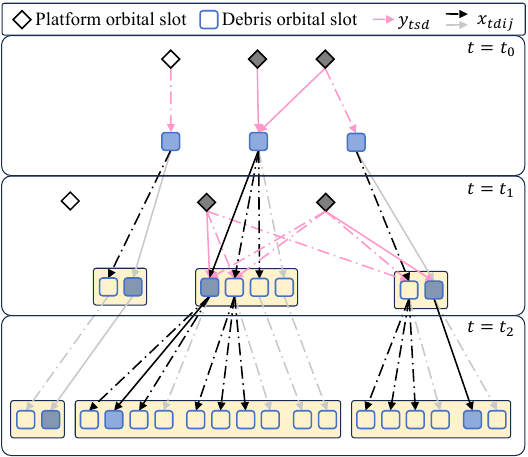}
    \caption{Illustration of the solution space for CLSP.}
    \label{fig:local}
\end{figure}

To illustrate CLSP, we present a small-scale instance with $S=3$, $T=3$, $D=3$, and $P=2$. Figure~\ref{fig:local} outlines the tree structure of the problem with platform orbital slots $s$ represented with rhombus, and debris orbital slots $i,j$ with squares. The magenta lines represent feasible decision variables $y_{tsd}$, the gray lines are debris $d$ feasible relocation variables $x_{tdij}$ at time step $t$ from orbital slot $i\in \mathcal{J}_{td}$ to the set of orbital slots $j \in\Tilde{\mathcal{J}}_{t+1,d}$, and the dark lines are debris feasible relocation variables $x_{tdij}$ to orbital slots $j \in \mathcal{J}_{t+1,d}\setminus{\Tilde{\mathcal{J}}_{t+1},d}$. Unfeasible decision variables are given a large negative reward $R_{tdij}$ and are not displayed in the figure.

In Fig.~\ref{fig:local}, debris is initially located at orbital slot $i \in \mathcal{J}_{t_0,d}$ and has an associated virtual orbital slot $j \in\Tilde{\mathcal{J}}_{t_1,d}$, which represents no change in orbit if no L2D engagement occurs, consequently, $R_{t_{0},dij}=0$. The subset of debris orbital slots that represent new orbits depends on the number of feasible L2D engagements, for instance, while the leftmost debris can be engaged only by the leftmost platform orbital slot, debris located in the center orbital slot can be engaged by at most two platforms, hence it has three new orbits, all of them with its corresponding reward $R_{t_{0},dij}$. Given that, at time step $t_1$, the number of total debris orbital slots is eight and that the number of possible L2D engagements is seven, the total number of debris orbital slots at time step $t_2$ expands to 18.

The problem has 94 integer binary variables and 85 constraints. It is solved using the Gurobi optimizer 11, which retrieves an objective value (\textit{i.e.,} debris remediation capacity reward obtained) of \num{7.50} with a duality gap of \SI{0.00}{\%}. In Fig.~\ref{fig:local}, the selected platform orbital slots are colored in dark gray, the debris orbital slots in blue, and decision variables $y_{tsd}=1$ and $x_{tdij}=1$ are straight lines, conversely, feasible decision variables equal to zero are represented with dashed lines. 

Even though the problem is presented for a small number of time steps and orbital slots, the number of decision variables and constraints grows exponentially. Figure~\ref{fig:local} accurately showcases the tree structure of the problem and how it rapidly expands. Scaling this formulation for larger mission horizons with more debris and platform orbital slots makes the problem computationally prohibitive, restricting its application to small-scale problems. To overcome this challenge, we propose two separate optimization formulations. First, we find the optimal laser constellation configuration design using MCLP. Then, we leverage it to optimize the debris remediation capacity using the L2D-ESP.

\subsection{Optimal laser constellation configuration design}\label{sec:MCLP}
The problem of designing an optimal laser constellation configuration is tackled under the assumption that debris does not relocate; consequently, decision variables $x_{tdij}$ are relaxed, that is, removed, along with constraints~\eqref{local:flow_relocate}, \eqref{local:flow_balance}, \eqref{local:engagements} and~\eqref{local:x_dkj}. The resulting formulation inherits the location-scheduling structure from CLSP with location decision variables $z_s$, L2D engagement decision variables $y_{tsd}$, and an orbital slot-independent reward, given as:
\begin{equation}
    R_{td} = {C}^0_{td} + {M}_{td}
    \label{eq:loc_reward}
\end{equation}

The constellation configuration reward $R_{td}$ accounts for feasible conjunctions between debris $d$ and any valuable asset in $\mathcal{K}$ leveraging parameter ${C}^0_{td}$ and for the mass of the engaged debris $d$ with parameter $M_{td}$, both introduced in Eq.~\eqref{eq:reward}. However, the new reward $R_{td}$ cannot assess whether that L2D potential interaction will effectively lower debris periapsis radius because the $\Delta h_{tdij}$ term from Eq.~\eqref{eq:reward} is relaxed. To address this problem, we introduce the Boolean parameter $W_{tsd}'$, defined as:
\begin{equation}
    W'_{tsd}=\begin{cases}
        1, &\text{if $W_{tsd}=1$ and $h_{t+1,dj} \le h_{tdi}$}\\
        0, &\text{otherwise}
    \end{cases}
\end{equation}

Given the nature of reward $R_{td}$, the new formulation can prioritize the selection of platform orbital slots closer to larger debris and neglect orbital slots that can engage small debris objects, since the change in debris periapsis radius is not considered in the reward. Further, the resulting configuration, concentrated over larger debris, can considerably reduce the L2D engagement opportunities during the L2D-ESP described in Sec.~\ref{sec:scheduling}. To tackle this problem, we propose to leverage MCLP to design an optimal laser constellation configuration. We refer the reader to~\ref{appendix:MCLP} for an in-depth justification. Table~\ref{tab:MCLP} presents the new parameters and decision variables introduced for the MCLP formulation.

\begin{table}[htb]
\caption{Parameters and variables specific to the MCLP formulation.}
\label{tab:MCLP}
\centering
\begin{tabular}{lll}
\hline
Type & Symbol & Description \\
\hline
Parameters
    &$R_{td}$ & Constellation configuration reward for debris $d$ at time step $t$ \\
    &$S_{td}$ & Engagement requirement for debris $d$ at time step $t$ ($S_{td}\in\mathbb{Z}_{\ge1}$) \\
    &$W'_{tsd}$ & $\begin{cases}
      1, & \text{if orbital slot $s$ can engage debris $d$ at time step $t$} \\
      0, & \text{otherwise}
    \end{cases}$ \\
Decision variables 
    &$x_{td}$ & $\begin{cases}
        1, &\text{if debris $d$ is engaged at time step $t$} \\
        0, &\text{otherwise}
    \end{cases}$\\
\hline
\end{tabular}
\end{table}

\subsubsection{Decision variables}
To formulate our problem as an MCLP, we keep platform location decision variables $z_s$ and introduce new decision variables $x_{td}$ defined as:
\begin{equation}
x_{td} = \begin{cases}
    1, & \text{if debris $d$ is engaged at time step $t$} \\
    0, & \text{otherwise}
\end{cases}
\end{equation}

These new decision variables are distinct from CLSP L2D engagement decision variables $y_{tsd}$ in a manner that they encode if debris is engaged or not and do not discriminate which platforms perform L2D ablation over debris.

\subsubsection{Constraints and objective function}
In light of the fact that no L2D engagement can occur between a platform located at orbital slot $s$ and debris $d$ if orbital slot $s$ is not taken, linking constraints~\eqref{loc:coupling} between decision variables $z_s$ and $x_{td}$ are imposed. Decision variables $x_{td}$ are activated if at least $S_{td} \in\mathbb{Z}_{\ge1}$ platforms engage debris $d$ at time step $t$.
\begin{equation}
\sum_{s\in\mathcal{S}} W'_{tsd}z_s\ge S_{td} x_{td}, \quad \forall t\in \mathcal{T}, \forall d \in \mathcal{D} \label{loc:coupling}
\end{equation}
Here, $S_{td}$ is a parameter inspired by the cardinality of set $\mathcal{S}_{tdj}$ used in constraints~\eqref{local:engagements}, which imposes the required platforms to engage debris $d$ in order to relocate it to orbital slot $j$. In this formulation, $S_{td}$ aims to impose a lower bound on the number of platforms required to engage debris $d$ at time step $t$. This requirement can be beneficial in scenarios where large debris is targeted and a large $\Delta \bm{v}$ is required to generate a significant change in debris orbit.

The number of platforms to be used in the constellation is enforced by constraints~\eqref{local:cardinality}, and the decision variables' domain in the location formulation is defined as:
\begin{subequations}
\begin{alignat}{2}
    & z_s \in \{0,1\}, & \forall s\in \mathcal{S} \label{loc:z_s} \\
    & x_{td} \in \{0,1\}, & \quad\forall t\in \mathcal{T}, \forall d \in \mathcal{D} \label{loc:y_td}
\end{alignat}
\end{subequations}

The MCLP formulation aims to maximize the realized constellation configuration reward during the mission time horizon $\mathcal{T}$, which is encoded in objective function~\eqref{loc:obj}:
\begin{equation}
   \pi = \sum_{t\in\mathcal{T}}\sum_{d\in \mathcal{D}}R_{td} x_{td} \label{loc:obj}
\end{equation}
If an L2D engagement exists for debris $d$ at time step $t$, then $x_{td}=1$ and reward $R_{td}$ associated is obtained.  It is important to mention that if $R_{td} = 1, \forall t \in \mathcal{T}, \forall d \in \mathcal{D}$, the problem will yield the total number of potential engagements that the constellation has at the beginning of the L2D-ESP.

\subsubsection{Mathematical formulation}
Piecing all constraints and the objective function together, the resulting MCLP optimization formulation that determines the location of the laser platforms to maximize the constellation configuration reward is given as:
\begin{alignat*}{2}
\max \quad & \text{Objective function~\eqref{loc:obj}} \\
\text{s.t.} \quad & \text{Constraints \eqref{local:cardinality}, \eqref{loc:coupling}, \eqref{loc:z_s} and \eqref{loc:y_td}}
\end{alignat*}

\subsection{Optimal laser-to-debris engagement scheduling problem}\label{sec:scheduling}
In this section, we build on the assumption that the location of platforms is given, and we maintain the engagement constraints from CLSP to derive the L2D-ESP. Relaxing constraints~\eqref{local:location_coupling} and \eqref{local:cardinality}, the L2D-ESP is:
\begin{alignat*}{2}
\max \quad & \text{Objective function~\eqref{local:obj}} \\
\text{s.t.} \quad & \text{Constraints \eqref{local:eng_limit}, \eqref{local:flow_relocate}, \eqref{local:flow_balance}, \eqref{local:engagements}, \eqref{local:y_sd} and \eqref{local:x_dkj}}
\end{alignat*}

Even though constraints involving location decision variables $z_s$ are relaxed, the scheduling problem continues to suffer from an exponential expansion of the solution space. To overcome this challenge, we propose a trade-off between the optimality of the solution and the computational runtime by implementing a myopic policy algorithm. The algorithm consists of breaking the problem into a set of coupled subproblems and solving them in a sequential manner such that the solution of a subproblem is used as initial conditions for the immediately subsequent subproblem. For the L2D-ESP, we parameterize the time step $t$ and partition it into $T$ coupled subproblems, where each subproblem solves the L2D-ESP for a time step $t$, therefore, each one of them has a solution space smaller than the original scheduling problem, making it easier to solve. Table~\ref{table: schedule} presents the sets, parameters, and decision variables specific to each subproblem.

\begin{table}[htb]
\caption{Specific sets, parameters, and variables for L2D-ESP subproblem $t$.}
\centering
\begin{tabular}{lll}
\hline
Type & Symbol & Description \\
\hline
Sets  
    &$\mathcal{P}_{dj}(t)$ & Set of laser platforms that generate orbital slot $j$\\
    &&for debris $d$ (index $p$; cardinality $P_{dj}(t)$) \\
Parameters     &$R_{dj}(t,i)$ & Debris remediation capacity reward\\
Decision variables &  $y_{pd}(t)$ & $\begin{cases}
                     1, &\text{if laser platform $p$ engages debris $d$} \\
                    0, &\text{otherwise}
                    \end{cases}$ \\
                    &  $x_{dj}(t, i)$ & $\begin{cases}
                     1, &\text{if debris $d$ locates to orbital slot $j$} \\
                    0, &\text{otherwise}
                    \end{cases}$ \\
\hline
\label{table: schedule}
\end{tabular}
\end{table}

\subsubsection{Decision variables}
To solve the L2D-ESP by leveraging myopic policy, we parameterize time step $t$ present in decision variables~\eqref{local:y_sd} and~\eqref{local:x_dkj}. Additionally, orbital slot index $i \in \mathcal{J}_{td}$ is relaxed from constraints~\eqref{local:x_dkj} since the initial orbital slot for each subproblem is known. The new decision variables for the subproblem are defined as:
\begin{subequations}
    \begin{alignat}{2}
    &y_{pd}(t) = \begin{cases}
        1, & \text{if platform $p$ engages debris $d$} \\
        0, & \text{otherwise}
    \end{cases}\\
    &x_{dj}(t, i) = \begin{cases}
        1, & \text{if debris $d$ relocates to orbital slot $j$} \\
        0, & \text{otherwise}
    \end{cases}
    \end{alignat}
\end{subequations}

\subsubsection{Constraints and objective function}
To define the engagement constraints, we parameterize time step $t$ in constraints~\eqref{local:eng_limit} and~\eqref{local:engagements}. Given that platform orbital slots are not considered in the formulation, we introduce the set of platforms $\mathcal{P}_{dj}(t) \subseteq \mathcal{P}$ with index $p$ and cardinality $P_{dj}(t)$, that generate orbital slot $j \in \mathcal{J}_d(t)$ for debris $d \in \mathcal{D}$. Then, the engagement and linking constraints between variables for the subproblem are:
\begin{subequations}
\begin{alignat}{2}
& \sum_{d\in\mathcal{D}}y_{pd}(t) \le 1, \quad& \forall p\in\mathcal{P}\label{policy:engagement}\\
& \sum_{p\in\mathcal{P}_{dj}(t)} y_{pd}(t) \ge {P}_{dj}(t)x_{dj}(t, i), \quad  &\forall d\in\mathcal{D}, \forall j \in \mathcal{J}_d(t) \label{policy:coupling}
\end{alignat}
\end{subequations}

In addition to constraints~\eqref{policy:coupling}, we need to enforce each debris to relocate to at most one new orbital slot. To achieve this, we introduce constraints~\eqref{policy:flow}.
\begin{equation}
\sum_{j\in\mathcal{J}_d(t)} x_{dj}(t, i) \le 1, \quad \forall d \in\mathcal{D} \label{policy:flow}
\end{equation}

The domains of the decision variables are given as:
\begin{subequations}
\begin{alignat}{2}
& x_{dj}(t,i) \in \{0,1\}, \quad& \forall d\in\mathcal{D}, \forall j\in\mathcal{J}_d(t) \label{policy:x_dj} \\
& y_{pd}(t) \in \{0,1\}, \quad& \forall p\in\mathcal{P}, \forall d\in\mathcal{D}\label{policy:y_pd}
\end{alignat}
\end{subequations}

To account for the debris remediation capacity of the constellation for each subproblem $t$, we parameterize indices $t$ and $i$ from Eq.~\eqref{eq:reward}. The reward for each subproblem is given as:
\begin{equation}
    R_{dj}(t, i) = {C}^0_{d}(t) + {C}_{dj}(t,i) + \alpha\Delta h_{dj}(t,i)  + \beta{M}_{d}(t) 
    \label{eq:reward_l2desp}
\end{equation}
albeit we are solving subproblem $t$, we let ${C}^0_{d}(t)={C}^0_{td}$ and ${C}_{dj}(t,i)={C}_{tdij}$ with ${C}^0_{td}$ and ${C}_{tdij}$ defined as in Eq.~\eqref{eq:reward}. Therefore, objective function~\eqref{policy:obj_t} encodes the debris remediation capacity for subproblem $t$.
\begin{equation}
V(t) = \sum_{d\in\mathcal{D}}\sum_{j\in\mathcal{J}_d(t)}R_{dj}(t,i)x_{dj}(t,i)\label{policy:obj_t}\\
\end{equation}

\subsubsection{Mathematical formulation}
The optimization problem for a single subproblem $t$ is formulated as:
\begin{alignat*}{2}
\max \quad & \text{Objective function~\eqref{policy:obj_t}}\\
\text{s.t.} \quad & \text{Constraints \eqref{policy:engagement}, \eqref{policy:coupling} \eqref{policy:flow}, \eqref{policy:x_dj} and \eqref{policy:y_pd}}
\end{alignat*}

As a result of implementing the myopic policy algorithm, we can obtain the debris remediation capacity over the entire mission horizon by summing up objective function~\eqref{policy:obj_t} for all subproblems, as defined in Eq.~\eqref{policy:obj}.
\begin{equation}
V = \sum_{t\in\mathcal{T}}V(t)\label{policy:obj}
\end{equation}

\section{Case studies} \label{sec:case_studies}
To demonstrate the extension and flexibility of the proposed formulations against different mission environments, we present three case studies involving debris remediation missions for (1) small, (2) large, and (3) mixed debris populations (comprised of both small and large debris) in the presence of 10 valuable assets. For each case study, we determine the optimal constellation configuration of 10 laser platforms based on the given debris field using the MCLP, obtain its constellation configuration reward, and then evaluate its debris remediation capacity using the L2D-ESP.

The results obtained for the optimal constellation of 10 space-based laser platforms are benchmarked against two baseline cases: (1) a single platform case and (2) a 10-platform case in an optimized Walker-Delta constellation. The first benchmark, against a single platform system, helps reveal the extent to which debris remediation capacity improves with additional platforms. The second benchmark is intended to compare the effectiveness of asymmetry in the constellation configuration in the debris remediation capacity. Walker-Delta~\citep{walker1984satellite} is a symmetrical constellation pattern that enforces satellites to be placed in circular orbits, all of them having the same altitude and inclination, making them beneficial for missions that require global coverage. The Walker-Delta pattern is denoted as $P/O/F$, where $P$ stands for the number of satellites (platforms), $O$ for the number of orbital planes, and $F$ for the phasing factor between satellites. To identify the optimized 10-platform Walker-Delta constellation, we performed an enumeration of every possible pattern given 10 platforms, resulting in 18 patterns. Additionally, we randomly selected 20 pairs of semi-major axes and inclinations from the set of platform's orbital slots and generated all possible combinations of Walker-Delta constellation configurations using 10 satellites, resulting in 360 constellation configurations. Using a brute-force algorithm, we obtain from the pool of generated Walker-Delta constellations the best-performing one (\textit{i.e.,} the one that collects the highest constellation configuration reward) and obtain its debris remediation capacity leveraging the L2D-ESP.

\subsection{Adopted space-based laser platform}
We adopt the laser and debris parameters from~\citep{phipps2014adroit} (except for the perigee threshold). The laser system considered is a laser-diode-pumped solid-state oscillator-amplifier with an amplifier medium of Nd:YAG or ND:glass, which operates using \SI{100}{ps} ultraviolet pulses at the 3rd harmonic of neodymium at a wavelength of \SI{355}{nm}. We assume that the debris material composition is aluminum and that the laser platform varies the pulse energy to maintain the desired fluence and $c_m$ values. The orbits of the laser platforms and debris are propagated considering up to the $J_2$ perturbation. The epoch $t_0$ is set to February 26, 2024, at 04:30:51 Coordinated Universal Time (UTC). After an L2D engagement, we assume that the debris periapsis radius remains constant, neglecting the influence of $J_2$ to reduce the computational runtime when ${C}_{tdij}$ is computed. Lastly, to avoid a computationally prohibitive optimization problem, we propose a trade-off between the number of time steps and debris considered, that is, a larger population will be simulated with fewer time steps and vice versa. Table~\ref{table:common_parameters} summarizes the parameters shared across all case studies.

\begin{table}[ht]
\caption{Common parameters to all case studies.}
\centering
\begin{tabular}{lll}
\hline
Parameter & Value &Unit\\
\hline
Primary mirror diameter, $D_{\text{eff}}$ &1.50 &\si{m}\\
Momentum coupling coefficient, $c_\text{m}$ &99& \si{N/MW}\\
Wavelength, $\lambda$   &355 &\si{nm}\\
Pulse length, $\tau$  & 100 &\si{ps}\\
Optimal fluence on debris, $\varphi$ & 8.50&\si{kJ/m^2}\\
Beam quality factor, $B^2$ & 2 &-\\
Perigee threshold, $h^*$  & 100 & \si{km}\\
Number of platforms, $P$ & 10&-\\
\hline
\label{table:common_parameters}
\end{tabular}
\end{table}

\subsection{Case study 1: small debris field}\label{sec:case_study_small}
We consider a seven-day mission time horizon, uniformly discretized with a time step size of \SI{130}{s}, adopted from~\citep{phipps2014adroit}. The time step size is determined by the L2D engagement duration of \SI{10}{s}, with a PRF of \SI{56}{Hz}, plus an additional \SI{120}{s} for the laser equipment's cooling. Assuming an uniform surface density of \SI{1}{kg/m^2}, the resulting per-pulse $\Delta v$ is \SI{0.425}{m/s}.
Table~\ref{table:small_parameters} presents the remaining laser and debris parameters.

\begin{table}[ht]
\caption{Small debris remediation specific parameters.}
\centering
\begin{tabular}{lll}
\hline
Parameter & Value & Unit\\
\hline
Laser range & $[175,325]$ & \si{km}\\
Debris density, $\rho$ & 1 & \si{kg/m^2}\\
Time step size & 130 & \si{s}\\
Length L2D engagement & 10 & \si{s}\\
PRF & 56 & \si{Hz}\\
Pulse efficiency, $\eta$ & 0.50 & -\\
\hline
\label{table:small_parameters}
\end{tabular}
\end{table}

\begin{figure}[htbp]
    \centering
    \begin{subfigure}[b]{0.47\textwidth}
        \centering
        \includegraphics[width=0.6\linewidth]{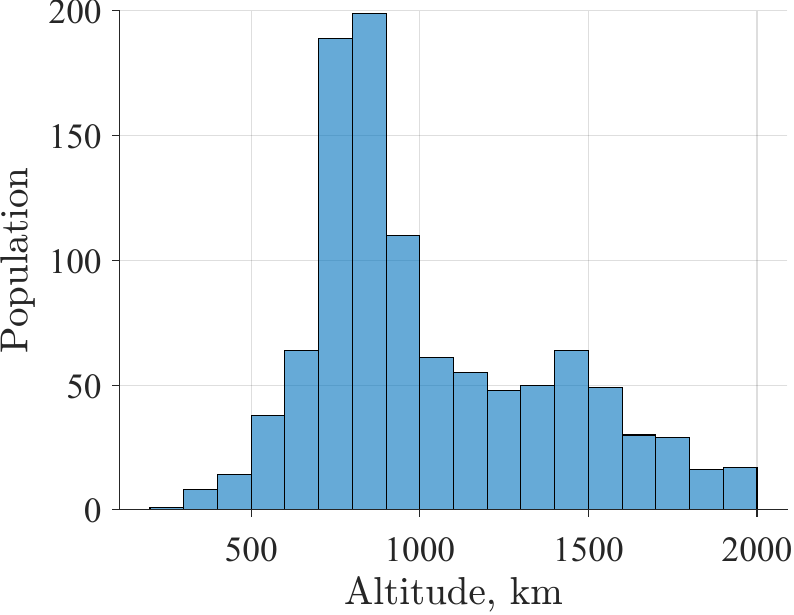}
        \caption{Small debris population.}
        \label{fig:sm_deb_alt_hist}
    \end{subfigure}
    \begin{subfigure}[b]{0.47\textwidth}
        \centering
        \includegraphics[width=0.6\linewidth]{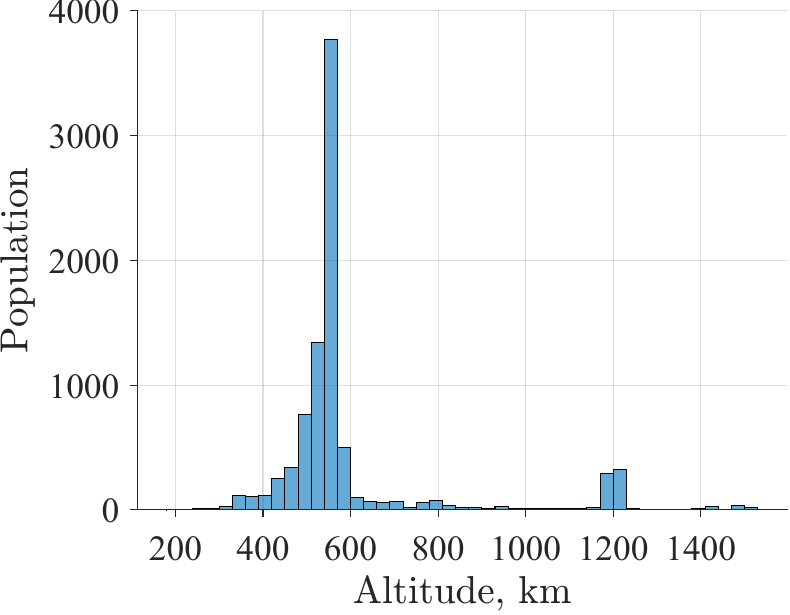}
        \caption{Active satellite population.}
        \label{fig:small_field_n}
     \end{subfigure}
    \caption{Small debris and active satellites population's distribution per altitude bins.}
    \label{fig:debris_vs_sats}
\end{figure}

The small debris field is initialized based on results from a simulation conducted using ESA's Meteoroid And Space debris Terrestrial Environment Reference (MASTER-8) model~\citep{MASTER8}. The simulation considers small debris over an altitude range of \SI{186}{km} to \SI{2000}{km}, which accounts for explosion/collision fragments, sodium-potassium droplets, solid rocket motor slag, and multi-layer insulation. The MASTER-8 outputs are provided as relative frequencies across \num{100} altitude bins that discretize the LEO altitude range of \SI{100}{km} to \SI{2000}{km} into bins having a width of \SI{18.14}{km}. To initialize a representative small debris field across the LEO altitude band, \num{820} samples are drawn from the respective relative frequencies for each altitude bin. Figure~\ref{fig:sm_deb_alt_hist} shows the small debris population per altitude bin. Each small debris object is assumed to be in a circular orbit with uniformly distributed values for the argument of periapsis, the argument of latitude, and the right ascension of the ascending node (RAAN), bounded between 0 and \SI{360}{deg}, and inclination ranging from 0 to \SI{180}{deg}.

Next, we generate a set of orbital slots that laser platforms can take by defining an altitude range from \SI{400}{km} to \SI{1100}{km}, equally spaced in nine altitude layers. We adopt this range as it embraces the peak of operating satellites, outlined in Fig.~\ref{fig:small_field_n}, and the small debris population. Further, we define 10 steps for the argument of latitude and RAAN, uniformly spaced between \SI{0}{deg} and \SI{360}{deg}. In addition, we define 9 uniformly spaced steps for the inclination between \SI{35}{deg} and \SI{90}{deg}. All orbits are assumed to have an eccentricity of zero at the epoch.

For this specific case study, we are not considering valuable assets operating in the environment; consequently, ${C}^0_{td} = {C}_{tdij} = 0$ for all $t \in \mathcal{T}, d \in \mathcal{D}, i \in \mathcal{J}_{td}, j \in \mathcal{J}_{t+1,d}$. Second, ${M}_{td}$ is assumed to be equal for all $d \in \mathcal{D}$ due to the assumption that all debris has equal density.

\subsubsection{Results and discussions}\label{sec:results_small}
\begin{figure}[htbp]
    \centering
     \begin{subfigure}[b]{0.48\textwidth}
        \centering
        \includegraphics[width=0.8\linewidth]{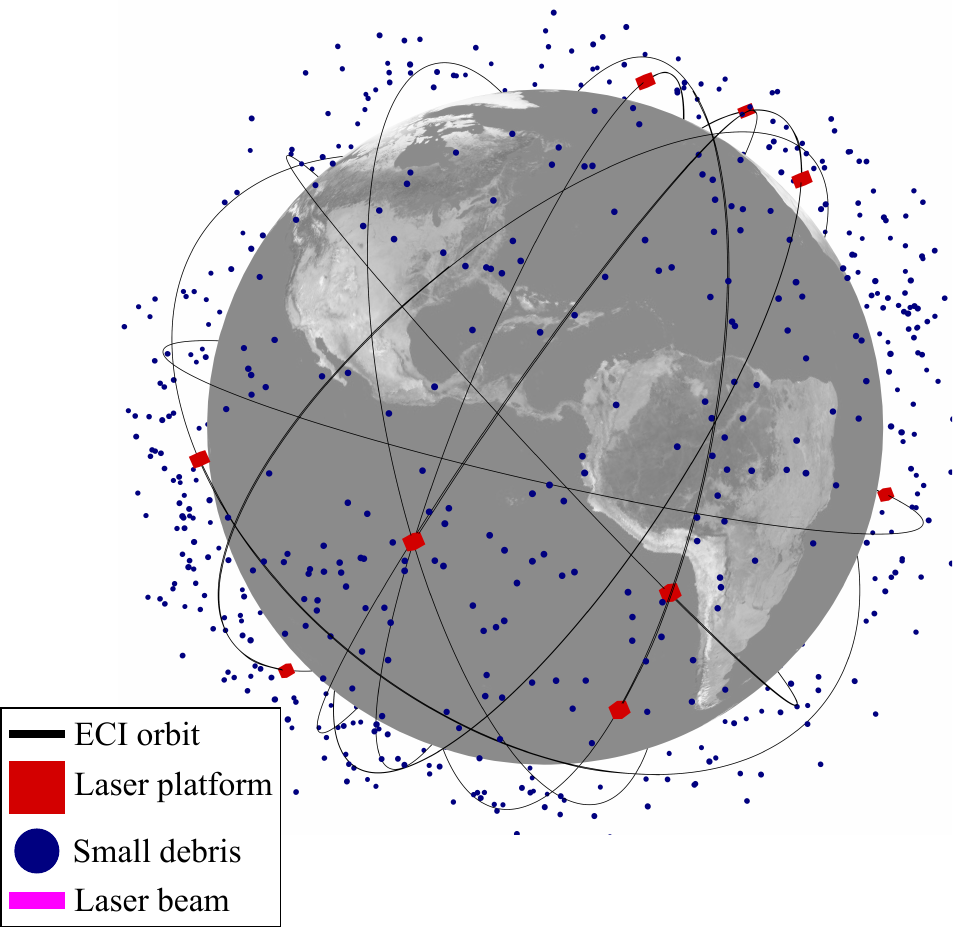}
        \caption{Small debris field and laser platforms at $t=t_0$.}
        \label{fig:small_field}
     \end{subfigure}
     \begin{subfigure}[b]{0.48\textwidth}
        \centering
        \includegraphics[width=0.8\linewidth]{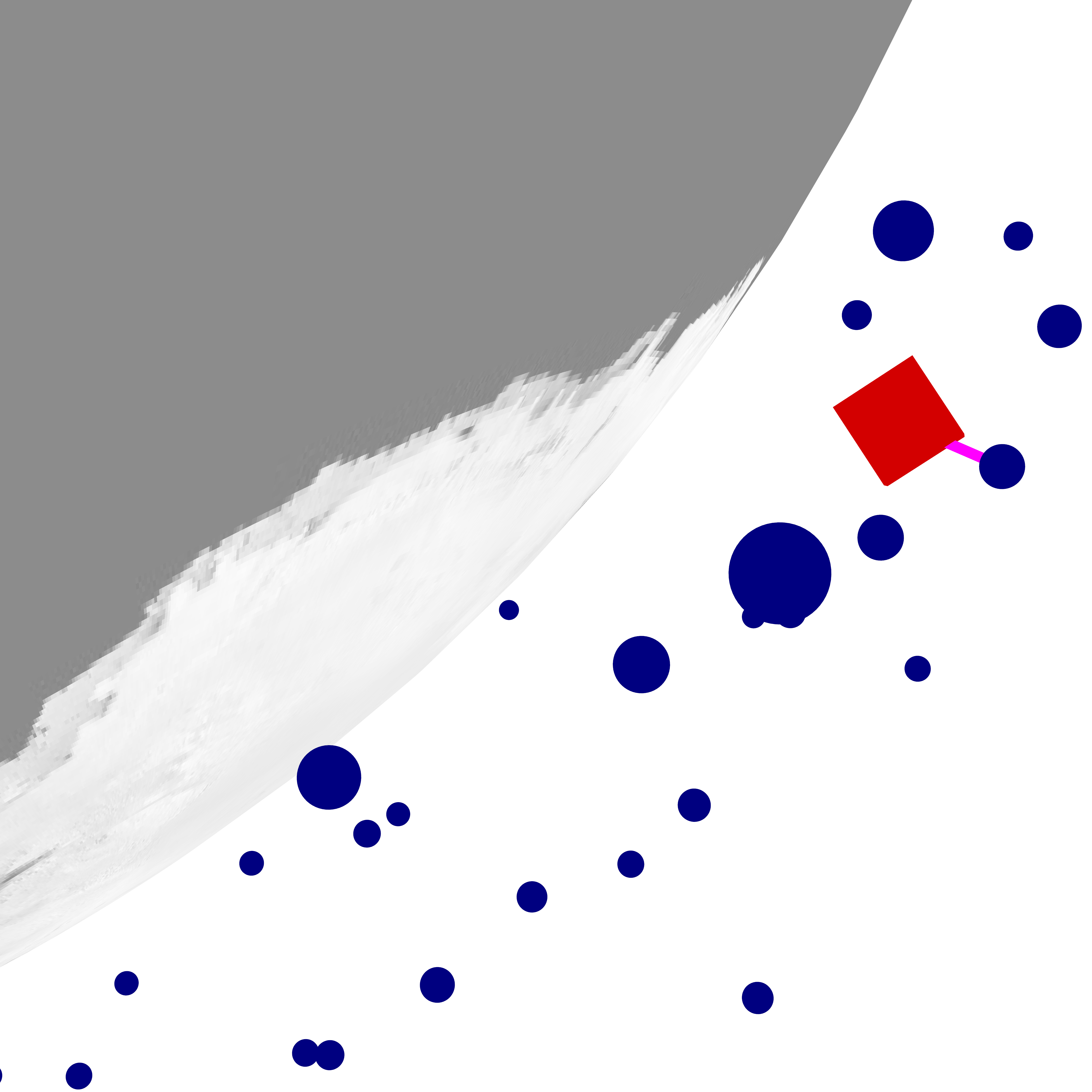}
        \caption{Detailed L2D engagement at $t=11$.}
        \label{fig:beam_small_field}
    \end{subfigure}
    \caption{Optimal 10-platform constellation for small debris remediation mission.}
\end{figure}

The constellation configuration reward obtained is $\pi^\ast= \num{9269}$, and the retrieved optimal constellation configuration consists of 10 platforms asymmetrically distributed. The debris remediation capacity reward is $V^\ast=\num{3132.34}$. Figure~\ref{fig:small_field} provides a visualization of the debris field with 10 laser platforms and their orbits in the Earth-centered inertial (ECI) frame at the epoch $t = t_0$, and Fig.~\ref{fig:beam_small_field} visualizes an L2D engagement over the Antarctic region at time step $t = 11$.

From the L2D optimization results, it is possible to derive certain metrics that are useful for characterizing the mission performance, even though they are not the objective of the optimization. First, the constellation engages \num{599} small debris throughout the seven-day mission time horizon, representing \SI{73.04}{\%} of the total population considered. Second, the number of debris deorbited during this time horizon is \num{422}, or \SI{51.46}{\%} of the total population. Third, we understand debris nudging as the sum of the differences between the periapsis radius at the epoch and the end of the simulation for all engaged but not deorbited debris, yielding a value for the optimal 10-platform constellation of \SI{64707.56}{km}.

The single platform case collects a constellation configuration reward of $\pi^\ast = 1038$ and has a debris remediation capacity of $V^\ast = 722.08$. The number of engaged debris is \num{313}, deorbiting a total of 29 and nudging it \SI{95321.31}{km}. The adoption of 10 platforms over one increases by \SI{88.80}{\%} and \SI{76.94}{\%} the constellation configuration reward $\pi$ and the debris remediation capacity $V$, respectively. Considering the derived metrics, increasing the number of platforms leads to \SI{93.12}{\%} more debris deorbited. The single platform mission nudges debris \SI{95321.31}{km} outperforming the 10-laser platform by \SI{43.31}{\%}. Given that the former is less successful in deorbiting objects, it has more debris contributing to this metric.

The best-performing Walker-Delta constellation has a pattern 10/1/0, an altitude of \SI{7303.14}{km} and inclination of \SI{48.75}{deg}. The constellation configuration reward obtained is $\pi= 7354$, and its debris remediation capacity is $V=2927.24$. The derived metrics for this constellation exhibit that 538 debris are engaged, 382 are deorbited, and it nudges debris a total of \SI{63801.59}{km}.

The debris remediation capacity of the constellation depends not only on the number of platforms but also on their distribution. Breaking the symmetry in the distribution of lasers leads to an increase of \SI{20.66}{\%} and \SI{6.54}{\%} in the constellation configuration reward and the debris remediation capacity, respectively. Additionally, the MCLP-based constellation outperforms the symmetrical Walker-Delta constellation on the derived metrics as it deorbits \SI{9.47}{\%} more debris, and presents an increase of \SI{1.40}{\%} in debris nudged.

\begin{figure}[htbp]
    \centering
    \includegraphics[width=\linewidth]{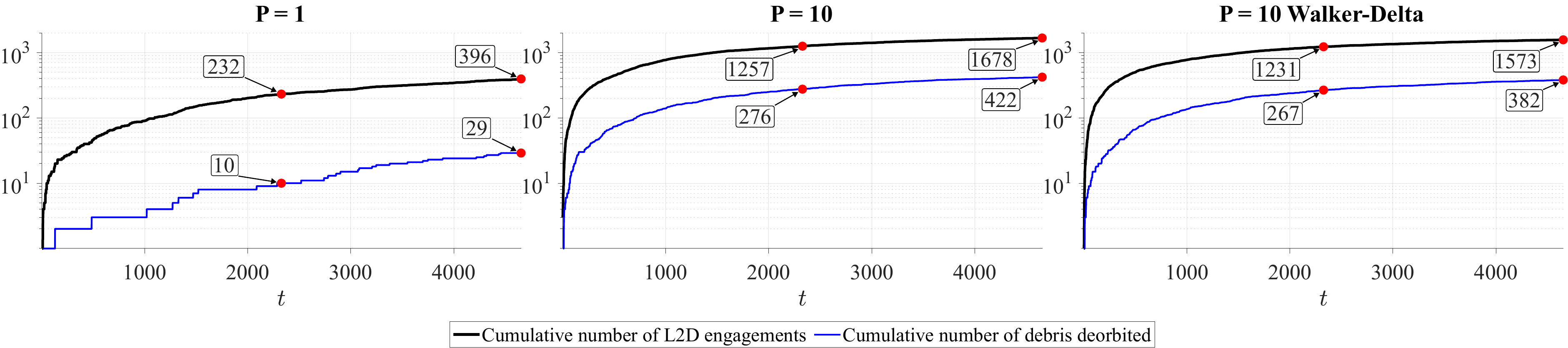}
    \caption{Cumulative number of L2D engagements and small debris deorbited for the three configurations. The highlighted points correspond to the cumulative number of L2D engagements and debris deorbited at the middle and end of the mission planning horizon.}
    \label{fig:periapsis_small}
\end{figure}

Figure~\ref{fig:periapsis_small} presents the cumulative number of L2D engagements and debris deorbited for the single platform and the two constellation configurations. For all three configurations, the number of debris engagements rapidly increases during the first time steps and then stagnates for the remainder of the mission planning horizon. Similarly, for both configurations adopting 10 platforms, the cumulative number of debris deorbited follows the number of L2D engagements with a negative offset, indicating that not all L2D engagements result in an instantaneous deorbit. Conversely, the single platform configuration presents a slower, but steady, increment in the number of debris deorbited, with a negative offset between the two lines.

\begin{figure}[htbp]
    \centering
    \begin{subfigure}[b]{0.32\textwidth}
        \centering
        \includegraphics[width=\linewidth]{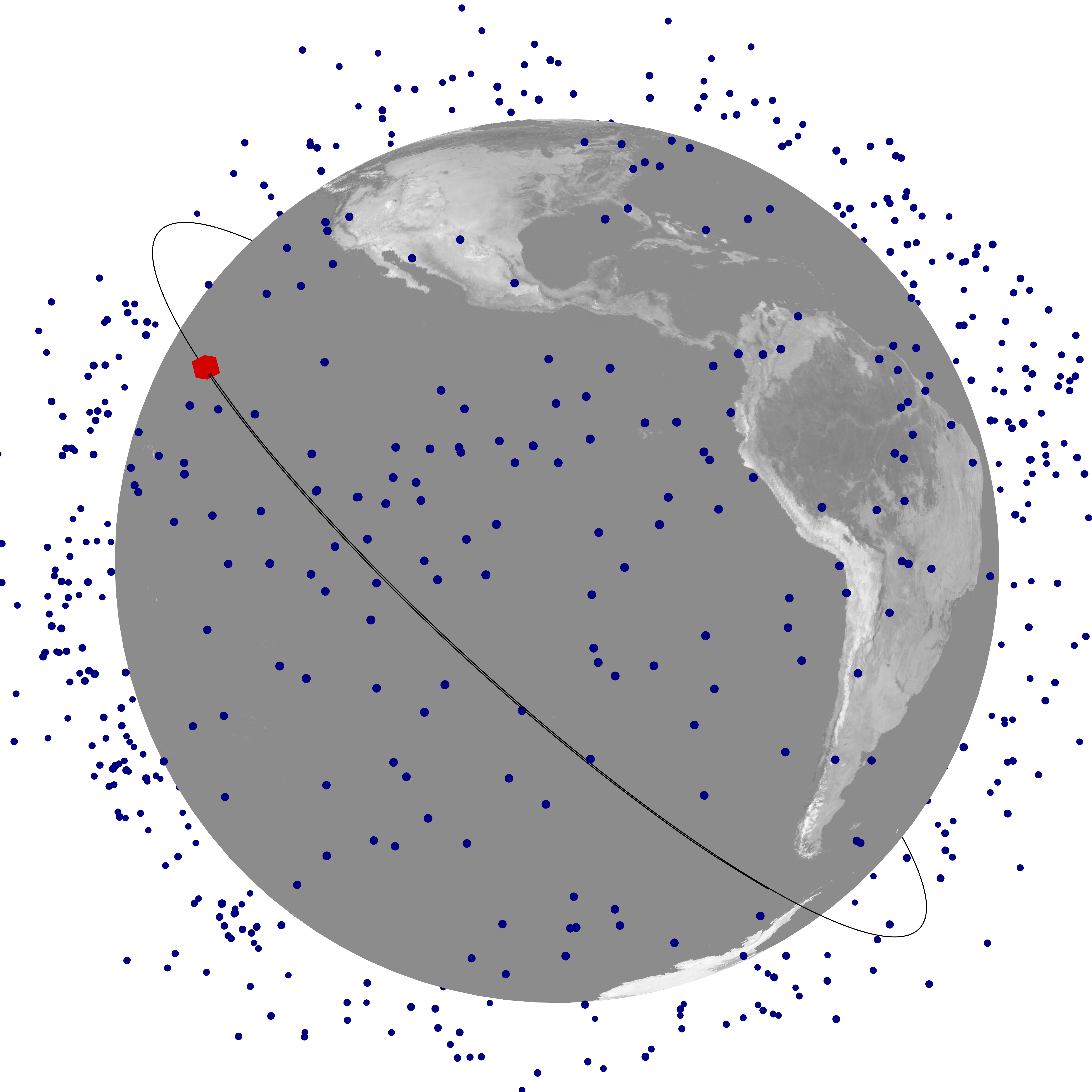}
        \caption{$P=1$.}
     \end{subfigure}
     \begin{subfigure}[b]{0.32\textwidth}
        \centering
        \includegraphics[width=\linewidth]{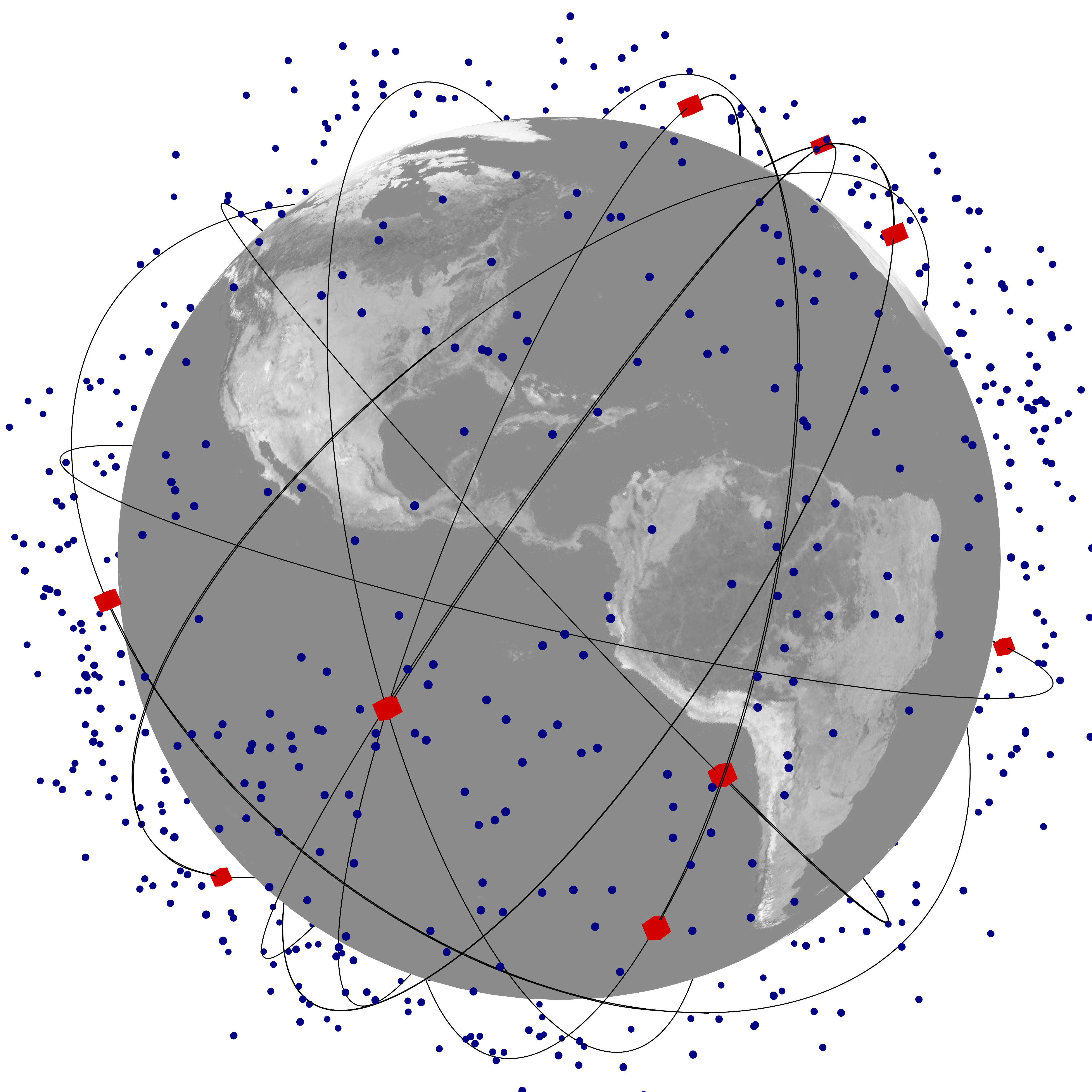}
        \caption{$P=10$.}
    \end{subfigure}
    \begin{subfigure}[b]{0.32\textwidth}
        \centering
        \includegraphics[width=\linewidth]{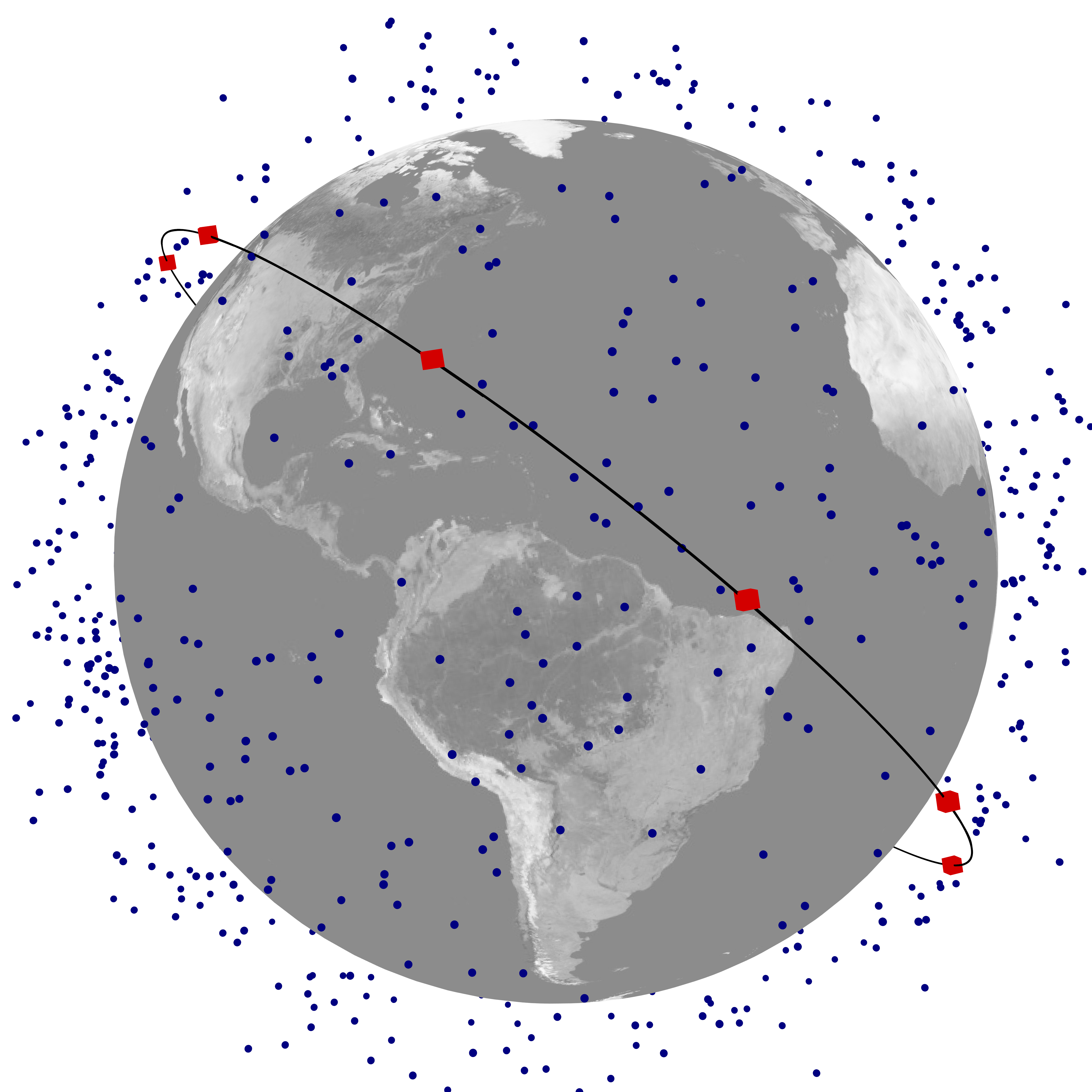}
        \caption{$P=10$ Walker-Delta.}
    \end{subfigure}
    \caption{Small debris field case study: 3D visualization of debris remediation constellations with orbits in ECI at $t=t_0$.}
    \label{fig:3d_small}
\end{figure}

Figure~\ref{fig:3d_small} outlines the small debris field, the two optimal constellation configurations for one and 10 space-based platforms, and the best-performing Walker-Delta constellation. Table~\ref{table:orbelements_small} presents the orbital elements of the three constellation configurations.

\begin{landscape}
    \begin{table}[htb]
\caption{Orbital elements of platforms, defined at epoch $t_0$, for the small debris field case study.}
\centering
\begin{tabular}{lrrrrr}
\hline
Constellation & Sat. index & SMA, km & Incl., deg. & RAAN, deg & Arg. of latitude, deg.\\
\hline
Single platform & $p_1$ & 7390.64 & 62.50 & 280 & 160 \\
10 platform & $p_1$ & 7040.64 & 48.75 & 120 & 80 \\ 
  & $p_2$ & 7040.64 & 69.37 & 320 & 200 \\ 
   &  $p_3$ & 7128.14 & 35 & 80 & 320 \\ 
  & $p_4$ & 7128.14 & 35 & 120 & 360 \\ 
  & $p_5$ & 7128.14 & 83.12 & 120 & 120 \\ 
 & $p_6$ & 7215.64 & 55.62 & 160 & 320 \\ 
  & $p_7$ & 7303.14 & 35 & 160 & 40 \\ 
   & $p_8$ & 7303.14 & 35 & 160 & 80 \\ 
  & $p_9$ & 7303.14 & 35 & 320 & 280 \\ 
   & $p_{10}$ & 7390.64 & 62.50 & 280 & 160 \\ 
10 platform Walker-Delta&$p_1$ & 7303.14 & 48.75 & 0 & 0 \\ 
 &$p_2$ & 7303.14 & 48.75 & 0 & 36 \\ 
 &$p_3$ & 7303.14 & 48.75 & 0 & 72 \\ 
&$p_4$ & 7303.14 & 48.75 & 0 & 108 \\ 
 &$p_5$ & 7303.14 & 48.75 & 0 & 144 \\ 
&$p_6$ & 7303.14 & 48.75 & 0 & 180 \\ 
 &$p_7$ & 7303.14 & 48.75 & 0 & 216 \\ 
&$p_8$ & 7303.14 & 48.75 & 0 & 252 \\ 
&$p_9$ & 7303.14 & 48.75 & 0 & 288 \\
 &$p_{10}$ & 7303.14 & 48.75 & 0 & 324 \\ 
\hline
\label{table:orbelements_small}
\end{tabular}
\end{table} 
\end{landscape}

\subsection{Case study 2: large debris field}\label{sec:study_large}
The mission time horizon is set to 31 days, uniformly discretized with a time step of \SI{160}{s}. Each L2D engagement lasts \SI{40}{s} followed by a cooling period of \SI{120}{s}. For the debris population, we adopt the 50 statistically most concerning debris objects as identified by~\citep{mcknight2021}. For simplicity, we assume that all debris objects have the same cross-sectional area while retaining their actual masses as reported in~\citep{mcknight2021}. With this approach, the magnitude of the delivered $\Delta v$ is not constant among the considered debris field since it is inversely proportional to debris mass. The maximum per-pulse $\Delta v$ is \SI{0.0011}{m/s} corresponding to the minimum mass debris object of \SI{800}{kg}. Conversely, the minimum per-pulse $\Delta v$ is \SI{0.94 e-5}{m/s} corresponding to the maximum mass debris object of \SI{9000}{kg}. To define the platform orbital slots, we maintain the values defined for Sec.~\ref{sec:case_study_small} except for the altitude range; given that the L2D range is larger, we extend the altitude upper bound to \SI{1400}{km}. The specific parameters used in this case study are reported in Table~\ref{table:big_parameters}.

\begin{table}[ht]
\caption{Large debris remediation specific parameters.}
\centering
\begin{tabular}{lll}
\hline
Parameter & Value & Unit\\
\hline
Laser range & $[300, 900]$ &  \si{km} \\
Time step size& 160 & \si{s}\\
Length L2D engagement & 40 & \si{s}\\
PRF  & 21 & \si{Hz}\\
Pulse efficiency, $\eta$ & 1 & -\\
\hline
\label{table:big_parameters}
\end{tabular}
\end{table}

\subsubsection{Results and discussions}
\begin{figure}[htbp]
    \centering
    \begin{subfigure}[b]{0.48\textwidth}
        \centering
        \includegraphics[width=0.8\linewidth]{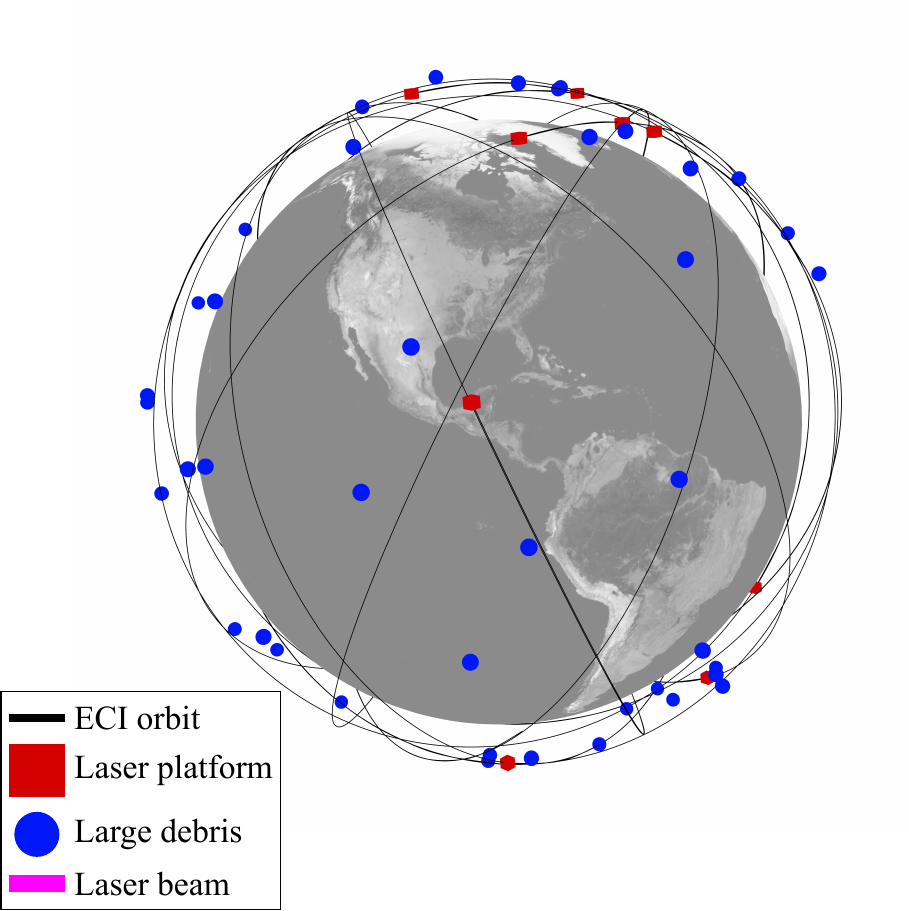}
        \caption{Large debris field and laser platforms at $t=t_0$.}
        \label{fig:large_field}
     \end{subfigure}
     \begin{subfigure}[b]{0.48\textwidth}
        \centering
        \includegraphics[width=0.8\linewidth]{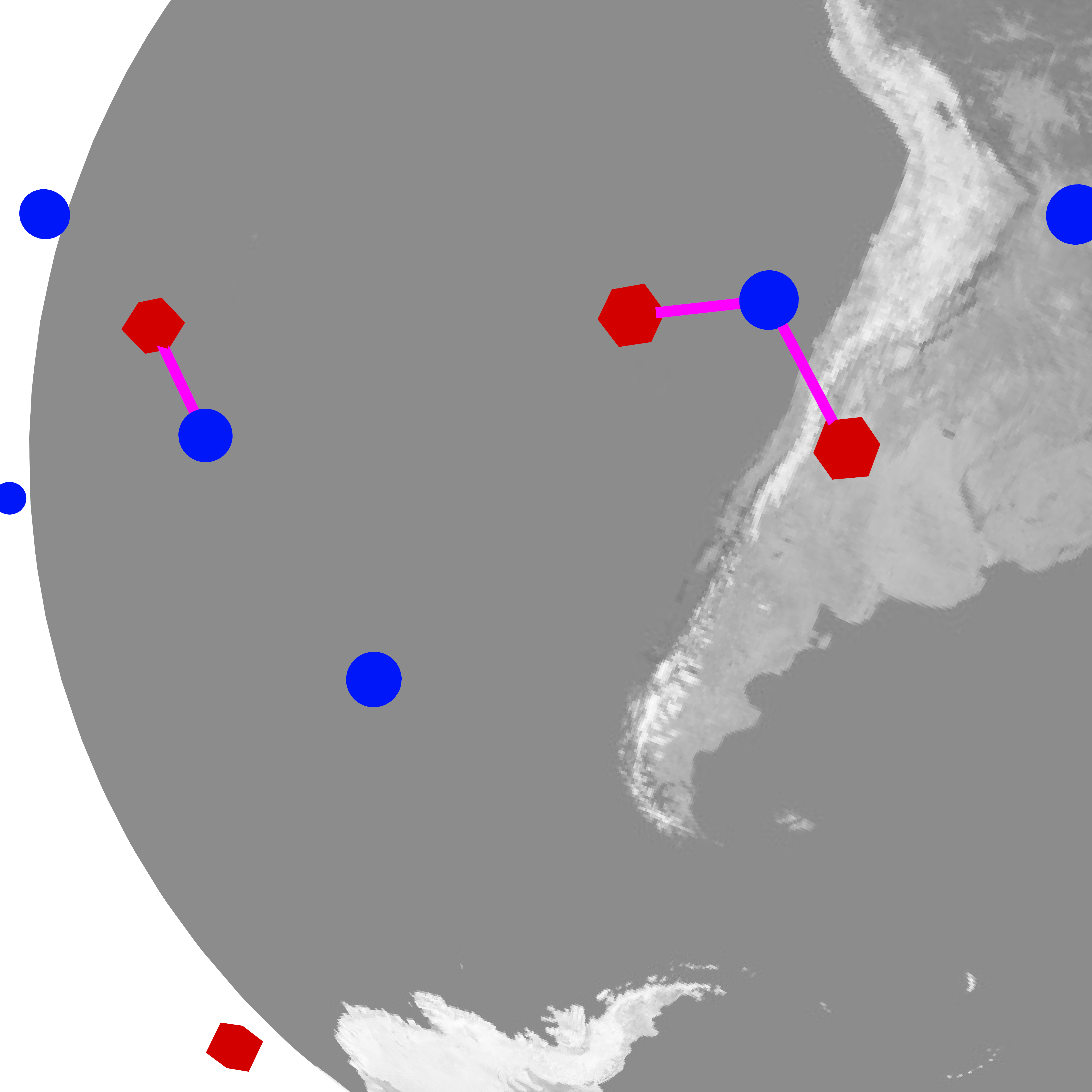}
        \caption{Detailed L2D engagement at $t=1314$.}
        \label{fig:large_field_zoomed}
    \end{subfigure}
    \caption{10-platform optimal constellation for large debris only.}
\end{figure}

The constellation configuration reward obtained is $\pi^\ast=\num{20330.23}$ with an asymmetrical platform distribution, and its debris remediation capacity is $V^\ast=\num{26826.56}$. Figure~\ref{fig:large_field} displays the large debris field and the 10 laser platforms at the epoch $t=t_0$, whereas Fig.~\ref{fig:large_field_zoomed} shows two L2D engagements at time step $t=1314$. The first one consists of a single L2D engagement from one platform to debris; the second one is a collaborative L2D engagement from two platforms to one debris, highlighting the significance of DVA and how the constellation leverages this framework to maximize rewards.

Similarly to the analysis carried out in Sec.~\ref{sec:results_small}, we derive metrics from the obtained results. The total number of L2D engagements during the operation of the mission is \num{20225}, and the mass of the most engaged debris is \SI{9000}{kg}. Given debris's large mass, it is less likely to induce significant periapsis radius reductions due to L2D engagements at each time step; consequently, the term that accounts for debris mass in Eq.~\eqref{eq:reward_l2desp} weighs the most in the reward. Further, the constellation nudges debris \SI{941.54}{km}.

Changing the cardinality of the constellation and imposing only one platform into the MCLP problem formulation retrieves a constellation configuration reward $\pi^\ast=2916.35$ and a debris remediation capacity of $V^\ast=6526.54$. The significant reduction in the number of platforms directly impacts the results by reducing \SI{85.65}{\%} and \SI{75.67}{\%} the constellation configuration reward and the debris remediation capacity, respectively. The single platform performs \num{4163} L2D engagements, with the same most engaged debris object as the 10-platform constellation has. The difference with respect to the 10-platform constellation is \SI{79.41}{\%}. The single platform nudges debris \SI{214.70}{km} signifying a decrement of \SI{77.19}{\%} compared to the 10-platform constellation.

The best-performing Walker-Delta constellation, with a pattern of 10/10/0, an altitude of \SI{7040.64}{km} and an inclination of \SI{76.25}{km}, achieves a constellation configuration reward of $\pi=\num{11255.10}$ and a debris remediation capacity $V=\num{22651.94}$. The Walker-Delta constellation underperforms compared to the asymmetrical MCLP-based constellation by \SI{44.63}{\%} and \SI{15.56}{\%} for the constellation configuration reward and debris remediation capacity, respectively. The total number of L2D engagements for the Walker-Delta constellation is \num{16947}, which is \SI{14.42}{\%} less than that of the optimal MCLP-based one, but it continues to have the most L2D engagements to a \SI{9000}{kg} debris object. Lastly, the Walker-Delta constellation nudges debris \SI{975.45}{km}, compared to the \SI{941.54}{km} obtained by the optimal constellation. This phenomenon can be associated with the fact that neither the MCLP formulation nor the L2D-ESP directly consider nudging as an optimization metric.

\begin{figure}[htbp]
    \centering
    \includegraphics[width=\linewidth]{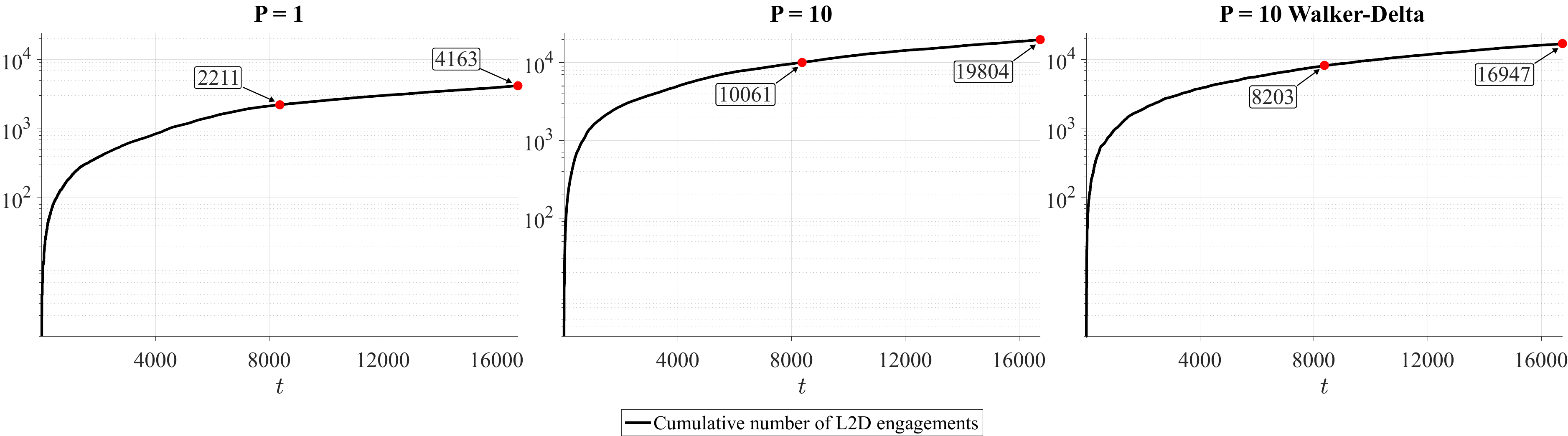}
    \caption{Cumulative number of L2D engagements for the large debris field. The highlighted points correspond to the cumulative number of L2D engagements at the middle and end of the mission planning horizon.}
    \label{fig:periapsis_large}
\end{figure}

\begin{figure}[htbp]
    \centering
    \begin{subfigure}[b]{0.32\textwidth}
        \centering
        \includegraphics[width=\linewidth]{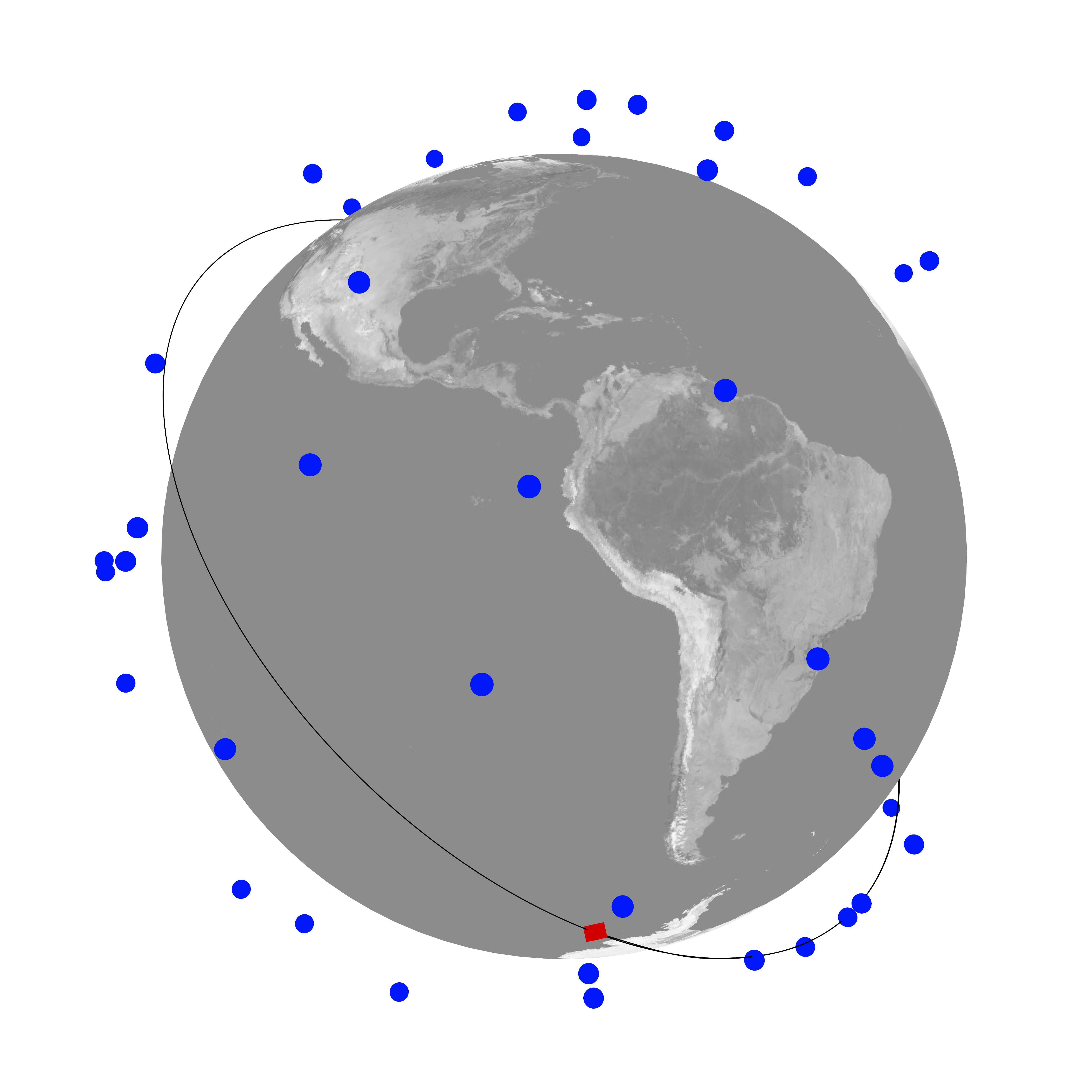}
        \caption{$P=1$.}
     \end{subfigure}
     \begin{subfigure}[b]{0.32\textwidth}
        \centering
        \includegraphics[width=\linewidth]{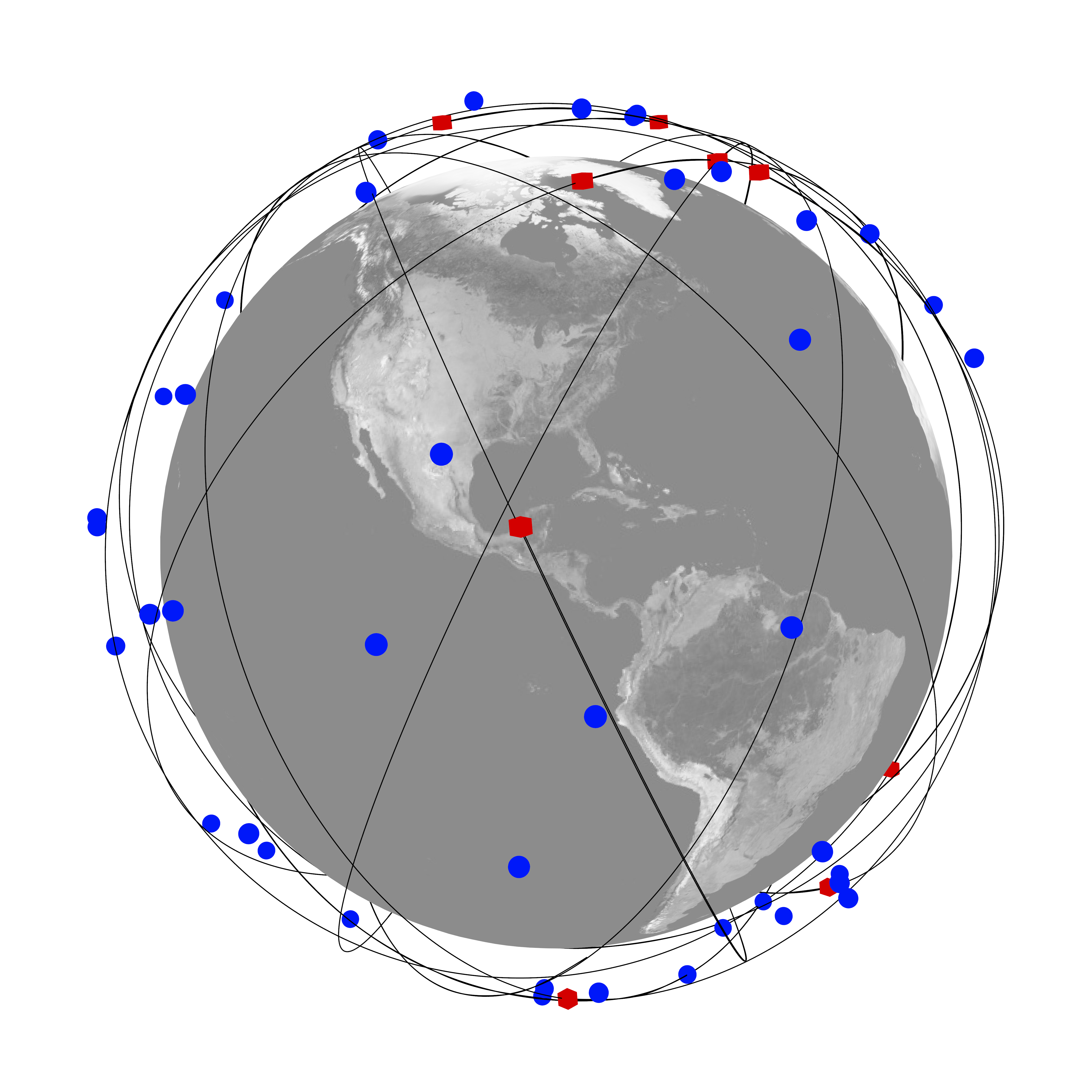}
        \caption{$P=10$.}
    \end{subfigure}
    \begin{subfigure}[b]{0.32\textwidth}
        \centering
        \includegraphics[width=\linewidth]{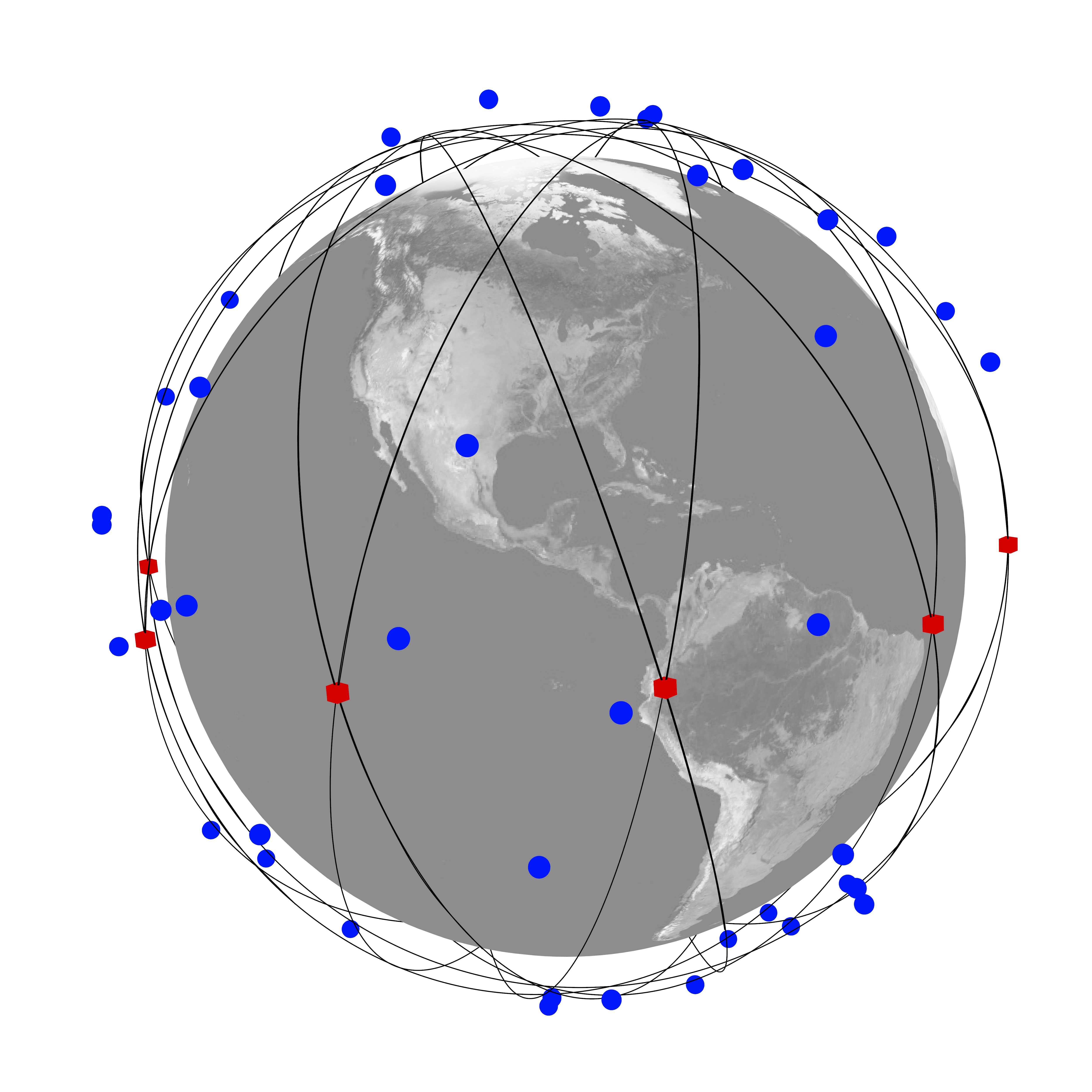}
        \caption{$P=10$ Walker-Delta.}
    \end{subfigure}
    \caption{Large debris field case study: 3D visualization of debris remediation constellations with orbits in ECI at $t=t_0$.}
    \label{fig:3d_large}
\end{figure}

Figure~\ref{fig:periapsis_large} presents the cumulative number of L2D engagements. All three configurations present a rapid increase in the number of debris engaged for the initial part of the mission planning horizon, and then their increase rate tapers off as the time step index increases. Figure~\ref{fig:3d_large} presents the large debris field, the two MCLP-based optimal configurations, and the best-performing Walker-Delta constellation at the epoch. The orbital elements of all configurations are listed in Table~\ref{table:orbelements_large}.

\begin{landscape}
   \begin{table}[htb]
\caption{Orbital elements of platforms, defined at epoch $t_0$, for the large debris field case study.}
\centering
\begin{tabular}{lrrrrr}
\hline
Constellation & Sat. index & SMA, km & Incl., deg. & RAAN, deg & Arg. of latitude, deg. \\
\hline
Single platform & $p_1$ & 7215.64 & 69.375 & 280 & 240\\
10 platform & $p_1$ & 7128.14 & 90 & 40 & 80 \\ 
&$p_2$ & 7215.64 & 62.50 & 80 & 80 \\ 
&$p_3$ & 7215.64 & 62.50 & 120 & 80 \\ 
 &$p_4$ & 7215.64 & 62.50 & 160 & 120 \\ 
&$p_5$ & 7215.64 & 69.375 & 200 & 80 \\ 
&$p_6$ & 7215.64 & 69.375 & 280 & 240 \\ 
&$p_7$ & 7215.64 & 69.375 & 320 & 160 \\ 
&$p_8$ & 7215.64 & 69.375 & 0 & 240 \\ 
 &$p_9$ & 7215.64 & 76.25 & 40 & 120 \\ 
&$p_{10}$ & 7215.64 & 83.125 & 240 & 320 \\ 
10 platform Walker-Delta&$p_1$ & 7040.64 & 76.25 & 0 & 0 \\ 
 &$p_2$ & 7040.64 & 76.25 & 36 & 0 \\ 
&$p_3$ & 7040.64 & 76.25 & 72 & 0 \\ 
 & $p_4$ & 7040.64 & 76.25 & 108 & 0 \\ 
 &$p_5$ & 7040.64 & 76.25 & 144 & 0 \\ 
&$p_6$ & 7040.64 & 76.25 & 180 & 0 \\ 
&$p_7$ & 7040.64 & 76.25 & 216 & 0 \\ 
&$p_8$ & 7040.64 & 76.25 & 252 & 0 \\ 
&$p_9$ & 7040.64 & 76.25 & 288 & 0 \\ 
&$p_{10}$ & 7040.64 & 76.25 & 324 & 0 \\  
\hline
\label{table:orbelements_large}
\end{tabular}
\end{table} 
\end{landscape}

\subsection{Case study 3: mixed debris field with valuable assets}\label{sec:case_study_alldebris}
In this case study, the mission horizon of seven days is uniformly discretized with a time step size of \SI{160}{s} to encompass the L2D engagement and cooling times for small and large debris. Further, we maintain the laser operational ranges for small and large debris as defined in Sec.~\ref{sec:case_study_small} and Sec.~\ref{sec:study_large}, respectively.

The small debris field is generated as in Sec.~\ref{sec:case_study_small}, and we define a population of 820 objects. The large debris field is the same as in Sec.~\ref{sec:study_large}, in addition to defunct satellite COSMOS 2221 (NORAD ID: 22236) with an associated reward of $G_0=\num{e6}$ for the MCLP formulation, and $G_0=\num{e4}$ for the L2D-ESP. Lastly, we define a set of 10 valuable assets $\mathcal{K} = \{k_1, \ldots, k_{10}\}$ and Table~\ref{table:catalog_number} presents their respective orbital elements defined at epoch $t_0$, with $\omega$ and $\nu$ as the argument of periapsis and true anomaly, respectively. The conjunction ellipsoid for all valuable assets is assumed to be a perfect sphere of radius \SI{10}{km}.

\begin{table*}[ht]
\caption{Orbital elements of valuable assets defined at epoch $t_0$.}
\centering
\begin{tabular}{lrrrrrr}
\hline
Index & SMA, km & $e$ & Incl., deg. & RAAN, deg & $\omega$, deg. & $\nu$, deg. \\
\hline
$k_1$ & 6862.80 & 0.0013 & 52.94 & 240.18 & 76.06 & 300.22 \\
$k_2$ & 6933.20 & 0.0011 & 53.20 & 109.69 & 106.38 & 63.51 \\
$k_3$ & 6921.01 & 0.0015 & 52.94 & 224.49 & 71.00 & 181.65 \\
$k_4$ & 6924.25 & 0.0014 & 53.35 & 66.92 & 68.18 & 316.11 \\
$k_5$ & 6874.50 & 0.0005 & 97.53 & 132.48 & 273.28 & 73.11 \\
$k_6$ & 7131.61 & 0.0021 & 86.46 & 16.74 & 77.71 & 170.80 \\
$k_7$ & 6917.09 & 0.0011 & 53.12 & 221.95 & 66.99 & 60.67 \\
$k_8$ & 7565.34 & 0.0010 & 88.04 & 81.73 & 122.65 & 66.20 \\
$k_9$ & 6955.44 & 0.0012 & 69.88 & 276.33 & 67.97 & 307.92 \\
$k_{10}$ & 6796.97 & 0.0006 & 51.70 & 150.06 & 78.94 & 50.82 \\
\hline
\label{table:catalog_number}
\end{tabular}
\end{table*}

In addition to the 10 valuable assets, the mixed debris field study includes a known close conjunction event on February 28, 2024. The conjunction event was between the defunct satellite 22236 and the non-maneuverable operational NASA TIMED satellite (NORAD ID: 26998) with a reported miss distance of \SI{20}{m} at the time of closest approach (TCA)~\citep{COSMOS2221_TIMED_CARA}. In the simulated conjunction event within the mixed debris field, a miss distance of \SI{2.30}{km} is observed at the TCA, corresponding to time step $t=1081$.

\subsubsection{Results and discussions}\label{sec:heterogeneous_results}
\begin{figure}[htbp]
    \centering
    \begin{subfigure}[b]{0.48\textwidth}
        \centering
        \includegraphics[width=0.9\textwidth]{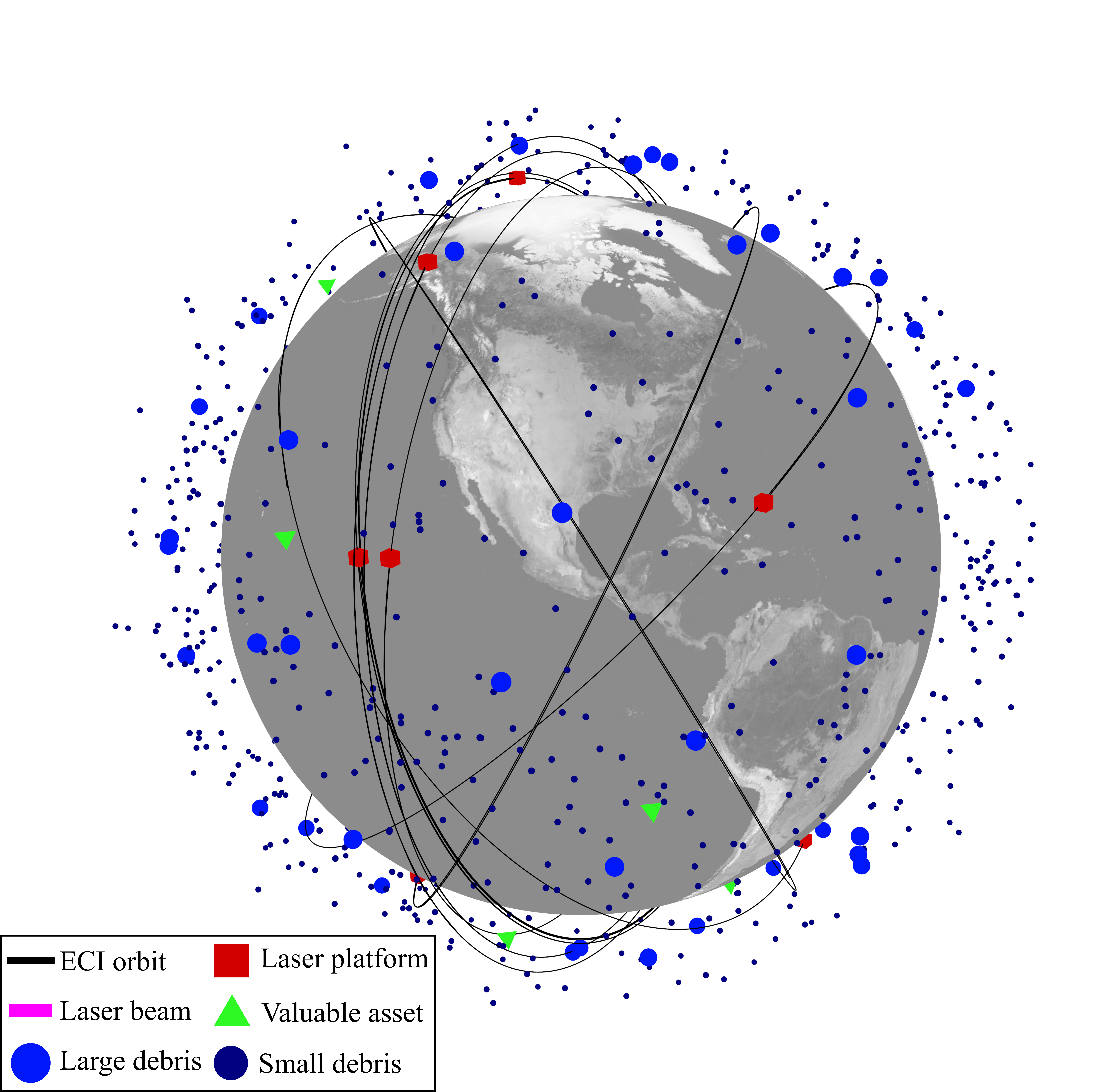}
        \caption{Mixed debris field and laser platforms at $t=t_0$.}
        \label{fig:all_field}
     \end{subfigure}
     \begin{subfigure}[b]{0.48\textwidth}
        \centering
        \includegraphics[width=0.9\textwidth]{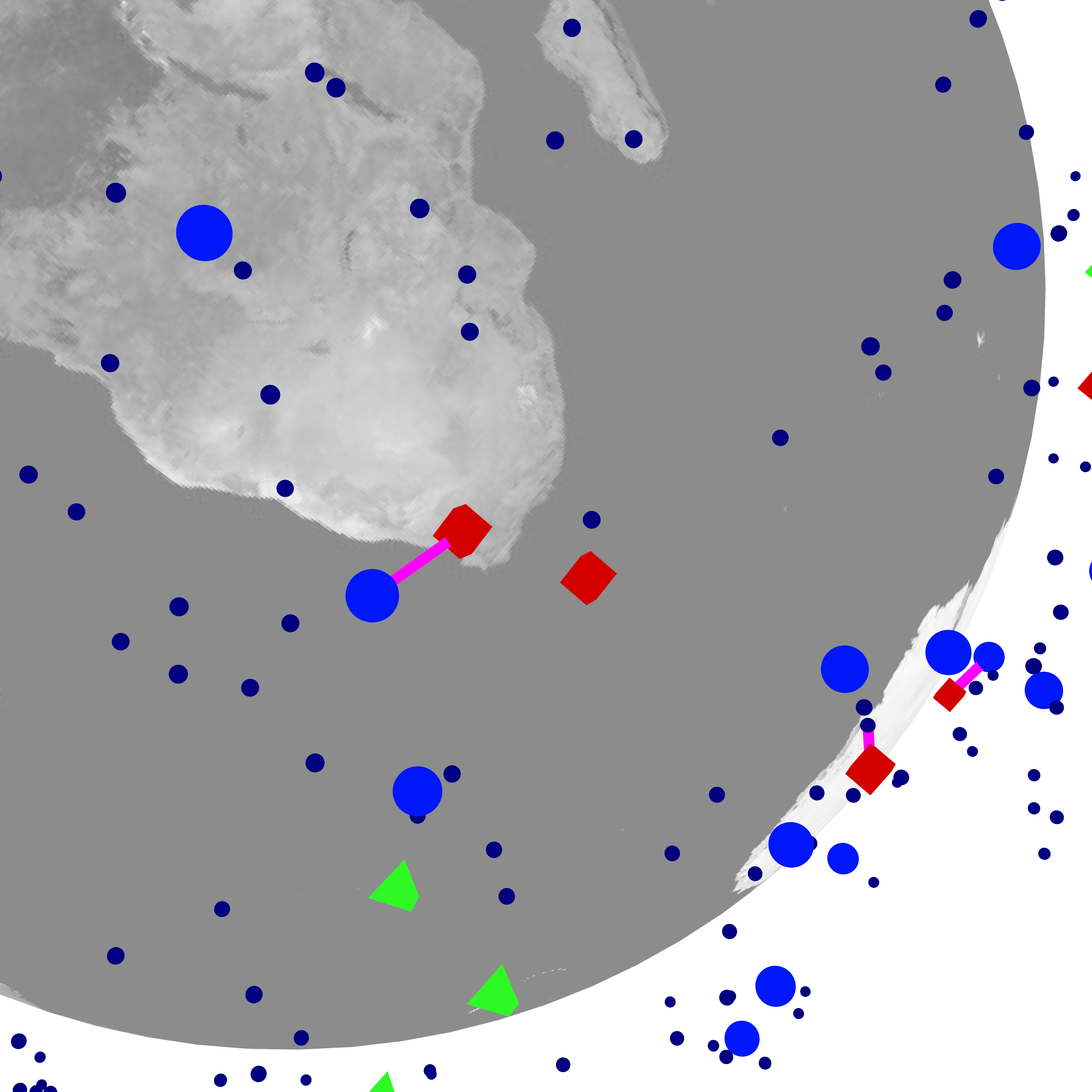}
        \caption{Detailed L2D engagement at $t=7$.}
        \label{fig:all_field_zoomed}
    \end{subfigure}
    \caption{10-platform constellation for small and large debris.}
\end{figure}

The optimal 10-platform constellation achieves a constellation configuration reward of $\pi^\ast=\num{7.10e7}$ and a debris remediation capacity of $V^\ast=\num{2.86e5}$. Figure~\ref{fig:all_field} displays the mixed debris field and the 10 laser platforms at the epoch $t=t_0$. Figure~\ref{fig:all_field_zoomed} is a detailed view that highlights three L2D engagements at time step $t=7$. First, two platforms engage two large debris objects; second, a single platform engages small debris from a closer range.

The optimal 10-platform constellation engages 546 debris objects of various sizes, successfully deorbiting 141 small debris, and nudging debris \SI{43510.04}{km}. The conjunction event expected to occur at time step $t=1081$ in this scenario is successfully avoided, with the miss distance between  NASA satellite 26998 and defunct satellite 22236 increasing from \SI{2.30}{km} to \SI{2104.47}{km} after 17 L2D engagements during the imposed time step window $t\in[500, 622]$. The new minimum miss distance between the two objects is \SI{545.91}{km} at time step $t=103$ when the debris has not been engaged yet.

The baseline single platform case collects a constellation configuration reward of $\pi^\ast=\num{1.60e7}$ and a debris remediation capacity of $V^\ast= \num{69211.24}$. Compared to the optimal 10-platform constellation, these values represent significant reductions of \SI{77.46}{\%} in constellation configuration reward and \SI{75.83}{\%} in debris remediation capacity. The single platform engages with 228 debris, successfully deorbits 24 small debris, and nudges debris \SI{25733.06}{km}. These metrics present reductions of \SI{82.97}{\%} in the number of total debris deorbited and \SI{40.85}{\%} in the total distance nudged compared to the 10-platform constellation. The new closest approach between defunct satellite 22236 and satellite 26998 has a distance of \SI{57.43}{km} at time step $t=936$, and a new miss distance of \SI{411.71}{km} at the original TCA corresponding to $t=1081$.

The best-performing Walker-Delta constellation has a pattern 10/5/2 with an orbital altitude of \SI{6953.14}{km} and an inclination of \SI{76.25}{deg}. It achieves a constellation configuration reward $\pi=\num{2.98e5}$ and a debris remediation capacity of $V=\num{1.23e5}$. When the same number of laser platforms is used but constrained to fit the symmetrical Walker-Delta pattern, there is a substantial reduction in both the constellation configuration reward and debris remediation capacity by \SI{99.58}{\%} and \SI{56.90}{\%}, respectively. The Walker-Delta constellation engages with 441 debris objects from which it successfully deorbits 129 of them, resulting in \SI{8.51}{\%} fewer than the optimal 10-platform constellation. The minimum miss distance between the active and defunct satellite is \SI{545.94}{km} at time step $t=103$. Even though the Walker-Delta does not engage the defunct satellite at the imposed time step window $t\in[500, 622]$, it does so before it. Consequently, the relative distance between the two objects at time step $t=1081$ is \SI{2605.02}{km}.

\begin{figure}[htbp]
    \centering
    \includegraphics[width=\linewidth]{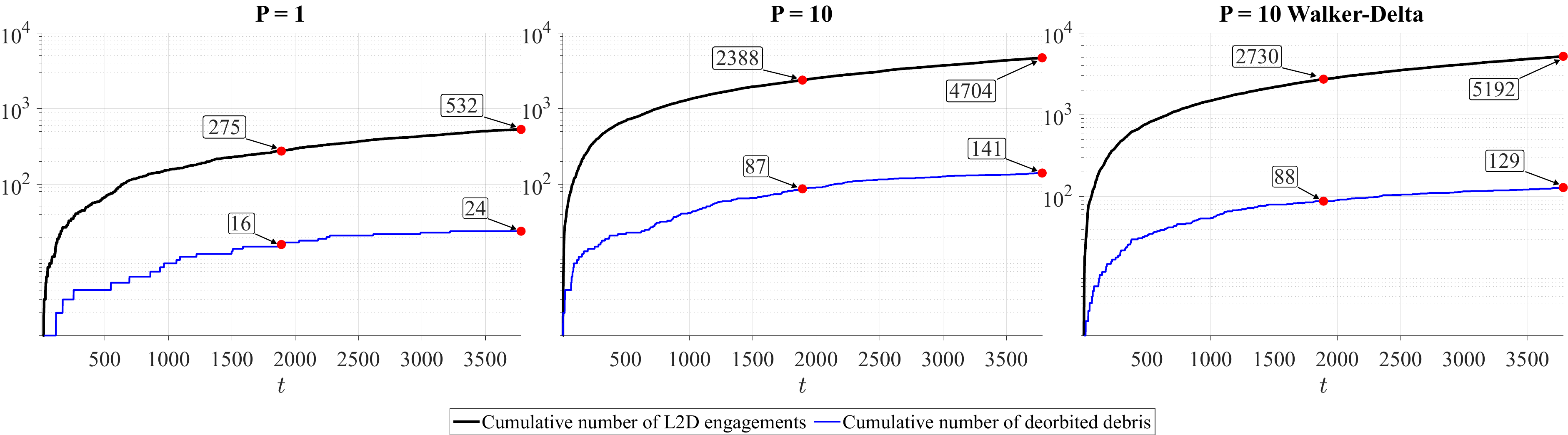}
    \caption{Cumulative number of L2D engagements and debris deorbited for the mixed debris field. The highlighted points correspond to the cumulative number of L2D engagements and debris deorbited at the middle and end of the mission planning horizon.}
    \label{fig:periapsis_all}
\end{figure}

\begin{figure}[htbp]
    \centering
    \begin{subfigure}[b]{0.32\textwidth}
        \centering
        \includegraphics[width=\linewidth]{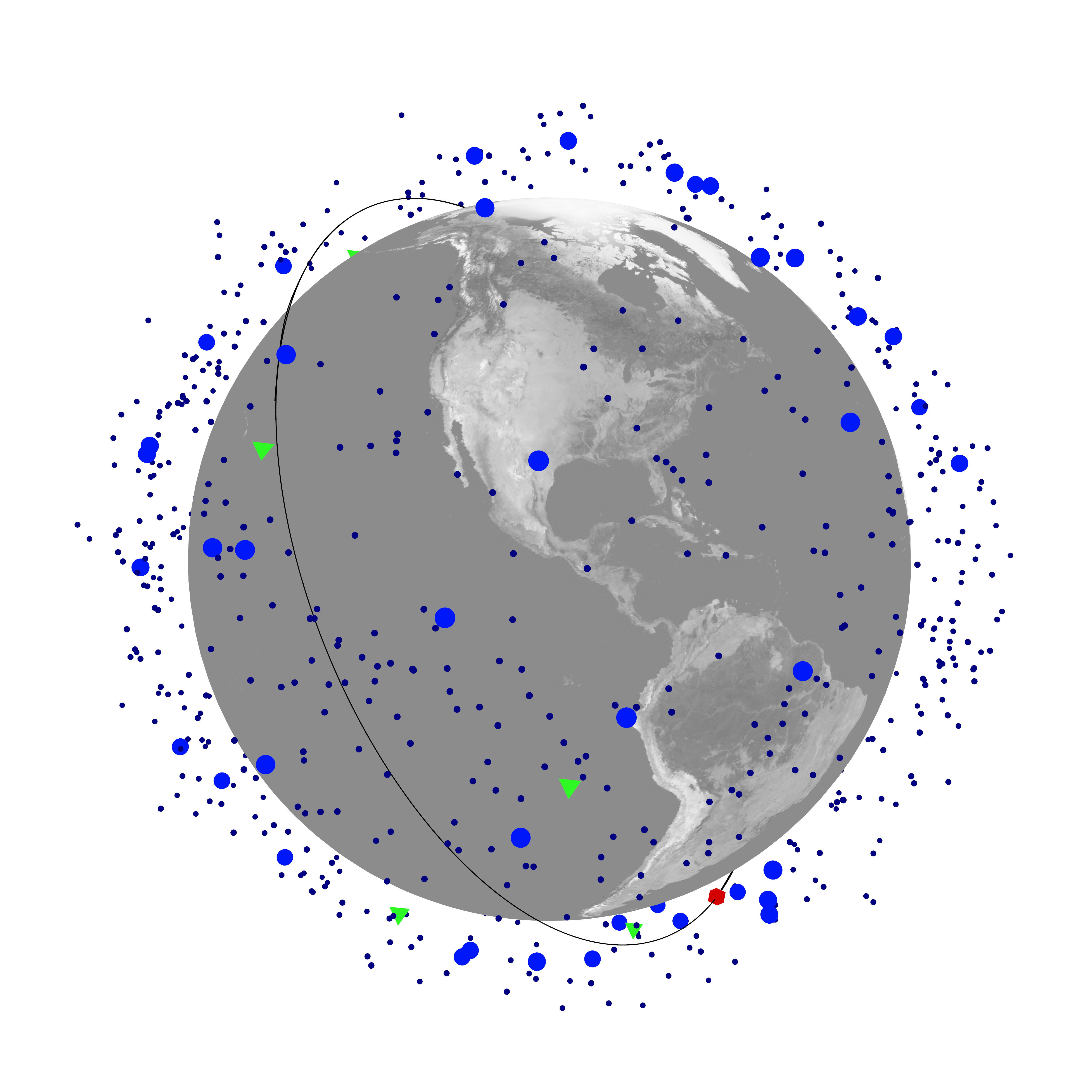}
        \caption{$P=1$.}
     \end{subfigure}
     \begin{subfigure}[b]{0.32\textwidth}
        \centering
        \includegraphics[width=\linewidth]{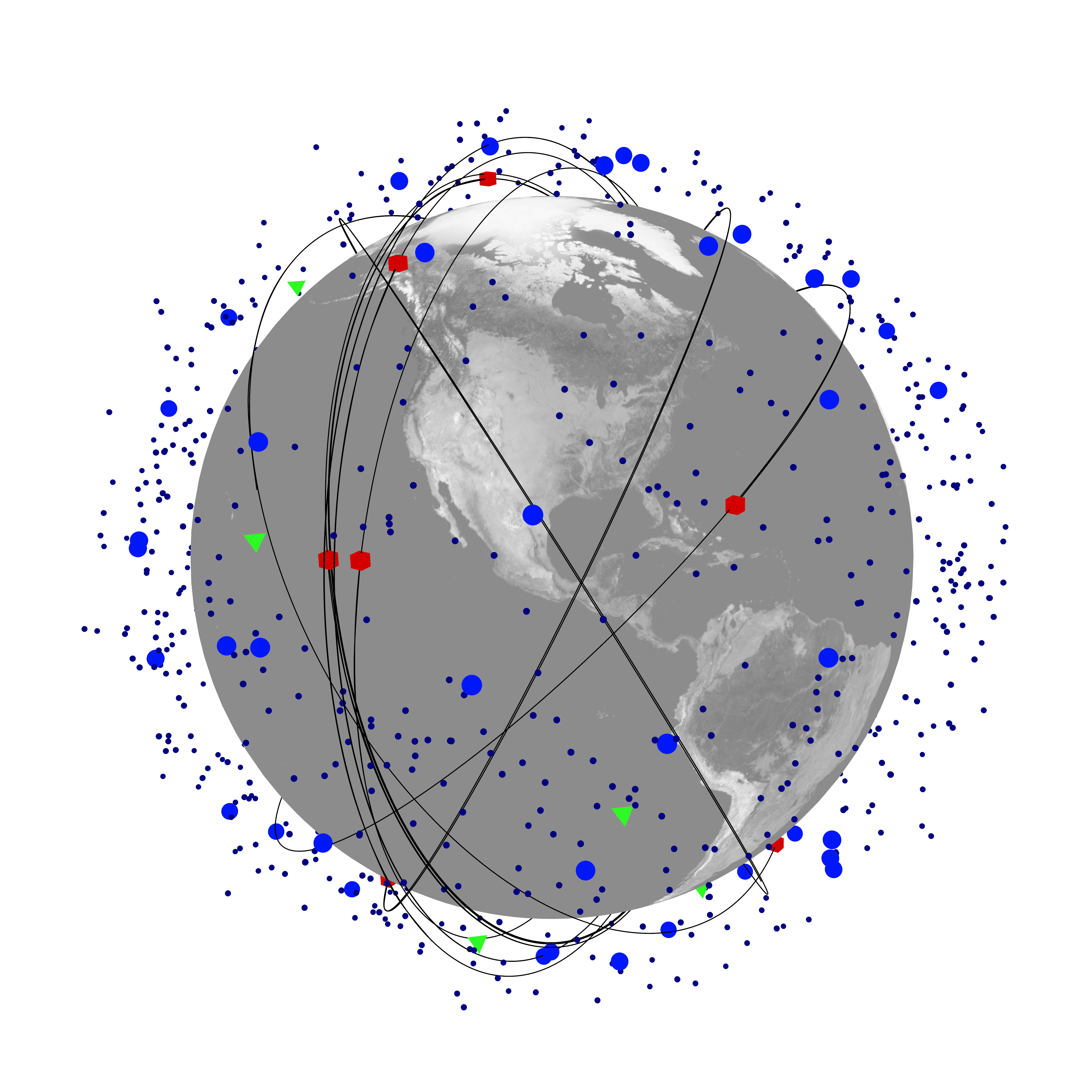}
        \caption{$P=10$.}
    \end{subfigure}
    \begin{subfigure}[b]{0.32\textwidth}
        \centering
        \includegraphics[width=\linewidth]{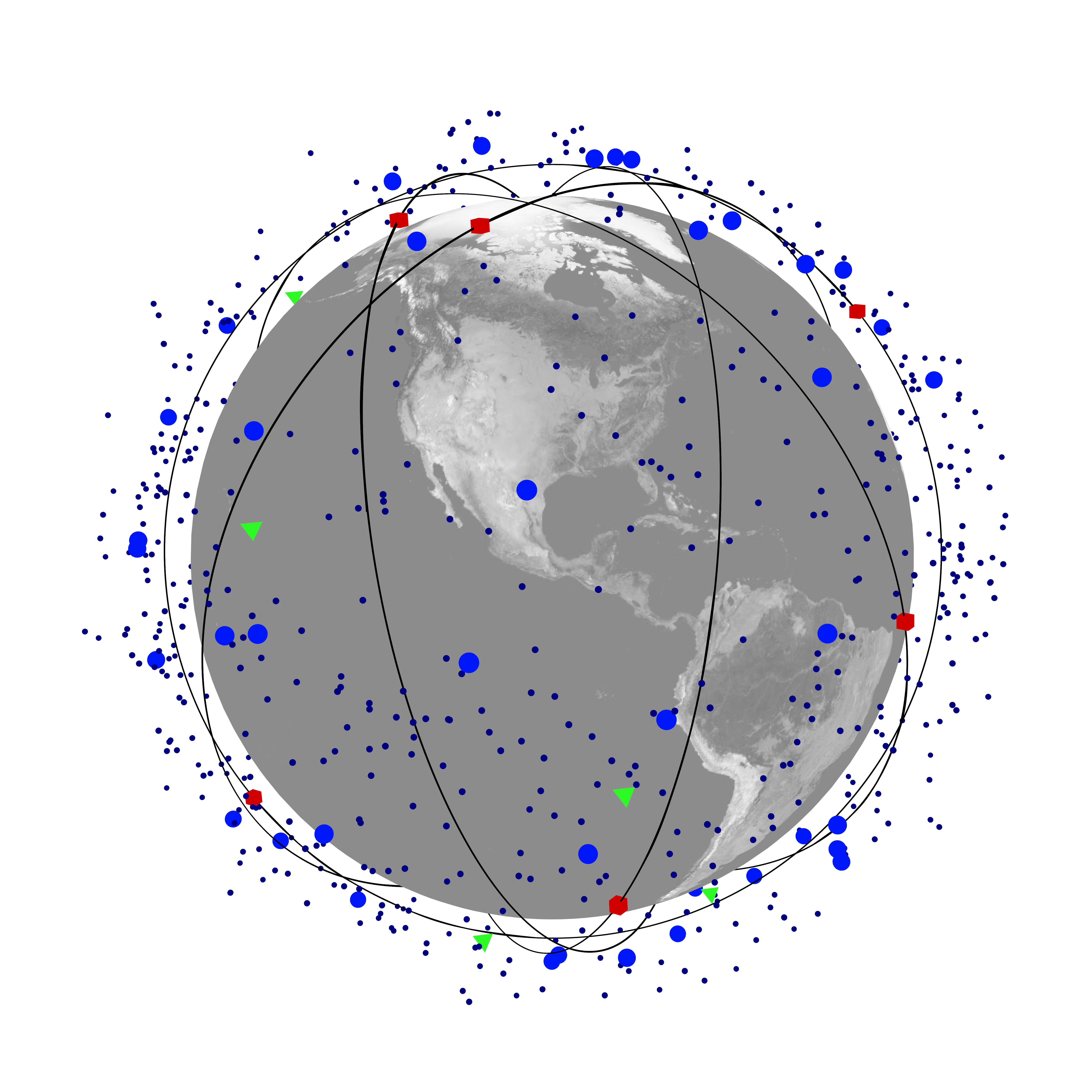}
        \caption{$P=10$ Walker-Delta.}
    \end{subfigure}
    \caption{Mixed debris field case study: 3D visualization of debris remediation constellations with orbits in ECI at $t=t_0$.}
    \label{fig:3d_all}
\end{figure}

Figure~\ref{fig:periapsis_all} illustrates the cumulative number of L2D engagements and debris deorbited for the three configurations. All three of them present a surge in the number of L2D engagements and debris deorbited for the first time steps and a slower increase rate for the remainder of the mission planning horizon. Further, the significant offset between the two lines originates as a consequence of the consideration of large debris in the simulation. The latter ones are engaged multiple times, and their periapsis radius is decreased; however, these engagements are unable to deorbit them. Figure~\ref{fig:3d_all} illustrates the three constellation configurations obtained, the mixed debris field, and the valuable assets. The platforms' orbital slots for the three configurations are presented in Table~\ref{table:orbelements_all}.

\begin{landscape}
    \begin{table}[htb]
\caption{Orbital elements of platforms, defined at epoch $t_0$, for the mixed debris field case study.}
\centering
\begin{tabular}{lrrrrr}
\hline
Constellation & Sat. index & SMA, km & Incl., deg. & RAAN, deg & Arg. of latitude, deg. \\
\hline
Single platform & $p_1$ & 6953.14 & 55.625 & 280 & 280 \\
10 platform & $p_1$ & 6778.14 & 76.25 & 280 & 40 \\ 
& $p_2$ & 6778.14 & 76.25 & 280 & 80 \\ 
 & $p_3$ & 6865.64 & 62.50 & 120 & 280 \\ 
& $p_4$ & 6865.64 & 76.25 & 280 & 160 \\ 
&$p_5$ & 6865.64 & 90 & 280 & 160 \\ 
&$p_6$ & 6953.14 & 41.875 & 120 & 40 \\ 
&$p_7$ & 6953.14 & 55.625 & 280 & 280 \\ 
&$p_8$ & 7040.64 & 62.50 & 320 & 0 \\ 
&$p_9$ & 7128.14 & 83.125 & 280 & 120 \\ 
&$p_{10}$ & 7390.64 & 83.125 & 280 & 40 \\ 
10 platform Walker-Delta&$p_1$ & 6953.14 & 76.25 & 0 & 0 \\ 
&$p_2$ & 6953.14 & 76.25 & 0 & 180 \\ 
&$p_3$ & 6953.14 & 76.25 & 72 & 72 \\ 
&$p_4$ & 6953.14 & 76.25 & 72 & 252 \\ 
&$p_5$ & 6953.14 & 76.25 & 144 & 144 \\ 
&$p_6$ & 6953.14 & 76.25 & 144 & 324 \\ 
&$p_7$ & 6953.14 & 76.25 & 216 & 216 \\ 
&$p_8$ & 6953.14 & 76.25 & 216 & 36 \\ 
&$p_9$ & 6953.14 & 76.25 & 288 & 288 \\ 
&$p_{10}$ & 6953.14 & 76.25 & 288 & 108 \\  
\hline
\label{table:orbelements_all}
\end{tabular}
\end{table}
\end{landscape}

\subsection{General discussions}
Even though the case studies presented in this paper do not span all possible mission scenarios, we handpicked the most representative ones for debris remediation missions. Each mission scenario is uniquely characterized by a set of variables that condition the outcome of the space-based laser debris remediation mission. For the three different scenarios, we tested a triplet of constellations, two of them optimally locate their one and 10 platforms exploiting the MCLP formulation considering the debris field distribution, and the remaining one is an optimized Walker-Delta constellation from a pool of 360.

Increasing the number of platforms from one to 10 leads to a significant rise in the debris remediation capacity of the constellation and its derived metrics. In all case studies, the debris remediation capacity improves close to \SI{75}{\%} when the number of platforms increases. The results obtained in all three case studies enable us to conclude that further increments in the debris remediation capacity and derived metrics can be achieved using the same number of platforms, but breaking the symmetry of the constellation and determining the optimal location of platforms leveraging MCLP to deliver the additional increments. The MCLP provides flexibility in the location of laser platforms as it accounts for the mission environment defined and seeks to select the orbital slots, which are not constrained to follow a symmetrical pattern, that maximizes the constellation configuration reward. The outcome of this formulation significantly improves the obtained debris remediation capacity given by the L2D-ESP, as the platforms are better located relative to each other, the debris field considered, and, if any, the satellite-debris conjunction events. However, the lack of modeling debris orbit changes in the MCLP yields a larger number of L2D engagements and debris deorbited for the first time steps, and then their rate starts to plateau when the debris orbits start to diverge from the considered ones in the constellation configuration design problem. This phenomenon is clearly illustrated in Figs.~\ref{fig:periapsis_small}, \ref{fig:periapsis_large}, and \ref{fig:periapsis_all}.

\section{Sensitivity analysis} \label{sec:sensitivity}
This section performs a sensitivity analysis considering the same mission environment adopted in Sec.~\ref{sec:case_study_alldebris} for different optimal constellations while varying the number of platforms $P$ from one to 10. Additionally, a 10-platform Walker-Delta constellation is included in the study.

\begin{figure}[htbp]
    \centering
    \includegraphics[width=\linewidth]{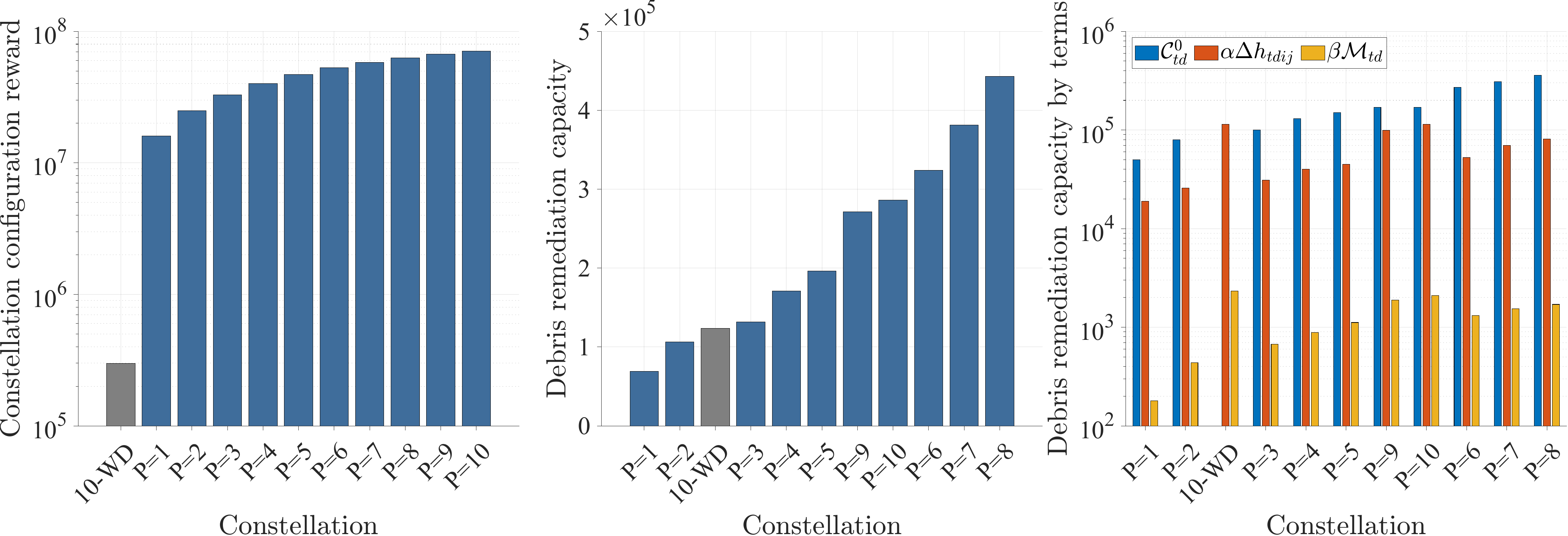}
    \caption{Obtained optimization rewards per constellation.}
    \label{fig:cost_benefit_rewards}
\end{figure}

Figure~\ref{fig:cost_benefit_rewards} highlights the rewards obtained by each constellation in the optimization problems, where 10-WD stands for the 10-platform Walker-Delta constellation. The leftmost bar plot of Fig.~\ref{fig:cost_benefit_rewards} presents the constellation configuration reward $\pi^{\ast}$, demonstrating a positive correlation between the increment in the number of platforms and the achieved reward. Further, all MCLP-based optimal constellations obtain a higher reward with respect to the 10-platform Walker-Delta constellation. The center of Fig.~\ref{fig:cost_benefit_rewards} outlines the total debris remediation capacity reward $V^\ast$ obtained in the L2D-ESP. Concerning this metric, the Walker-Delta constellation is outperformed by MCLP-based constellations with at least three platforms. Focusing on MCLP-based constellations only, the obtained $V^\ast$ does not maintain a monotonic increase with the increment in the number of platforms, since adopting 8 platforms retrieves the highest reward. The rightmost bar plot of Fig.~\ref{fig:cost_benefit_rewards} presents each term of the total debris remediation capacity reward defined in Eq.~\eqref{eq:reward}; the ${C}_{tdij}$ term has been omitted since no constellation performs L2D engagements that would trigger the penalty. First, ${C}^0_{td}$, the reward term that accounts for conjunctions if no L2D engagements occur, is obtained by all MCLP-based constellations as a consequence of the flexibility gained in the design space,  allowing the constellation to engage defunct satellite 22236 in the specified time window; conversely, the Walker-Delta constellation is not able to collect any ${C}^0_{td}$ reward. Further, the latter outperforms MCLP-based constellations analyzing individual reward terms $\alpha\Delta h_{tdij}$, which accounts for the ratio between debris periapsis radius after the L2D engagement and periapsis radius threshold, and $\beta{M}_{td}$, the term that accounts for engaged debris mass. However, it is not conveyed in the total reward given its lack to obtain ${C}^0_{td}$, the term that weighs the most in Eq.~\eqref{eq:reward}.

\begin{figure}[htbp]
    \centering
    \includegraphics[width=0.85\linewidth]{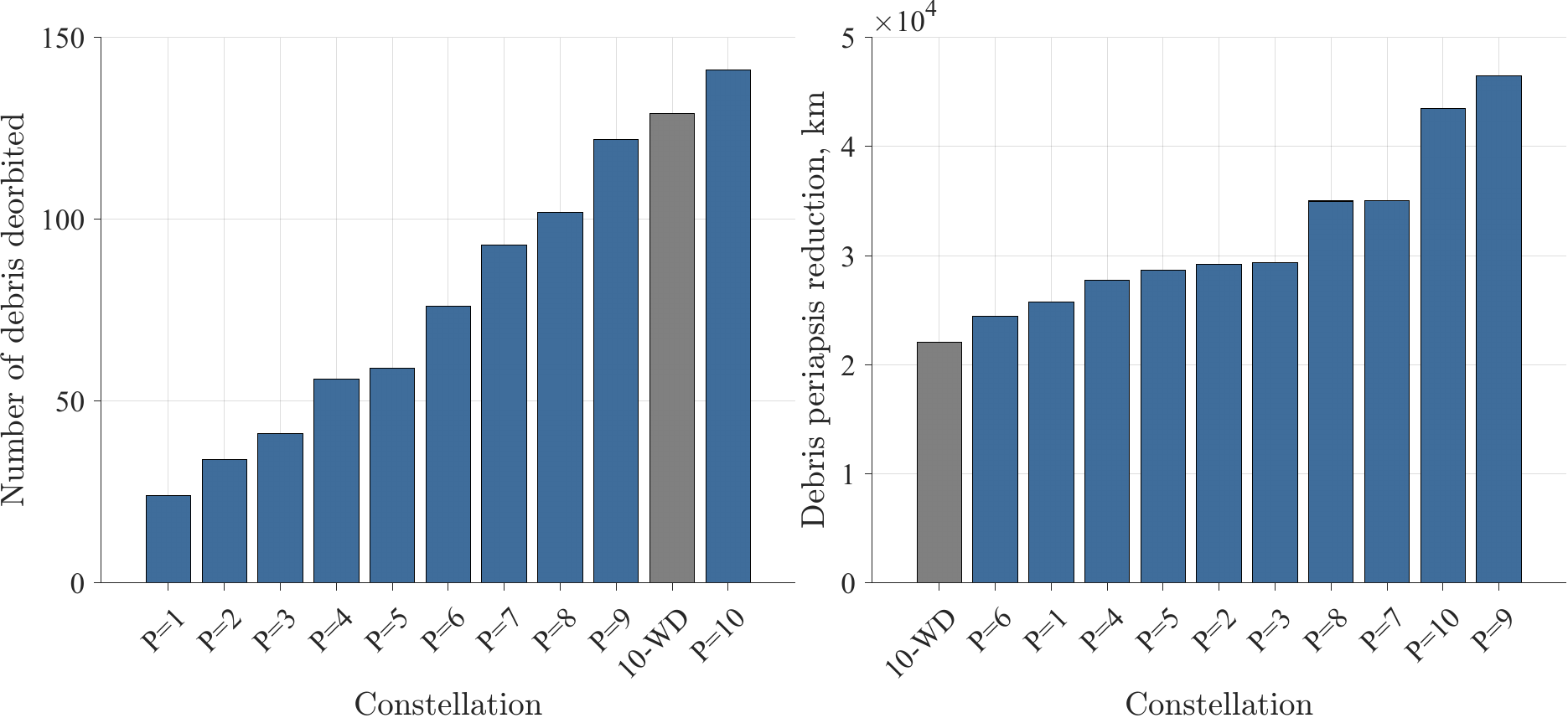}
    \caption{Optimization derived metrics per constellation.}
    \label{fig:cost_benefit_metrics}
\end{figure}

Figure~\ref{fig:cost_benefit_metrics} presents relevant mission metrics derived from the L2D-ESP. The left bar plot of Fig.~\ref{fig:cost_benefit_metrics} presents the number of small debris deorbited per constellation, exhibiting a positive correlation between the number of deorbited objects and the number of platforms. The right bar plot of Fig.~\ref{fig:cost_benefit_metrics} reports the magnitude of nudged debris per constellation, where the Walker-Delta constellation is outperformed by all MCLP-based constellations, and the 9-platform constellation provides the highest debris nudging value of \SI{46484.58}{km}. 

In summary, this section outlines the performance of different mission architectures varying the number of space-based lasers used and their constellation configuration. The results demonstrate the necessity of breaking the symmetry of Walker-Delta constellations. MCLP-based optimal constellations with at least three platforms outperform these symmetrical constellations in terms of debris remediation capacity. Adopting MCLP expands the optimization design space, allowing for the retrieval of probably optimal constellation configurations that do not adhere to the symmetrical Walker-Delta pattern.

\section{Limitations} \label{sec:limitations}

From the results obtained throughout this paper, certain key points need to be mentioned. First, this paper does not seek to endorse any specific laser platform system through the case studies presented. Instead, the case studies demonstrate the extension and capabilities of the proposed optimization formulations. Additionally, the optimization formulations can be used as a testbed to determine the optimal location and scheduling for any specific laser system. Second, increasing the number of platforms during the L2D-ESP does not guarantee a solution that is at least as good as one obtained with fewer platforms. This is attributed to the use of a myopic policy algorithm, which cannot obtain certified globally optimal solutions. Furthermore, due to the sequential nature of the algorithm, all missions included in the sensitivity analysis share the same debris field at the initial epoch. However, after the first L2D engagement, the debris field is likely to differ across the cases. Third, the proposed reward function is designed to capture what we consider to be the most critical figures of merit in a debris remediation mission. However, it is important to note that this representation is not exhaustive, nor is it the sole possible formulation of these figures of merit for such missions.

\section{Conclusions} \label{sec:conclusions}
This paper presents an optimization framework to determine the configuration and L2D scheduling of a constellation of cooperative space-based lasers, capable of deorbiting and nudging debris, as well as performing just-in-time collision avoidance based on any user-defined laser platform specifications. We propose the CLSP formulation to concurrently locate and schedule a constellation of space-based lasers; however, due to the rapid expansion in the solution space, we decouple the location and scheduling problems.
The adoption of the MCLP, for the first time, to design the constellation configuration of the debris remediation constellation enables us to consider critical variables of the mission environment, such as the characteristics of the targeted debris field, the number and specifications of the laser platforms, and the operation of valuable assets in space. The L2D-ESP maximizes the debris remediation capacity of the constellation during the operation of the mission. It leverages DVA, a collaborative engagement framework that enables achieving a higher degree of control over the L2D engagements.

This paper presents two sets of experiments. The first set comprehends case studies that exhibit the applicability of the optimization formulations proposed against two baseline cases, \textit{i.e.,} a single laser platform and a 10-laser Walker-Delta constellation. The case studies enable us to conclude that the debris remediation capacity can be increased if (1) the number of platforms used in the constellation is increased and (2) given the same number of laser platforms, breaking the symmetry of Walker-Delta constellations by leveraging the MCLP formulation. The second set is a sensitivity analysis that ranges the number of platforms from one to 10 and includes a 10-platform Walker-Delta constellation. Our findings highlight the benefits of adopting optimization into the placement and scheduling of space-based laser constellations to achieve scalable and effective solutions to remediate debris of different sizes. Further, the results obtained demonstrate that breaking the symmetry in the constellation configuration has a greater impact on the debris remediation capacity than increasing the number of platforms.

There are several avenues for future research. First, to include collision avoidance between laser platforms in the optimization. Second, explore various laser parameters, such as, but not limited to, input electrical power, PRF, pulse energy, and operational range, suitable for operation on small spacecraft. Consequently, a trade-off study, in the form of a cost-benefit analysis, can be conducted to determine the optimal location and scheduling of platforms by incorporating their size, mass, and laser parameters as decision variables. Third, perform a comparative analysis to determine the impact of the mission planning horizon's length on the constellation configuration design. To this extent, new solution methods (\textit{e.g.,} decomposition-based) will be explored. Lastly, to investigate different approaches to formulate the debris remediation capacity reward, and compare them to determine the most effective one.

\section*{Acknowledgment}
This work was supported by an Early Career Faculty grant from NASA’s Space Technology Research Grants Program under award No. 80NSSC23K1499. The authors would like to thank Gavin Baker (West Virginia University), Hao Chen (Stevens Institute of Technology), and Koki Ho (Georgia Institute of Technology), and the anonymous reviewers for their valuable comments and suggestions, which have helped improve the clarity and quality of this paper.

\appendix
\section{Full concurrent location-scheduling platform location problem formulation}
The full CLSP formulation introduced in Sec.~\ref{sec:allinone} is given as:
\begin{subequations}
\begin{alignat}{2}
\max \quad &\sum_{t\in\mathcal{T}\setminus\{t_{T-1}\}}\sum_{d\in\mathcal{D}}\sum_{i\in\mathcal{J}_{td}}\sum_{j\in\mathcal{J}_{t+1,d}} R_{tdij}x_{tdij} \tag{\ref{local:obj}}\\
    \text{s.t.} \quad
    &W_{tsd}z_s \ge y_{tsd}, &\quad \forall t\in\mathcal{T}, \forall s \in \mathcal{S},\forall d\in\mathcal{D}\tag{\ref{local:location_coupling}}\\
    &\sum_{d\in\mathcal{D}} y_{tsd} \le 1, &\quad \forall t\in\mathcal{T}, \forall s \in \mathcal{S}  \tag{\ref{local:eng_limit}}\\
    &\sum_{j\in\mathcal{J}_{t_1,d}} x_{t_{0},d,i_{0},j} = 1, &\quad \forall d \in \mathcal{D}\tag{\ref{local:flow_relocate}}\\
    &\sum_{j\in\mathcal{J}_{t+1,d}} x_{tdij} -  \sum_{\upsilon\in\mathcal{J}_{t-1,d}} x_{t-1,d\upsilon i} = 0, &\quad \forall t\in\mathcal{T}\setminus\{t_0,t_{T-1}\}, \forall d \in \mathcal{D}, \forall i \in \mathcal{J}_{td}\tag{\ref{local:flow_balance}}\\
    &\sum_{s\in\mathcal{S}_{tdj}} y_{tsd} \ge S_{tdj}x_{tdij}, &\quad  \forall t\in\mathcal{T}\setminus\{t_{T-1}\},\forall d\in\mathcal{D}, \forall i \in \mathcal{J}_{td},\forall j \in \mathcal{J}_{t+1,d}\setminus{\Tilde{\mathcal{J}}_{t+1,d}}\tag{\ref{local:engagements}}\\
    &\sum_{s\in\mathcal{S}} z_{s} = P \tag{\ref{local:cardinality}}\\
    & z_s \in \{0,1\}, &\quad\forall s\in \mathcal{S} \tag{\ref{local:x_s}} \\
    & y_{tsd} \in \{0,1\}, &\quad\forall t\in \mathcal{T},\forall s\in \mathcal{S}, \forall d \in \mathcal{D} \tag{\ref{local:y_sd}}\\
    & x_{tdij} \in \{0,1\}, &\quad\forall t\in \mathcal{T}\setminus\{t_{T-1}\},\forall d \in \mathcal{D},\forall i \in \mathcal{J}_{td},\forall j \in \mathcal{J}_{t+1,d} \tag{\ref{local:x_dkj}}
\end{alignat}
\end{subequations}

\section{Constellation configuration optimization formulation derivation}\label{appendix:MCLP}
In this section, we present a step-by-step justification for the adoption of the MCLP to tackle the optimal design of the space-based laser constellation configuration. As described in Sec.~\ref{sec:allinone}, the tree structure in the problem's solution space makes it computationally prohibitive. To overcome this problem, we adopt a series of assumptions that impact decision variables, parameters, and sets to obtain a space-based platform location optimization formulation.

The core assumption to overcome the tree structure is that debris maintains its initial orbit defined at $t_0$ throughout the entire mission. The consequences of this assumption are three-fold; first, debris cannot be deorbited, which makes the debris field constant over time. Second, debris orbital slot indices $i,j$ can be dropped; consequently, constraints~\eqref{local:flow_relocate} and~\eqref{local:flow_balance} are relaxed. Third, given that orbital slot indices are dropped, the CLSP decision variables $x_{tdij}$ are now defined as:
\begin{equation}
    x_{td} = \begin{cases}
    1, & \text{if debris $d$ is engaged at time step $t$} \\
    0, & \text{otherwise}
\end{cases}
\end{equation}

As a result of dropping the relocation indices, reward $R_{tdij}$ defined in Eq.~\eqref{eq:reward} is recast as $R_{td}={C}^0_{td} + {M}_{td}$. Removing $\Delta h_{tdij}$ prevents the optimization formulation from being informed about the post-L2D engagement behavior of debris (\textit{i.e.} whether the debris periapsis radius increases or decreases). To overcome this situation, we redefine $W_{tsd}$ as $W_{tsd}'$, where the latter incorporates an additional requirement that the debris periapsis radius after an L2D engagement must be lower than its initial value.

With the described assumptions and relaxations, the derived formulation is:    
\begin{subequations}
\begin{alignat}{2}
\max \quad &\sum_{t\in\mathcal{T}}\sum_{d\in\mathcal{D}}R_{td}x_{td} \tag{\ref{loc:obj}}\\
    \text{s.t.} \quad
    &W'_{tsd}z_s \ge y_{tsd}, &\quad \forall t\in\mathcal{T}, \forall s \in \mathcal{S},\forall d\in\mathcal{D}\label{eq:deriv1_coupling_z_y}\\
    &\sum_{d\in\mathcal{D}} y_{tsd} \le 1, &\quad \forall t\in\mathcal{T}, \forall s \in \mathcal{S}\tag{\ref{local:eng_limit}}\\
    &\sum_{s\in\mathcal{S}_{td}} y_{tsd} \ge S_{td}x_{td}, &\quad  \forall t\in\mathcal{T},\forall d\in\mathcal{D}\label{eq:deriv1_coupling_z_x}\\
    &\sum_{s\in\mathcal{S}} z_{s} = P,\tag{\ref{local:cardinality}}\\
    & z_s \in \{0,1\}, & \forall s\in \mathcal{S} \tag{\ref{local:x_s}} \\
    & y_{tsd} \in \{0,1\}, & \quad\forall t\in \mathcal{T},\forall s\in \mathcal{S}, \forall d \in \mathcal{D} \tag{\ref{local:y_sd}}\\
    & x_{td} \in \{0,1\}, & \quad\forall t\in \mathcal{T},\forall d \in \mathcal{D}\tag{\ref{loc:y_td}}
\end{alignat}
\end{subequations}

Even though relocation constraints are relaxed, the structure of the problem inherits the location-scheduling structure of the CLSP, since $y_{tsd}$ are the L2D engagement scheduling decision variables that dictate whether an L2D engagement occurs from a platform located at slot $s$ to debris $d$ at time step $t$. Consequently, despite the fact that the new formulation does not have a tree structure, the dimension of decision variables $y_{tsd} \in \{0,1\}^{T \times S \times D}$ can make the problem computationally prohibitive.
On top of that, given that the reward does not account for debris orbit changes, the formulation will tend to assign to each platform those orbital slots that maximize engagement with debris of larger mass.
Given two debris $d_1$ and $d_2$ with masses $m_{d_1} < m_{d_2}$ and ${C}^0_{td_1}={C}^0_{td_2}=0$, then ${M}_{td_1} < {M}_{td_2}$ implies that $R_{td_1} < R_{td_1}$. However, throughout the L2D-ESP where changes in debris periapsis are considered, if $\Delta h_{t,d_1,i,j} \gg \Delta h_{t,d_2,i,j}$, then:
\begin{subequations}
    \begin{alignat}{2}
    & \alpha\Delta h_{td_1,i,j}  + \beta{M}_{td_1} > \alpha\Delta h_{t,d_1,i,j}  + \beta{M}_{td_1}\\
    & R_{t,d_1,i,j} > R_{t,d_2,i,j}
    \end{alignat}
\end{subequations}
outlining that the outcome, which corresponds to the optimal constellation that maximizes reward $R_{td}$, can overlook platform orbital slots that collect higher debris remediation rewards during the L2D-ESP. To tackle this problem, and considering that it is not proper to perform L2D assignments if debris change in orbit is neglected, we propose to relax CLSP L2D engagement scheduling decision variables $y_{tsd}$. Consequently, constraints~\eqref{local:location_coupling},~\eqref{local:eng_limit},~\eqref{local:engagements} and~\eqref{local:y_sd} are dropped. The resulting optimization problem is given as:
\begin{subequations}
\begin{alignat}{2}
\max \quad &\sum_{t\in\mathcal{T}}\sum_{d\in\mathcal{D}}R_{td}x_{td} \tag{\ref{loc:obj}}\\
    \text{s.t.} \quad
    &\sum_{s\in\mathcal{S}} z_{s} = P,\tag{\ref{local:cardinality}}\\
    & z_s \in \{0,1\}, & \forall s\in \mathcal{S} \tag{\ref{local:x_s}}\\
    & x_{td} \in \{0,1\}, & \quad\forall t\in \mathcal{T},\forall d \in \mathcal{D}\tag{\ref{loc:y_td}}
\end{alignat}
\end{subequations}

The new problem lacks coupling constraints that link location decision variables $z_s$ and decision variables $x_{td}$ such that L2D engagements with debris $d$ occur only if a platform is occupying an orbital slot such that $W_{tsd}=1$. Inspired by constraints~\eqref{eq:deriv1_coupling_z_y} and~\eqref{eq:deriv1_coupling_z_x}, we introduce the new coupling constraints given as:
\begin{equation}
    \sum_{s\in\mathcal{S}} W'_{tsd}z_s\ge S_{td} x_{td}, \quad \forall t\in \mathcal{T}, \forall d \in \mathcal{D}\tag{\ref{loc:coupling}}
\end{equation}
where $x_{td}$ are activated if at least $S_{td}$ platforms perform L2D ablation on debris $d$ at time step $t$. Hence, the formulation imposes a minimum number of platforms required to obtain debris' reward $R_{td}$ at time step $t$. Consequently, the optimal constellation configuration design problem is given as:
\begin{subequations}
\begin{alignat}{2}
    \max \quad &\sum_{t\in\mathcal{T}}\sum_{d\in \mathcal{D}}R_{td} x_{td}\tag{\ref{loc:obj}}\\
    \text{s.t.} \quad
    &\sum_{s\in\mathcal{S}} z_{s} = P\tag{\ref{local:cardinality}}\\
    &\sum_{s\in\mathcal{S}} W'_{tsd}z_s\ge S_{td} x_{td}, &\forall t\in \mathcal{T}, \forall d \in \mathcal{D}\tag{\ref{loc:coupling}}\\
    & z_s \in \{0,1\}, & \forall s\in \mathcal{S}\tag{\ref{loc:z_s}}\\
    & x_{td} \in \{0,1\}, & \quad\forall t\in \mathcal{T}, \forall d \in \mathcal{D}\tag{\ref{loc:y_td}}
\end{alignat}
\end{subequations}
where its structure resembles the well-known MCLP proposed by~\citep{church1974maximal}.

\section{Case study with range-dependent laser fluence}\label{appendix:new_case_study}
In this Appendix, we aim to outline the model's flexibility by presenting a case study where the pulse energy is constant and the delivered fluence on target is range-dependent. Hence, in this Appendix, the magnitude of the L2D $\Delta v$ is different from Secs.~\ref{sec:case_studies} and~\ref{sec:sensitivity}. We obtain an optimal 10-platform constellation configuration, and then compare its constellation configuration reward and debris remediation capacity against an optimally located single laser platform and an optimized 10-platform Walker-Delta constellation.

The case study considers 395 small debris randomly generated as described in Sec.~\ref{sec:case_study_small}. The epoch is defined as February 26, 2024, at 04:30:51 UTC, and the mission time horizon is seven days uniformly discretized in time steps of \SI{160}{s} size. The laser platform orbital slots are generated from a list of 9 altitude and inclination steps uniformly distributed in the range between \SI{400}{km} and \SI{1100}{km}, and between \SI{35}{deg} and \SI{90}{deg}, respectively. The RAAN and argument of latitude are uniformly distributed in 9 steps spanning \SI{360}{deg}. The laser parameters are adopted from~\citep{pieters_asr_2023}, and the delivered fluence is calculated with Eq.~\eqref{eq:fluence}. We assume that debris has a surface density of \SI{0.2}{kg/m^2}. Table~\ref{table:appendix_parameters} presents the laser parameters used.

\begin{table}[htbp]
\caption{Case study laser parameters.}
\centering
\begin{tabular}{lll}
\hline
Parameter & Value &Unit\\
\hline
Pulse energy, $E$ & 300 &\si{J}\\
Effective beam diameter, $D_{\text{eff}}$ &2 &\si{m}\\
Sys. loss factor, $T_\text{tot}$ & 0.9 &-\\
Beam quality factor, $B^2$ & 2 &-\\
Diffraction constant $\zeta$ & 1.27 &-\\
Wavelength, $\lambda$   &335 &\si{nm}\\
Laser range, $u$ & $[30,200]$ & \si{km}\\
PRF & 66.6 & \si{Hz}\\
Length of interaction & 10 & \si{s}\\
Momentum coupling coefficient, $c_\text{m}$ &30& \si{N/MW}\\
\hline
\label{table:appendix_parameters}
\end{tabular}
\end{table}

\begin{figure}[htbp]
    \centering
     \begin{subfigure}[b]{0.48\textwidth}
        \centering
        \includegraphics[width=0.8\linewidth]{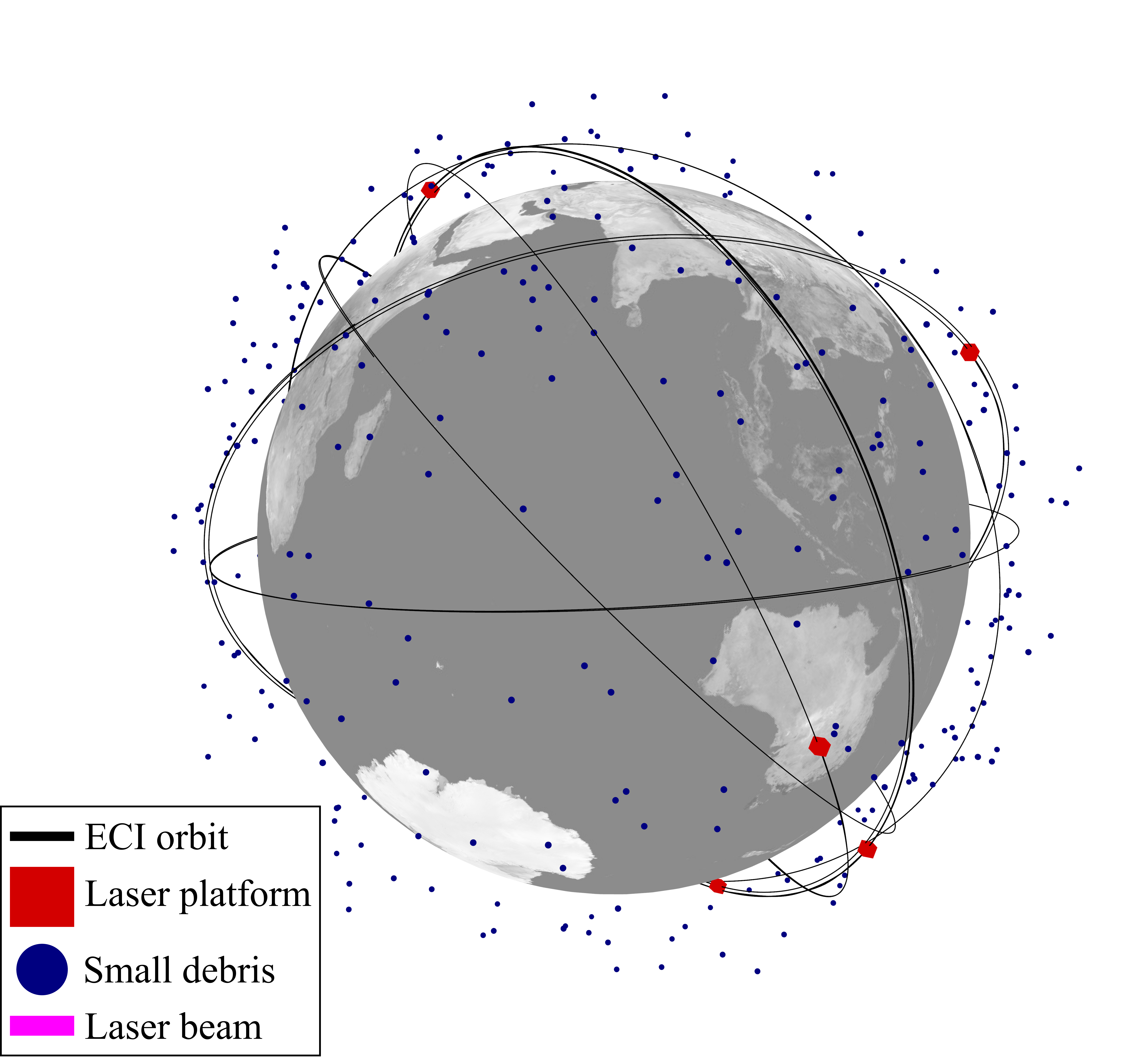}
        \caption{Small debris field and laser platforms at $t=t_0$.}
        \label{fig:app_small_field}
     \end{subfigure}
     \begin{subfigure}[b]{0.48\textwidth}
        \centering
        \includegraphics[width=0.8\linewidth]{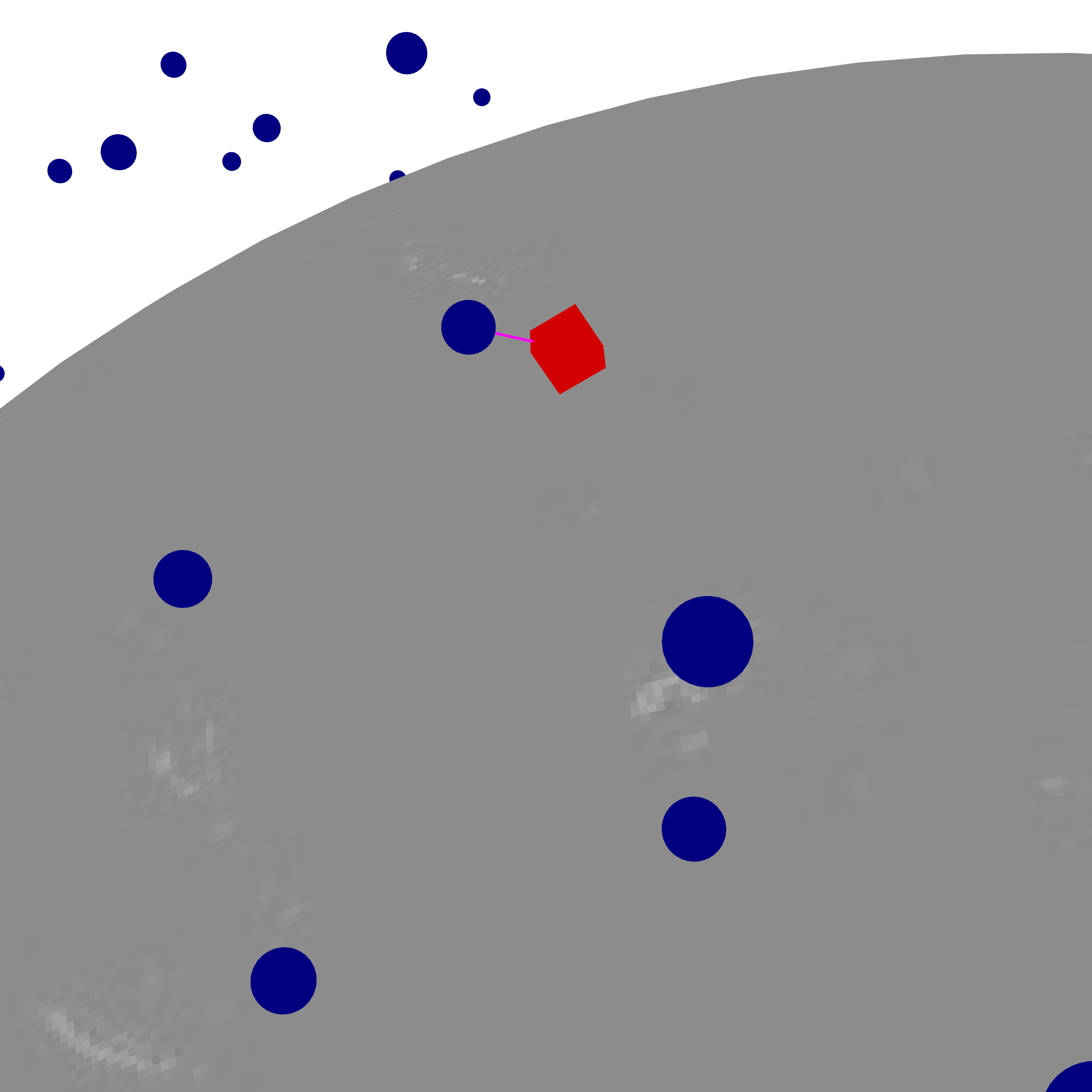}
        \caption{Detailed L2D engagement at $t=235$.}
        \label{fig:app_beam_small_field}
    \end{subfigure}
    \caption{Optimal 10-platform constellation for small debris remediation mission.}
\end{figure}

The optimal 10-platform constellation collects a constellation configuration reward of $\pi^\ast=1804$ and has a debris remediation capacity of $V^\ast=$ \num{22042.71}. Figure~\ref{fig:app_small_field} illustrates the constellation configuration at the epoch $t=t_0$ with the orbits in the ECI frame, the small debris field, and the laser platforms. Figure~\ref{fig:app_beam_small_field} provides a detailed view of a single laser platform engaging debris at the time step $t={235}$. The constellation configuration deorbits 156 debris and nudges \SI{16775.68}{km}.

The single laser platform obtains a constellation configuration reward of $\pi^\ast=259$ and has a debris remediation capacity of $V^\ast=$ \num{2642.21}. Representing a decrement of \SI{85.64}{\%} and \SI{88.01}{\%} of constellation configuration reward and debris remediation capacity, respectively. Further, the single platform deorbits 21 debris and nudges \SI{3005.15}{km}, outlining decrements of \SI{86.53}{\%} and \SI{82.08}{\%} for debris deorbited and nudged, respectively.

\begin{figure}[htbp]
    \centering
    \includegraphics[width=\linewidth]{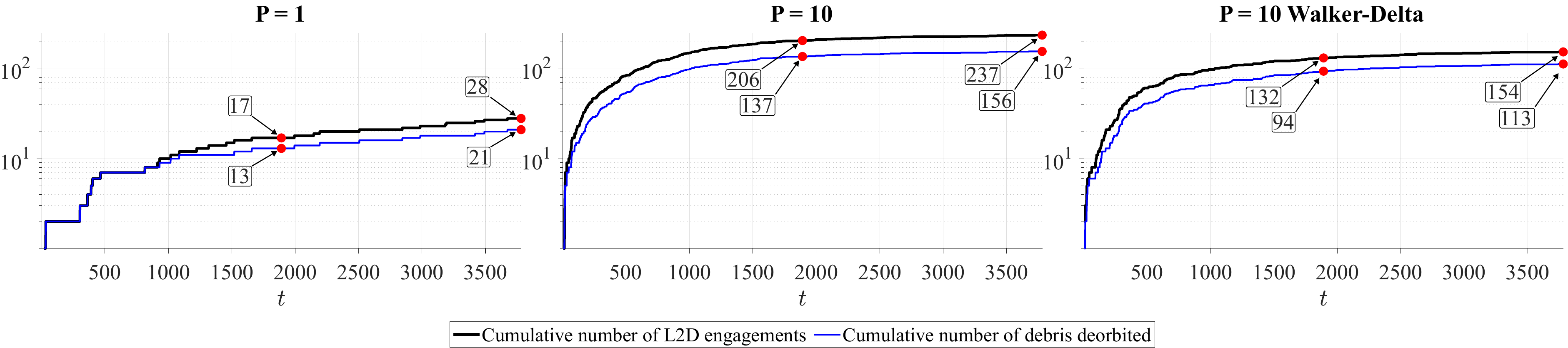}
    \caption{Cumulative number of L2D engagements and small debris deorbited. The highlighted points correspond to the cumulative number of L2D engagements and debris deorbited at the middle and end of the mission planning horizon.}
    \label{fig:periapsis_appendix}
\end{figure}

\begin{figure}[htbp]
    \centering
    \begin{subfigure}[b]{0.32\textwidth}
        \centering
        \includegraphics[width=\linewidth]{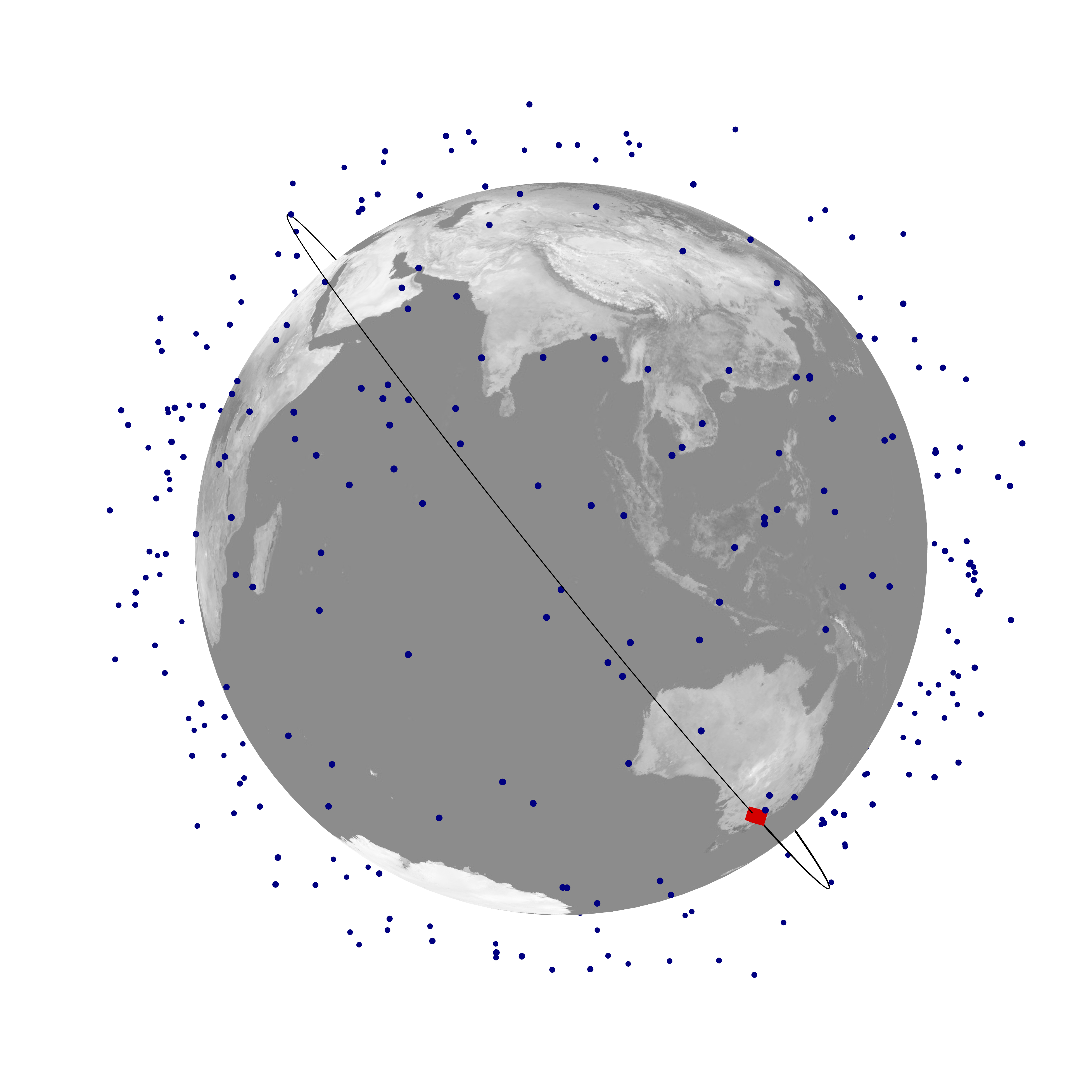}
        \caption{$P=1$.}
     \end{subfigure}
     \begin{subfigure}[b]{0.32\textwidth}
        \centering
        \includegraphics[width=\linewidth]{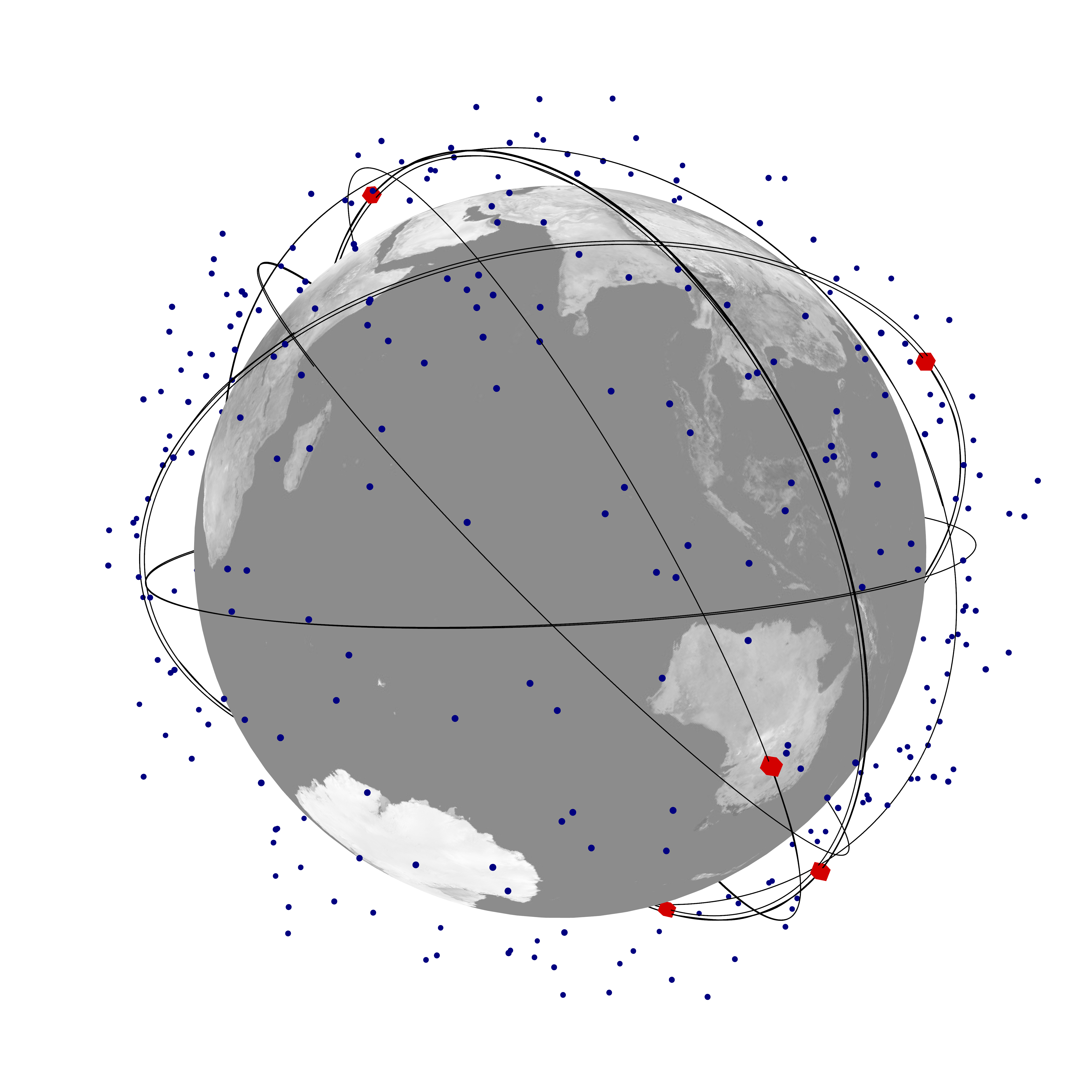}
        \caption{$P=10$.}
    \end{subfigure}
    \begin{subfigure}[b]{0.32\textwidth}
        \centering
        \includegraphics[width=\linewidth]{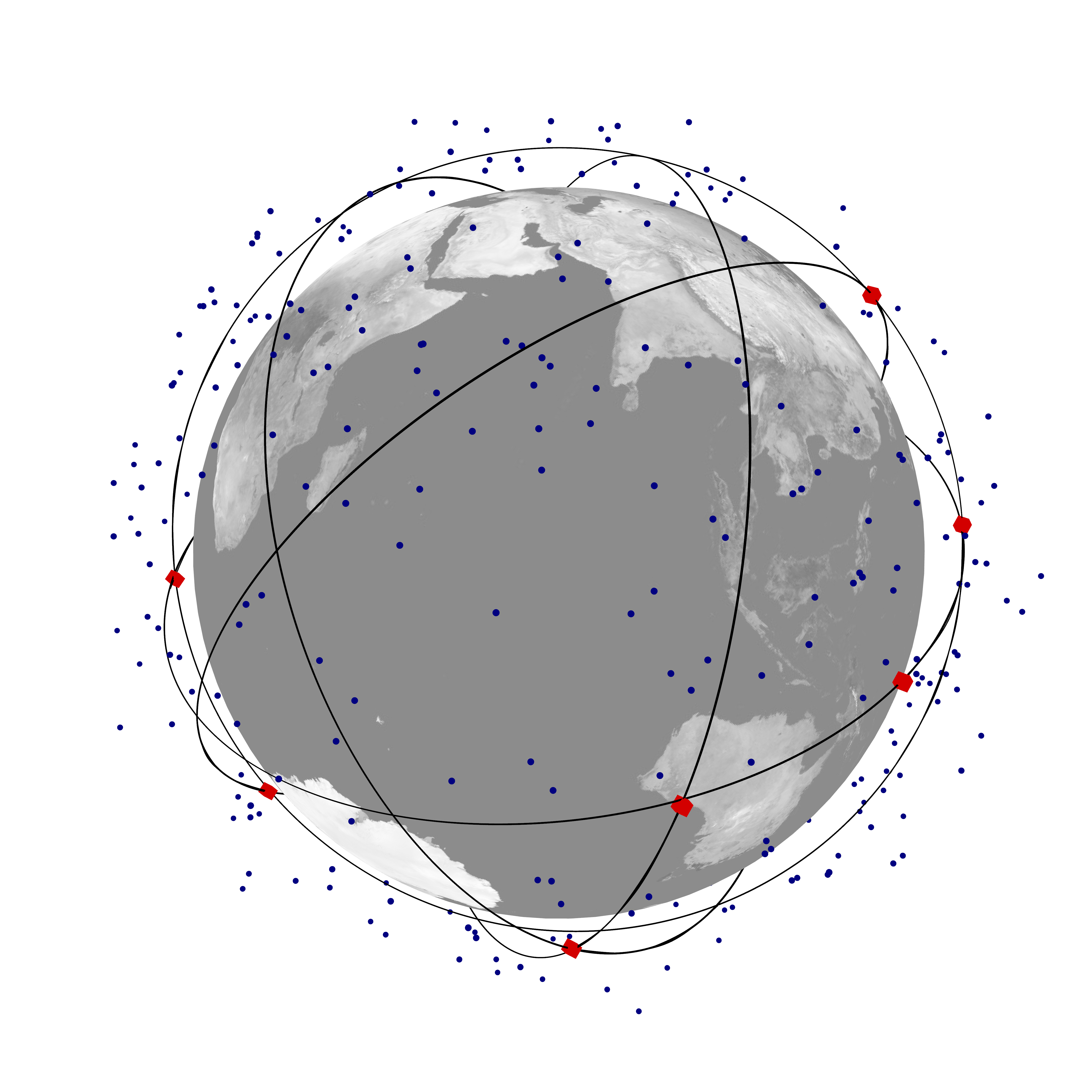}
        \caption{$P=10$ Walker-Delta.}
    \end{subfigure}
    \caption{Small debris field with range-dependent fluence case study: 3D visualization of debris remediation constellations with orbits in ECI at $t=t_0$.}
    \label{fig:app_3d_small}
\end{figure}

The best-performing Walker-Delta constellation collects a constellation configuration reward of $\pi^\ast=801$ and has a debris remediation capacity of $V^\ast=$ \num{14874.32}. The percentage decrements with respect to the optimal 10-platform constellation are \SI{55.59}{\%} and \SI{32.52}{\%} for the constellation configuration reward and debris remediation capacity, respectively. Furthermore, the optimized Walker-Delta constellation deorbits 113 small debris, \SI{27.56}{\%} less than the optimal 10-platform constellation, and nudges \SI{11642.80}{km}, \SI{30.59}{\%} less than the optimal 10-platform constellation. Figure~\ref{fig:periapsis_appendix} showcases the cumulative number of L2D engagements and small debris deorbited. The single-platform configuration successfully deorbits all debris engaged for the first time steps, and then the later L2D engagements can only decrease debris periapsis radius without deorbiting them. Similarly, the 10 platform constellation configurations have a similar trend but for a shorter number of time steps. Further, both lines plateau as the time step index increases while maintaining a similar offset between them. Lastly, Fig.~\ref{fig:app_3d_small} presents side-by-side the single platform and the two constellation configurations analyzed, and Table~\ref{table:app_orbelements_small} outlines their orbital elements.

\begin{landscape}
    \begin{table}[htb]
\caption{Orbital elements of platforms, defined at epoch $t_0$, for the small debris field with range-dependent fluence case study.}
\centering
\begin{tabular}{lrrrrr}
\hline
Constellation & Sat. index & SMA, km & Incl., deg. & RAAN, deg & Arg. of latitude, deg.\\
\hline
1 platform & $p_1$  & 7478.14 & 35     & 120 & 240  \\
10 platform & $p_1$  & 7040.64 & 35     & 160 & 280   \\
&$p_2$  & 7215.64& 35 & 80  &0\\
&$p_3$  & 7215.64& 35 &160&80\\
&$p_4$  & 7215.64 & 35 & 160 & 240 \\
&$p_5$  & 7215.64 &35& 160& 320\\
&$p_6$  & 7215.64& 35& 200&320\\
&$p_7$  & 7215.64&35&280&80\\
&$p_8$  & 7215.64&41.875& 0 & 160\\
&$p_9$  & 7303.14& 35 & 280 & 160\\
&$p_{10}$ & 7478.14&35 & 120 &240\\
10 platform Walker-Delta& $p_{1}$  & 7040.64 & 62.5 & 0   & 0\\
&$p_{2}$  & 7040.64 & 62.5 & 0   & 180 \\
&$p_{3}$  & 7040.64 & 62.5 & 72   & 108 \\
&$p_{4}$  & 7040.64 & 62.5 & 72  & 288 \\
&$p_{5}$  & 7040.64 & 62.5 & 144   & 216 \\
&$p_{6}$  & 7040.64 & 62.5 & 144   & 36 \\
&$p_{7}$  & 7040.64 & 62.5 & 216 & 324 \\
&$p_{8}$  & 7040.64 & 62.5 & 216 & 144 \\
&$p_{9}$  & 7040.64 & 62.5 & 288 & 72 \\
&$p_{10}$ & 7040.64 & 62.5 & 288 & 252 \\
\hline
\label{table:app_orbelements_small}
\end{tabular}
\end{table} 
\end{landscape}

\clearpage
\newpage
\bibliographystyle{jasr-model5-names}
\biboptions{authoryear}
\bibliography{references}
\end{document}